\documentclass[12pt,leqno]{article}
\usepackage{amsfonts}
\usepackage{amsmath, amsthm, amsfonts, amssymb, color}
\usepackage{mathrsfs}
\usepackage{color}
\usepackage[notcite,notref]{showkeys}

\theoremstyle{definition}

\newcommand{\scr}[1]{\mathscr #1}
\definecolor{wco}{rgb}{0.5,0.2,0.3}

\numberwithin{equation}{section} \theoremstyle{remark}

\newcommand{\ua}{\uparrow}

\title{{\bf  Convergence in Wasserstein Distance   for  Empirical Measures of   Non-Symmetric Subordinated Diffusion Processes}\footnote{Supported in
 part by   the National Key R\&D Program of China (No. 2022YFA1006000, 2020YFA0712900) and NNSFC (11921001).} }
\author{{\bf  Feng-Yu Wang   }\\
\footnotesize{ Center for Applied Mathematics, Tianjin University, Tianjin 300072, China}\\
  \footnotesize{  Department of Mathematics,
Swansea University, Bay Campus, SA1 8EN, United Kingdom}\\
\footnotesize{    wangfy@tju.edu.cn}}
\begin{document}
\allowdisplaybreaks
\def\R{\mathbb R}  \def\ff{\frac} \def\ss{\sqrt} \def\B{\mathbf
B}
\def\N{\mathbb N} \def\kk{\kappa} \def\m{{\bf m}}
\def\ee{\varepsilon}\def\ddd{D^*}
\def\dd{\delta} \def\DD{\Delta} \def\vv{\varepsilon} \def\rr{\rho}
\def\<{\langle} \def\>{\rangle}
  \def\nn{\nabla} \def\pp{\partial} \def\E{\mathbb E}
\def\d{\text{\rm{d}}} \def\bb{\beta} \def\aa{\alpha} \def\D{\scr D}
  \def\si{\sigma} \def\ess{\text{\rm{ess}}}\def\s{{\bf s}}
\def\beg{\begin} \def\beq{\begin{equation}}  \def\F{\scr F}
\def\Ric{\mathcal Ric} \def\Hess{\text{\rm{Hess}}}
\def\e{\text{\rm{e}}} \def\ua{\underline a} \def\OO{\Omega}  \def\oo{\omega}
 \def\tt{\tilde}\def\[{\lfloor} \def\]{\rfloor}
\def\cut{\text{\rm{cut}}} \def\P{\mathbb P} \def\ifn{I_n(f^{\bigotimes n})}
\def\C{\scr C}      \def\aaa{\mathbf{r}}     \def\r{r}
\def\gap{\text{\rm{gap}}} \def\prr{\pi_{{\bf m},\nu}}  \def\r{\mathbf r}
\def\Z{\mathbb Z} \def\vrr{\nu} \def\ll{\lambda}
\def\L{\scr L}\def\Tt{\tt} \def\TT{\tt}\def\II{\mathbb I}
\def\i{{\rm in}}\def\Sect{{\rm Sect}}  \def\H{\mathbb H}
\def\M{\mathbb M}\def\Q{\mathbb Q} \def\texto{\text{o}} \def\LL{\Lambda}
\def\Rank{{\rm Rank}} \def\B{\scr B} \def\i{{\rm i}} \def\HR{\hat{\R}^d}
\def\to{\rightarrow} \def\gg{\gamma}
\def\EE{\scr E} \def\W{\mathbb W}
\def\A{\scr A} \def\Lip{{\rm Lip}}\def\S{\mathbb S}
\def\BB{\scr B}\def\Ent{{\rm Ent}} \def\i{{\rm i}}\def\itparallel{{\it\parallel}}
\def\g{{\mathbf g}}\def\Sect{{\mathcal Sec}}\def\T{\mathcal T}\def\BB{{\bf B}}
\def\f\ell \def\g{\mathbf g}\def\BL{{\bf L}}  \def\BG{{\mathbb G}}
\def\Bd{{D^E}} \def\BdP{D^E_\phi} \def\Bdd{{\bf \dd}} \def\Bs{{\bf s}} \def\GA{\scr A}
\def\Bg{{\bf g}}  \def\Bdd{\psi_B} \def\supp{{\rm supp}}\def\div{{\rm div}}
\def\ddiv{{\rm div}}\def\osc{{\bf osc}}\def\1{{\bf 1}}\def\BD{\mathbb D}
\def\H{{\bf H}}\def\gg{\gamma} \def\n{{\mathbf n}}\def\GG{\Gamma}\def\HAT{\hat}
\def\SU{{\bf SU}}
\maketitle

\begin{abstract} By using the spectrum of the underlying symmetric diffusion operator,
the convergence in $L^p$-Wasserstein distance $\W_p (p\ge 1)$ is characterized  for the empirical measure $\mu_t$ of non-symmetric subordinated diffusion processes in an abstract framework. The main results are applied to  the subordinations of several typical models, which include the (reflecting) diffusion processes on compact manifolds, the conditional diffusion processes, the Wright-Fisher diffusion process, and hypoelliptic diffusion processes
on {\bf SU}(2).
In particular, for the (reflecting) diffusion processes on a compact Riemannian manifold with invariant probability measure $\mu$: \beg{enumerate}
\item[(1)] the sharp limit of $t\W_2(\mu_t,\mu)^2$ is derived in $L^q(\P)$ for concrete $q\ge 1,$
which provides a precise characterization on  the physical observation that a divergence-free  perturbation accelerates the convergence in $\W_2$;
\item[(2)] the sharp convergence rates are presented for   $(\E[\W_{2p}(\mu_t,\mu)^{q}])^{\ff 1 q} (p,q\ge 1)$, where   a critical phenomenon appears  with the critical rate $t^{-1}\log t$ as $t\to\infty$.
\end{enumerate} 

  \end{abstract} \noindent

 \vskip 2cm

 \section{Introduction}

{\bf A. Background of the study.} In statistical physics, the empirical measure is a fundamental object to simulate the stationary distribution (Gibbs measure). Since  the ground breaking series work  \cite{DV} (1975-1983) where Donsker and Varadhan developed  their   celebrated  larger deviation principle,   the   long time behavior of empirical measures has become a key research topic in the study of Markov processes,   see   \cite{Wu,Wu2} for criteria on the central limit theorem and large deviations for hyperbounded Markov processes.

On the other hand, the Wasserstein distance is  intrinsic   in  the theory of optimal transport and calculus on Wasserstein space, see   \cite{AM,V} and references therein.   So, it is crucial and interesting to study the convergence in Wasserstein distance for the empirical measure  of Markov processes.

Moreover, it was observed in \cite{GG} that a divergence-free perturbation to  symmetric stochastic systems
may accelerate the algorithm of Gibbs measures. This has been confirmed in several papers for the convergence of Markov semigroups to stationary distributions, see \cite{HM0,HM, HM2} and references therein. It is interesting to provide a sharp characterization on the acceleration for the convergence  of empirical measures in Wasserstein distance.

In   recent years, the sharp convergence rate in the second moment of the $L^2$-Wasserstein distance has been derived in \cite{W1, W2, W3, W4, WZ} for empirical measures of symmetric diffusion processes. In particular, in lower dimensions the precise limit is explicitly formulated by using eigenvalues and eigenfunctions of the generator. These results have been extended to subordinated processes in \cite{WW,LW1,LW2,LW3} and the fractional Brownian motion on torus in \cite{HT}, see also \cite{DU} for the study of McKean-Vlasov SDEs.

\ \newline
{\bf B. Purpose of the present work.}
Based on the above background,   this paper   investigates the convergence of empirical measures for \emph{non-symmetric subordinated} diffusion processes in an \emph{abstract framework}, describes the \emph{acceleration of   the convergence} for  divergence-free perturbations to symmetric systems, and illustrates the main results
by   typical examples.

To  figure out a clear picture of  our  general results (see Section 2 for details), in the following we only consider non-symmetric diffusion processes on a compact manifold. See Section 5 for applications of the general results to three more examples including the subordinated conditional diffusion process,
the subordinated Wright-Fisher process, and the subordinated subelliptic diffusion process on $\SU(2)$.

\ \newline
{\bf C. A picture for non-symmetric diffusion processes on compact manifolds.}
Let $M$ be an $n$-dimensional compact connected Riemannian manifold possibly with a boundary $\pp M$. Let $\scr P$ be the space of all probability measures on $M$, let $\rr$ be the Riemannian distance, and for any $p\ge 1,$ let $\W_p$ be the $L^p$-Wasserstein distance induced by $\rr$, cf. \eqref{WP} below.

Let $\mu(\d x):=\e^{V(x)}\d x\in \scr P,$   where $V\in C^2(M)$  and $\d x$ is the volume measure on $M$, and let $Z$ be a  $C^1$-vector field with   ${\rm div}_\mu Z=0$, i.e.
$$\mu(Zf):= \int_M \<Z,\nn f\>\d\mu=0,\ \ f\in C^1(M).$$   Then the spectrum of  $\hat L
:=\DD+\nn V$ (with Neumann boundary if $\pp M$ exists) is discrete, and all eigenvalues $\{\ll_i\}_{i\ge 0} $  of $-\hat L$
listed in the increasing order counting multiplicities satisfy
$$c_1 i^{\ff n 2}\le \ll_i\le c_2 i^{\ff n 2},\ \ \ i\ge 0$$
for some constants $c_1,c_2>0,$ see for instance \cite{Chavel}. Let $\{\phi_i\}_{i\ge 0}$ with $\phi_0\equiv 1$ being the corresponding unitary eigenfunctions in $L^2(\mu)$.

 Let $X_t$   be the   diffusion process  on $M$ generated by $$L:= \DD+\nn V +Z, $$     with reflecting boundary if $ \pp M$ exists.
We consider the empirical measure
$$\mu_t:= \ff 1 t \int_0^t \dd_{X_s}\d s,\  \ \ t>0,$$ where $\dd_{X_s}$ is the Dirac measure at $X_s$.
By the central limit theorem (see  \cite{Wu}), for any
$$f\in L_0^2(\mu):=\big\{f\in L^2(\mu):\ \mu(f)=0\big\},$$ we have
\beq\label{CLT}  \lim_{t\to \infty}\ss t \mu_t(f) = \lim_{t\to \infty} \ff 1 {\ss t}\int_0^t f(X_s)\d\s =N(0,{\bf V}(f))\ \ \text{in\ law},\end{equation}
where $N(0,{\bf V}(f))$ is the centered normal distribution on $\R$ with variance
\beq\label{L-3}{\bf V}(f):= \int_0^\infty \mu(fP_sf)\d s= \mu(|\nn L^{-1} f|^2).\end{equation}
For any $k,R\ge 1$, let
\beq\label{KR}\scr P_{k,R}:= \big\{\nu\in \scr P:\ \d\nu=h\d\mu, \|h\|_k\le R\big\}, \end{equation}
where $\|\cdot\|_k$ is the norm in $L^k(\mu).$
For any $\nu\in \scr P$, let $\E^\nu$ be the expectation for the diffusion process $X_t$ with initial distribution $\nu$.

We first consider the long time behavior of $\W_2(\mu_t,\mu)$.
The following  result shows that when  $n\le 3$ and $\pp M$ is either empty or convex, for long time
 $t\W_2(\mu_t,\mu)^2$ behaves as
$$\Xi(t):= \sum_{i=1}^\infty \ff 1 {\ll_i} |\psi_i(t)|^2,\ \ \ \psi_i(t):= \ff 1 {\ss t}\int_0^t \phi_i(X_s)\d s,$$ so that  uniformly in $\nu\in \scr P$,
$t\E^\nu  [\W_2(\mu_t,\mu)^2]$ converges to
\beq\label{ZZ}  \eta_Z:= \sum_{i=1}^\infty \ff{2{\bf V}(\phi_i)}{\ll_i} = \sum_{i=1}^\infty \ff{2}{\ll_i^2}\Big(1-\ff 1 {\ll_i} {\bf V}(Z\phi_i)\Big), \end{equation}
where the second equality follows from Lemma \ref{PH} below.

\beg{thm}\label{TN}   There exists a constant $\kk  \ge 1$ with $\kk=1$ when $\pp M$ is either empty or convex, such that the following assertions hold.
\beg{enumerate}
\item[$(1)$]
When $n\le 2,$  for any $q\in [1,\ff{2n}{(3n-4)^+})$,
\beq\label{UL} \lim_{t\to\infty} \sup_{\nu\in \scr P} \E^\nu\Big[ \big|\big\{t\W_2(\mu_t,\mu)^2- \Xi(t)\big\}^+
 + \big\{\Xi(t)- \kk t\W_2(\mu_t,\mu)^2\big\}^+\big|^q\Big] =0, \end{equation}
 so that when $\pp M$ is either empty or convex,
 \beq\label{UL'} \lim_{t\to\infty} \sup_{\nu\in \scr P} \E^\nu \Big[\big|t\W_2(\mu_t,\mu)^2- \Xi(t)|^q \Big]
  =0.  \end{equation}
\item[$(2)$]  When $n=3$, for any $R\in [1,\infty)$, $k\in (\ff 3 2,\infty]$ and  $q\in [1,\ff 65),$
\beq\label{UL2} \lim_{t\to\infty} \sup_{\nu\in \scr P_{k,R}} \E^\nu\Big[ \big|\big\{t\W_2(\mu_t,\mu)^2- \Xi(t)\big\}^+
 + \big\{\Xi(t)- \kk t\W_2(\mu_t,\mu)^2\big\}^+\big|^q\Big] =0,\end{equation}
 so that when $\pp M$ is either empty or convex,
\beq\label{UL2'} \lim_{t\to\infty} \sup_{\nu\in \scr P_{k,R}} \E^\nu\Big[ \big|t\W_2(\mu_t,\mu)^2- \Xi(t)|^q\Big]
  =0.\end{equation}
 \item[$(3)$] For $n\le 3,$ we have  $\eta_Z\in (0,\infty)$ and
\beq\label{UL0} \lim_{t\to\infty}  \sup_{\nu\in \scr P}  \Big(\big\{t\E^\nu[W_2(\mu_t,\mu)^2]-\eta_Z\big\}^+
+ \big\{\eta_Z-\kk t\E^\nu[W_2(\mu_t,\mu)^2]\big\}^+\Big)=0.\end{equation}
In particular, when $\pp M$ is either empty or convex,
\beq\label{UL0'} \lim_{t\to\infty}  \sup_{\nu\in \scr P}   \big|t\E^\nu[W_2(\mu_t,\mu)^2]-\eta_Z\big|
 =0.\end{equation}
  \item[$(4)$] For $n=4$, there exist constants $c_1,c_2,t_0>0$ such that
  \beq\label{S4}   \ff{ c_1}{ t} \log (1+t) \le  \inf_{\nu\in \scr P}  \E^\nu[\W_2(\mu_t,\mu)^2] \le    \sup_{\nu\in \scr P}   \E^\nu [\W_2(\mu_t,\mu)^{2}]  \le  \ff{c_2}{ t }\log (1+t),\ \ t\ge t_0. \end{equation}
 \item[$(5)$] For $n\ge 5$,  there exist constants $c_1,c_2,t_0>0$ such that
 \beq\label{S5} \beg{split}  c_1 t^{-\ff 2 {d-2}}  \le  \inf_{\nu\in \scr P}  \big(\E^\nu[\W_1(\mu_t,\mu)] \big)^2  \le \sup_{\nu\in \scr P}   \E^\nu [\W_2(\mu_t,\mu)^{2}]  \le c_2 t^{-\ff 2 {d-2}},\ \ t\ge 1.  \end{split}\end{equation}
\end{enumerate}\end{thm}

\paragraph{Remark 1.1.} (1) By \eqref{ZZ} we have  $\eta_Z<\eta_0$ for $Z\ne 0$, so \eqref{UL0} and \eqref{UL0'} provide a precise characterization on  the acceleration of 
   a divergence-free perturbation $Z$ for  the convergence of   empirical measures in $\W_2$.

(2) When $Z=0$ (i.e. the symmetric case), \eqref{UL0}, \eqref{UL0'}, \eqref{S5} and the upper bound in \eqref{S4} have been presented in \cite{WZ},  which are covered by Theorem \ref{TN}. The $L^q$-convergence \eqref{UL}-\eqref{UL2'} appear here for the first time, which together with the lower bound in
 \eqref{S4} are  new also in the symmetric case.

(3)  It is proved in \cite{WZ} that for $M$ being the $4$-dimensional torus and $L=\DD$, there exists a constant $c>0$ such that
$$\inf_{\nu\in \scr P} \big(\E^\nu[\W_1(\mu_t,\mu)]\big)^2\ge c t^{-1}\log(1+t),\ \ t\ge 1.$$
We hope that this estimate also holds for general non-symmetric diffusions on $4$-dimensional compact manifolds, such that the lower bound estimate in \eqref{S4} is strengthened with $\W_1$ replacing $\W_2$.

\

In the next result,  we estimate $\big(\E[\W_{2p}(\mu_t,\mu)^{2 q}]\big)^{\ff 1 q}$   for  all $p,q\in [1,\infty)$. Besides   the critical phenomenon  in Theorem \ref{TN} with   the critical convergence rate $t^{-1}\log t$ for dimension $n=4$,  the critical rate   also appears to
dimensions $n=2,3$ with different $(p,q)$.

 \beg{thm}\label{T1} There exist     $c,t_0\in (0,\infty)$ and $\kk:  [1,\infty)\times [1,\infty)\to (0,\infty)$, such that the following assertions hold.
  \beg{enumerate}
 \item[$(1)$] When   $n=1$,  for any
 $(p,q)\in [1,\infty)\times [1,\infty),$
 \beq\label{S1}   \ff c t\le  \inf_{\nu\in \scr P} \big(\E^\nu[\W_1(\mu_t,\mu)]\big)^2
  \le    \sup_{\nu\in \scr P}   \big\{\E^\nu [\W_{2p}(\mu_t,\mu)^{2q}]\big\}^{\ff 1 q} \le \ff {\kk_{p,q}}t,\ \
    \ t\ge t_0.\end{equation}
 \item[$(2)$] Let $n=2$. Then   $\eqref{S1}$ holds for any    $p\in [1,\infty)$ and $q\in [1,\ff p{p-1})$. Next, for any $ p\in (1,\infty)$ and $ q=\ff p {p-1},$
   \beq\label{S2}   \sup_{\nu\in \scr P}   \big\{\E^\nu [\W_{2p}(\mu_t,\mu)^{2q}]\big\}^{\ff 1 q} \le
  \ff {\kk_{p,q}}t\log (1+t), \ \ t\ge t_0.\end{equation}
  Finally, for   any $ p\in (1,\infty)$ and $q\in (\ff p {p-1},\infty),$
  \beq\label{S3}   \sup_{\nu\in \scr P}   \big\{\E^\nu [\W_{2p}(\mu_t,\mu)^{2q}]\big\}^{\ff 1 q} \le \kk_{p,q}t^{-\ff {2}{n(3-p^{-1}-q^{-1})-2}},  \ \ t\ge 1.   \end{equation}
 \item[$(3)$] Let $n=3$.  Then  $\eqref{S1}$ holds for any    $p\in [1,\ff 32)$ and $q\in [1,\ff {3p}{5p-3})$;
  $\eqref{S2}$ holds for $p\in [1,\ff 32)$ and $q=\ff {3p}{5p-3}$;
  and $\eqref{S3}$ holds for any  $ p\in [1,\infty)$ and $q\in (\ff {3p}{5p-3},\infty)\cap [1,\infty).$
\item[$(4)$] Let $n=4$. Then $\eqref{S2}$ holds for $p=q=1$, and  $\eqref{S3}$ holds for any
       $ (p,q)\in   [1,\infty)\times [1,\infty)\setminus \{(1,1)\}.$
 \item[$(5)$] When $n\ge 5$, $\eqref{S3}$ holds for any   $(p,q)\in   [1,\infty)\times [1,\infty).$
\end{enumerate}\end{thm}

 \ \newline
{\bf D. Structure of the paper.}
 In Section 2, we state our main results for non-symmetric subordinated diffusion processes in an abstract framework.  In Sections 3 and 4, we prove the main results on upper and lower bound estimates respectively.  In Section 5, we apply the main results to some concrete models, where the result for the first model covers Theorems \ref{TN} and \ref{T1} as direct consequences with
 $B(\ll)=\ll$ (hence, $\aa=1$).

\section{Main results in an abstract framework}

We first introduce the framework of the study, then state the main results on  the Wasserstein distance of the empirical measures for non-symmetric  subordinated diffusion processes.

\subsection{The framework}

{\bf A. State space.}   Let  $(M,\rr)$ be a  length space,   let $\scr P$ be the set of all probability measures on $M$, let $\B_b(M)$ be  the class of bounded measurable functions on $M$, and let $C_{b,L}(M)$ be the set of all bounded Lipschitz continuous functions on $M$.
 For any
$p\in [1,\infty)$, the $L^p$-Wasserstein distance is defined as
\beq\label{WP}\W_p(\nu_1,\nu_2):= \inf_{\pi\in \C(\nu,\gg)} \bigg(\int_{M\times M} \rr(x,y)^p\pi(\d x,\d y)\bigg)^{\ff 1 p},
\ \ \nu_1,\nu_2\in \scr P,\end{equation}
where   $\C(\nu_1,\nu_2)$ is the  set of all couplings for $\nu_1$ and $\nu_2$.

\ \newline
{\bf B. Symmetric diffusion process.}
Let $\hat X_t$ be a reversible Markov process on $M$ with the unique invariant probability measure $\mu\in \scr P$ having full support. For any $q\ge p\in [1,\infty]$, let $\|\cdot\|_p$ be the norm in $L^p(\mu)$, and let $\|\cdot\|_{p\to q}$ the operator norm from $L^p(\mu)$ to $L^q(\mu)$.
Throughout the paper, we simply denote $\mu(f)=\int_M f\d\mu$ for $f\in L^1(\mu).$

The Markov semigroup $\hat P_t$ is formulated as
$$\hat P_t f(x)=\E^x[f(\hat X_t)],\ \ t\ge 0, x\in M, f\in \B_b(M),$$
where  and in the sequel,   $\E^x$ stands for the expectation
for the underlying Markov process starting at point $x$. In general, for any $\nu\in \scr P$,  $\E^\nu$ is  the expectation for the underlying   Markov process with initial distribution $\nu$. 

Let $(\hat \EE,\D(\hat \EE))$ and
$(\hat L, \D(\hat L))$ be, respectively,  the associated symmetric Dirichlet form and self-adjoint generator in $L^2(\mu)$.
We assume that $C_{b,L}(M)$ is a dense subset of $\D(\hat\EE)$ under the   $\hat \EE_1$-norm
$\|f\|_{\hat\EE_1}:= \ss{\mu(f^2)+\hat\EE(f,f)},$
and 
$$\hat\EE(f,g)=\int_M \GG(f,g)\d\mu,\ \ f,g\in C_{b,L}(M)$$ holds for
a symmetric  local square field (champ de carr\'e)
$$\GG: C_{b,L}(M)\times C_{b,L}(M)\to \B_b(M),$$
  such that for any $f,g,h\in C_{b,L}(M)$ and $\phi\in C_b^1(\R),$  we have
\beg{align*} &\ss{\GG(f,f)(x)}=|\nn f(x)|:=\limsup_{y\to x} \ff{|f(y)-f(x)|}{\rr(x,y)},\ \ x\in M,\\
&\GG(fg,h)= f\GG(g,h)+ g\GG(f,h),\ \ \ \GG(\phi(f), h)=  \phi'(f) \GG(f,h).\end{align*}
We also assume that $\hat L$ satisfies the chain rule
$$\hat L\Phi(f)= \Phi'(f)\hat Lf +\Phi''(f) |\nn f|^2,\ \ \ f\in \D(\hat L)\cap C_{b,l}(M), \Phi\in C^2(\R).$$

\ \newline {\bf C. Non-symmetric perturbation.} Let
$$Z: C_{b,L}(M)\to \B_b(M)$$
be a bounded vector field  with ${\rm div}_\mu Z=0$, i.e. it satisfies
\beg{align*} &Z(fg)= fZg+gZ f,\ \   Z(\phi(f))= \phi'(f) Zf,\ \ \ f,g  \in C_{b,L}(M),\ \phi\in C^1(\R),\\
  & \|Z\|_\infty:= \inf\big\{K\ge 0: \ |Zf|\le K |\nn f|,\ f\in C_{b,L}(M)\big\}<\infty,\\
 & \mu(Zf):=\int_M (Z f) \d\mu=0,\ \ \ f\in C_{b,L}(M). \end{align*}
 Consequently, $Z$ uniquely extends to a bounded linear operator from $\D(\hat \EE)$ to $L^2(\mu)$ with
\beq\label{RM} \mu(Zf) =0,\ \ \ f\in \D(\hat \EE),  \end{equation}
and
$$\EE(f,g):=\hat \EE(f,g)+\mu(f Zg),\ \ \ f,g\in \D(\EE)=\D(\hat\EE)$$
is a (non-symmetric) conservative   Dirichlet form with generator
$$L:=\hat L+Z,\ \ \ \D(L)=\D(\hat L),$$ which satisfies the chain rule
$$L\Phi(f)= \Phi'(f)Lf +\Phi''(f) |\nn f|^2,\ \ \ f\in \D(\hat L)\cap C_{b,l}(M), \Phi\in C^2(\R).$$
Assume that  $L$ generates a unique diffusion process  $X_t$ on $M$, such that the associated Markov semigroup is given by
$$P_t f(x)= \E^x[f(X_t)],\ \ x\in M, t\ge 0, f\in \B_b(M).$$
By Duhamel's formula,
\beq\label{00} P_tf= \hat P_tf+\int_0^t P_s\{Z\hat P_{t-s}f\}\d s,\ \ f\in \D(\hat \EE), t\ge 0.\end{equation}

{\bf C.  Subordination.}  Let {\bf B} be the set of Bernstein functions $B$ satisfying $B(0)=0$ and $ B(r)>0$ for $r>0.$
For each $B\in {\bf B}$, there exists a unique  stable increasing process $S_t^B$ on $[0,\infty)$ with Laplace transform
\beq\label{LT} \E[\e^{-r S_t^B}]= \e^{-B(r)t},\ \ \ t,r\ge 0.\end{equation}
Let $S_t^B$ be independent of $X_t$. We consider the subordinated diffusion process
$$X_t^B:= X_{S_t^B},\ \ t\ge 0,$$
and study  the convergence to $\mu$ in $\W_p$ ($p\ge 1$) for  the   empirical measure
$$\mu_t^B:=\ff 1 t \int_0^t \dd_{X_s^B}\d s,\ \ t>0.$$

We will mainly consider $\aa$-stable type time change  for $\aa\in [0,1]$, i.e. the Bernstein function $B$ is  in the   classes
\beg{align*} {\bf B}^\aa:=\Big\{B\in {\bf B}:\ \liminf_{r\to\infty}B(r)r^{-\aa}>0\Big\},\ \
 {\bf B}_\aa:=\Big\{B\in {\bf B}:\ \limsup_{r\to\infty}B(r)r^{-\aa}<\infty\Big\}.\end{align*}

\subsection{Upper bound estimates}
We make the following assumption, where  \eqref{A10} implies that the spectrum of $-\hat L$ is discrete and
all eigenvalues $\{\ll_i\}_{i\ge 0}$ listed in the increasing order counting multiplicities satisfy
$$\ll_i \ge c i^{\ff 2 d},\ \ i\in \mathbb Z_+$$ for some constant $c>0$, where $\ll_1\ge \ll$, see for instance   \cite{Davies}. In general,
  $\ll_i$ may increase faster than $i^{\ff 2 d}$, see for instance Subsection 5.2 where $d=n+2$ but $\ll_i\sim i^{\ff 2n}$,
we make the additional assumption $\eqref{A1-1}$.

\beg{enumerate} \item[$(A_1)$] Let $B\in {\bf B}^\aa$ for some $\aa\in [0,1].$ There exist  constants $c,\ll >0, d\ge d'\ge 1$ and a map $k: (1,\infty)\to (0,\infty)$ such that
\beq\label{A10} \|\hat P_t-\mu\|_{1\to\infty}\le c   t^{- \ff d 2} \e^{-\ll t},\ \ t>0,\end{equation}
\beq\label{A1-1} \ll_i\ge c i^{\ff 2 {d'}},\ \ \ i\in \mathbb Z_+,\end{equation}
\beq\label{A120} |\nn   \hat P_t f| \le   k(p)  ( \hat P_t |\nn f|^p)^{\ff 1 p},\ \ t\in [0,1],  p\in (1,\infty), f\in C_{b,L}(M).\end{equation}
\end{enumerate}
We will also need the following condition on  the continuity of $\hat X_t$.

\beg{enumerate} \item[$(A_2)$] For any $p\in [1,\infty)$ there exists a constant $c(p)>0$ such that
\beq\label{A13} \E^\mu\big[\rr(\hat X_0,\hat X_t)^p\big]\le c(p) t^{\ff p 2},\ \ t\in [0,1].\end{equation}
 \end{enumerate}

 Let $\{\phi_i\}_{i\ge 0}$ with $\phi_0\equiv 1$ be the   unitary eigenfunctions for $\{\ll_i\}_{i\ge 0},$ i.e.
 \beq\label{A1-3} \hat L\phi_i=-\ll_i\phi_i,\ \  \hat P_t\phi_i=\e^{-\ll_i t} \phi_i,\ \ \mu(\phi_i\phi_j)=1_{\{i=j\}},\ \ \ i,j\in \mathbb Z_+, t\ge 0.\end{equation}
 Let
 \beq\label{XIB}\Xi^B(t):=\sum_{i=1}^\infty \ff{\psi_i^B(t)^2}{\ll_i},\ \ \ \psi_i^B(t):=\ff 1 {\ss t}\int_0^t \phi_i(X_s^B)\d s,\ \ t>0, i\in\mathbb N.\end{equation}
 Let $\scr P_{k,R}$  be in \eqref{KR}, let
 \beq\label{QA} q_{\aa}:=     \ff{2d}{(2d+d' -2-2\aa)^+},\end{equation}
and  denote the integer parts of $q\in [1,\infty)$ by
 $$\i(q):= \sup\big\{i\in\mathbb N:  i\le q\big\}.$$
 The first main result of the paper is the following.

 \beg{thm}\label{TU1} Assume $(A_1)$ and $(A_2)$  with  $d'< 2(1+\aa).$
 \beg{enumerate} \item[$(1)$] If $q_\aa>\ff d{2\aa}$ and
$$\aa>\aa(d,d'):=\ff 1 4\Big(\ss{(2+d-d')^2+4d(d+d'-2)}+d'-d-2 \Big),$$
  then
 \beq\label{BU01} \lim_{t\to\infty} \sup_{\nu\in \scr P} \E^\nu\Big[\Big|\big\{t\W_2(\mu_t^B,\mu)^2-\Xi^B(t)\big\}^+\Big|^q\Big]=0,\ \ q\in [1, q_\aa).\end{equation}
 \item[$(2)$] For any $q\in [1,q_\aa)$ and $k\in (\ff d {2\aa \i(q)},\infty]\cap [1,\infty]$, where we set $(\ff d {2\aa \i(q)},\infty]=\{\infty\}$ if $\aa=0$,
 \beq\label{BU02} \lim_{t\to\infty} \sup_{\nu\in \scr P_{k,R}} \E^\nu\Big[\Big|\big\{t\W_2(\mu_t^B,\mu)^2-\Xi^B(t)\big\}^+\Big|^q\Big]=0,\ \   R\in (0,\infty).\end{equation}\end{enumerate}
 \end{thm}

 To estimate $ \E[\W_2(\mu_t^B,\mu)^{2}],$ we let
 $$\eta_Z^B:= \sum_{i=1}^\infty \ff{2{\bf V}_B(\phi_i)}{\ll_i},\ \ \ \ {\bf V}_B(\phi_i):= \int_0^\infty \mu(\phi_iP_s^B\phi_i)\d s.$$

 \beg{thm}\label{TU2} Assume $(A_1)$ and $(A_2)$. Let $q\in [1,\infty).$ 
 \beg{enumerate} \item[$(1)$] If $d'<2(1+\aa)$, then $\eta_Z^B<\infty$ and
 \beq\label{BU02'} \limsup_{t\to\infty}\sup_{\nu\in \scr P} t \E^\nu[\W_2(\mu_t^B,\mu)^2]\le \eta_Z^B.\end{equation}
 \item[$(2)$] Let  $d'\ge 2(1+\aa).$ Then there exists a constant $c>0$ such that
 for any $t\ge 1$,
 \beq\label{BUa} \sup_{\nu\in \scr P}  \E^\nu[\W_{2}(\mu_t^B, \mu)^{2}] \le \beg{cases} ct^{-1}\log (1+t),\ &\text{if}\ d'=2(1+\aa),\\
 c t^{-\ff{2}{d'-2\aa}},\ &\text{if}\ d'>2(1+\aa).\end{cases} \end{equation}\end{enumerate}
 \end{thm}

 To estimate $\W_{p}(\mu_t^B, \mu)$, we will use the $L^p$-boundedness of the Riesz transform $\nn (a_0-\hat L)^{-\ff 1 2}$ for some $a_0\ge 0$.
 According to \cite{AC}, together with  the non-degeneracy condition,  the volume doubling condition and the scaled   Poincar\'e inequality, \eqref{A120}   implies
 \beq\label{RS}  \|\nn (a_0-\hat L)^{-\ff 1 2}\|_p<\infty \ \text{for\ some }\ a_0\in [0,\infty) \ \text{and \ all} \ p\in (2,\infty).\end{equation}
 Under assumption $(A_1)$, let
 $$\gg_{\aa, p,q}:= \ff{d'}2+\ff d 2\big(2-p^{-1}-q^{-1}\big)-\aa-1,\ \
 p,q\in [1,\infty), \aa\in [0,1].$$

 \beg{thm}\label{TU3} Assume $(A_1)$ and $(A_2)$.
 \beg{enumerate}
 \item[$(1)$] If   $\gg_{\aa,p,q}<0$, then there exists a constant $c>0$ such that
 \beq\label{BUb} \sup_{\nu\in \scr P} \big(\E^\nu[\W_{2p}(\mu_t^B, \mu)^{2q}\big)^{\ff 1 q}\le c t^{-1},\ \ t\ge 1.\end{equation}
  \item[$(2)$]   If  $\gg_{\aa,p,q}\ge 0$, then for any $\gg>\gg_{\aa,p,q},$ there exists a constant $c>0$ such that
      \beq\label{BUc} \sup_{\nu\in \scr P} \big(\E^\nu[\W_{2p}(\mu_t^B, \mu)^{2q}]\big)^{\ff 1 q}\le c t^{-\ff 1{1+\gg}},\ \ t\ge 1.\end{equation}
   \item[$(3)$]  Let   $\eqref{RS}$ hold. If     $\gg_{\aa,p,q}\ge 0$, then  there exists a constant $c>0$ such that for any $t\ge 1$,
        \beq\label{BUd} \sup_{\nu\in \scr P} \big(\E^\nu[\W_{2p}(\mu_t^B, \mu)^{2q}]\big)^{\ff 1 q}\le \beg{cases}c t^{-1}\log(1+t),\ &\text{if}\ \gg_{\aa,p,q}=0,\\
        c t^{-\ff 1{1+\gg_{\aa,p,q}}},\ &\text{if}\ \gg_{\aa,p,q}>0.\end{cases}\end{equation}
 \end{enumerate} \end{thm}

 \subsection{Lower bound estimate}
To derive sharp lower bound for $\E[\W_2(\mu_t^B,\mu)^2]$, we make the following assumption.
 \beg{enumerate} \item[$(B)$] $(M,\rr)$ is a geodesic space, there exist constants $\theta,K>0$ and $m\ge 1$ such that
 \beq\label{B*} |\nn \hat P_t\e^f|^2\le (\hat P_t \e^f)\hat P_t(|\nn f|^2\e^f)+ K t^\theta \|\nn f\|_\infty^2 (\hat P_t\e^{2mf})^{\ff 1 m},
 \ \ t\in [0,1], f\in C_{b,L}(M),\end{equation}
 and there exists   a function $h\in C([0,1];[1,\infty))$ such that
\beq\label{B22}\W_2(\nu \hat P_r,\mu)^2\le h(r)\W_2(\nu,\mu)^2,\ \ \nu\in\scr P, r\in [0,1].\end{equation}
\end{enumerate}

When $M$ is a Riemannian manifold without boundary or with convex boundary, if the Bakry-Emery curvature of $\hat L$ is bounded below by a constant $-K$,
then $(B)$ holds for $m=1$ and $h(r)=\e^{2Kr}$, see for instance \cite[Theorem 2.3.3$(2')(9)$]{Wbook} or \cite{ST}.

 \beg{thm}\label{TL1} Assume $(A_1)$  and $(B)$  with  $d'< 2(1+\aa).$
 \beg{enumerate} \item[$(1)$] If $\aa>\aa(d,d')$ and $q_\aa>\ff d{2\aa},$ then
$$\lim_{t\to\infty} \sup_{\nu\in \scr P} \E^\nu\Big[\Big|\big\{th(0)\W_2(\mu_t^B,\mu)^2- \Xi^B(t)\big\}^-\Big|^q\Big]=0,\ \ q\in [1, q_\aa).$$
\item[$(2)$] For any $q\in [1,q_\aa)$ and $k\in (\ff d {2\aa \i(q)},\infty]\cap [1,\infty)$, where we set $(\ff d {2\aa \i(q)},\infty]=\{\infty\}$ if $\aa=0$,  
 $$ \lim_{t\to\infty} \sup_{\nu\in \scr P_{k,R}} \E^\nu\Big[\Big|\big\{th(0)\W_2(\mu_t^B,\mu)^2- \Xi^B(t)\big\}^-\Big|^q\Big]=0,\ \ q\in [1,q_\aa), R\in [1,\infty).$$
 \item[$(3)$]  
 $\liminf\limits_{t\to\infty} \sup\limits_{\nu\in \scr P} t\E^\nu [\W_2(\mu_{t}^B,\mu)^2]\ge h(0)^{-1} \eta_Z^B.$
 \end{enumerate}
 \end{thm}

The next result manages the critical case where the convergence rate of $\E[\W_2(\mu_t^B,\mu)^2]$
is at most $t^{-1}\log t$, correspondingly to \eqref{BUa} on the upper bound estimate.

 \beg{thm}\label{TL2} Assume $\eqref{A10},\eqref{A120}$, $(A_2)$,  $(B)$ and that
 \beq\label{A2'} k' i^{\ff 2 {d'}}\le \ll_i\le k i^{\ff 2 {d'}},\ \ \ i\in \mathbb N\end{equation} holds for some constants $k,k'>0$. If  $\aa':= \ff{d'}2-1\in [0,1]$ and  $B\in {\bf B}^\aa\cap {\bf B}_{\aa'}$ for  some   $\aa\in [0,\aa']\cap (\aa'-1,\aa'],$     then     there exist constants $c,t_0>0$ such that
 \beq\label{BLa} \inf_{\nu\in \scr P} \E^\nu[\W_2(\mu_t^B, \mu)^{2}]\ge  c t^{-1}\log (1+t),\ \ t\ge t_0.\end{equation}
   \end{thm}

Finally, we consider the lower bound estimate  on $\W_1$.

\beg{thm}\label{TL3} Let $B\in {\bf B}$. Then  the following assertions hold.
 \beg{enumerate}
 \item[$(1)$] Assume   $\eqref{A10}$, $\eqref{A120}$ and that  the completion $\bar M$ of $M$ is a Polish space.   Then
there exist    constants $c,t_0>0$ such that
\beq\label{LBD1} \inf_{\nu\in \scr P} \E^\nu[\W_1(\mu_t^B,\mu)]\ge c t^{-\ff 1 2},\ \ t\ge t_0.\end{equation}
\item[$(2)$]  Assume that  $\eqref{A13}$ holds for $p=1$, and  there exist constants $k,d''>0$ such that
\beq\label{MU} \sup_{x\in M} \mu(B(x,r))\le k r^{d''},\ \ r\ge 0,\end{equation}
where $B(x,r):=\{y\in M: \rr(x,y)\le r\}.$ If $B\in {\bf B}_\aa$ for some $\aa\in [0,1]$ with $d''>2(1+\aa)$, then there exist constants $c,t_0>0$ such that
\beq\label{LBD2} \inf_{\nu\in \scr P} \E^\nu[\W_1(\mu_t^B,\mu)]\ge c t^{-\ff 1 {d''-2\aa}},\ \ t\ge t_0.\end{equation}\end{enumerate}
 \end{thm}

 \section{Proofs of Theorems \ref{TU1}-\ref{TU3}}
For a density function $f$ with respect to $\mu$, let
$(f\mu)(A):=\int_A f\d\mu$ for a measurable set $A\subset M$.
 Recall that for any probability density functions $f,f_1,f_2\in L^2(\mu)$, we have
 \beq\label{AM0} \W_2(f\mu,\mu)^2\le \int_M\ff{|\nn \hat L^{-1} (f-1)|^2}{\scr M(f)}\,\d\mu,\ \ \scr M(f):= 1_{\{f>0\}}\ff{f-1}{\log f},\end{equation}
 \beq\label{Ledoux} \W_p(f_1\mu,f_2\mu)^p \le p^p \int_M \ff{|\nn \hat L^{-1} (f_1-f_2)|^p}{f_2^{p-1}}\,\d\mu.\end{equation}
 These estimates
   have been   presented in \cite{AMB} and \cite{Ledoux} respectively by using the Kantorovich dual formula  and   Hamilton-Jacobi equations, which are available when $(M,\rr)$ is a   length space as we assumed, see \cite{V}.

   Since the empirical measure $\mu_t^B$ is singular with respect to $\mu$, to apply these estimates we make the following regularization of $\mu_t^B$:
   \beq\label{GPP} \mu_{t,r}^B:= f_{t,r}^B\mu,\ \ f_{t,r}^B:=1+\ff 1 {\ss t} \sum_{i=1}^\infty   \e^{-\ll_i r} \psi_i^B(t)\phi_i,\ \ t,r>0,\end{equation}
   where $\psi_i^B(t):= \ff 1 {\ss t} \int_0^t \phi_i(X_s^B)\d s.$ Letting $\nu \hat P_r$ being the distribution of
   $\hat X_r$ with initial distribution $\nu$, by \eqref{A1-3} and the spectral representation
   $$\hat p_r(x,y)=1+\sum_{i=1}^\infty \e^{-\ll_i r}\phi_i(x)\phi_i(y)$$
   for the heat kernel $\hat p_r$  of $\hat P_r$ with respect to $\mu$, we have
   \beq\label {GPP'} \mu_{t,r}^B= \mu_t^B\hat P_r,\ \ \ t,r>0.\end{equation}
 So,   \eqref{AM0} implies
\beq\label{AM} \W_2(\mu_{t,r}^B,\mu)^2\le \int_M\ff{|\nn \hat L^{-1} (f_{t,r}^B-1)|^2}{\scr M(f_{t,r}^B)}\,\d\mu,\ \ t,r>0.\end{equation}
According to  \eqref{A10}, we have  $\lim_{t\to\infty}f_{t,r}^B\to 1$   so that
$\lim_{t\to\infty} \scr M(f_{t,r}^B)=1$. When the convergence is fast enough,
 \eqref{A1-3}  and \eqref{GPP} would imply that for large enough $t$, $t\W_2(\mu_{t,r}^B,\mu)^2$ is bounded above  by
\beq\label{BXR} \Xi_r^B(t):= t\mu(|\nn \hat L^{-1} (f_{t,r}^B-1)|^2) =\sum_{i=1}^\infty \ff{\e^{-2\ll_i r}}{\ll_i} \psi_i^B(t)^2,\ \ \ t,r>0.\end{equation}
On the other hand, by  \eqref{Ledoux} we have
\beq\label{Ledoux'} \W_p(\mu_{t,r}^B,\mu)^p \le p^p \mu\big(|\nn \hat L^{-1} (f_{t,r}^B-1)|^p\big),\ \ t,r>0, p\in [1,\infty).\end{equation}
With the above observations, and noting that 
$\W_p(\mu_t^B,\mu)\le \W_p(\mu_{t,r}^B,\mu)+ \W_p(\mu_t^B,\mu_{t,r}^B),$ to estimate $\W_p(\mu_t^B,\mu)$  we  present some lemmas on $\Xi_r^B(t)$, $\mu\big(|\nn \hat L^{-1} (f_{t,r}^B-1)|^p\big)$ and
  $\W_p(\mu_{t,r}^B,\mu_t^B)$ respectively. 

 \subsection{Some lemmas}

To apply \eqref{Ledoux'}, we need estimate $\|\nn \hat L^{-1}(f_{t,r}^B-1)\|_p$, see \eqref{A11} below. To this  end, and also for later use,
we first estimate  $\|P_t^B-\mu\|_{p\to q}$  and $\|P_t Z\|_{2p}$   for $q\ge p\ge 1$, where
$P_t^B$ is  the Markov semigroup for the subordinated diffusion process $X_t^B$   given by
\beq\label{SB} P_t^B f(x):= \E^x[f(X_t^B)]= \E[P_{S_t^B} f(x)],\ \ \ t\ge 0, x\in M, f\in \B_b(M).\end{equation}

 \beg{lem}\label{L0} Assume $\eqref{A10}$. Then there exists a possibly different constant $\ll\in (0,\ll_1]$ such that the following assertions hold.
 \beg{enumerate}\item[$(1)$] Let $B\in {\bf B}^\aa$ for some $\aa\in (0,1]$. Then there exists a constant   $k>0$   such that
 \beq\label{UAA} \|P_t^B-\mu\|_{p\to q} \le kt^{-\ff{d(q-p)}{2pq\aa}} \e^{-\ll t},\ \ t>0, q\ge p\in [1,\infty].\end{equation}
 \item[$(2)$] Let $\eqref{A120}$ hold. Then for any $p\in [1,\infty)$ there exists a constant $c(p)>0$ such that
     \beq\label{A12} \|\nn P_t f\|_p\le c(p)   t^{-\ff 1 2}\e^{-\ll t} \|f\|_p,\ \ t>0,f\in C_{b,L}(M),\end{equation}
 \beq\label{A12'}  \|P_t(Zf)\|_{2p}\le c(p) \|Z\|_\infty t^{-\ff 1 2} \e^{-\ll t}\|f\|_{2p},\ \ f\in \D(\hat\EE)\cap L^{2p}(\mu), t>0.\end{equation}
 Moreover, for any $\kk\in (0,\ff 1 2)$ there exists a constant $c(p,\kk)>0$ such that
 \beq\label{A11}\|\nn\hat L^{-1} f\|_{2p}\le c(p,\kk) \|(-\hat L)^{\ff{d(p-1)}{4p}-\kk}f\|_2,\ \ \mu(f)=0.\end{equation}
 \end{enumerate}
 \end{lem}

\beg{proof} (a) We will use some known results on functional inequalities which can be found in e.g. \cite{Wbooka}. Firstly,  since $\ll_1>0,$ we have the Poincar\'e inequality
\beq\label{PC} \mu(f^2)\le \ff 1 {\ll_1} \hat\EE(f,f),\ \ f\in\D(\hat\EE), \mu(f)=0.\end{equation}
By \eqref{RM} we have $\EE(f,f)=\hat\EE(f,f)$. So, \eqref{PC}  implies
\beq\label{EXP}  \|P_t-\mu\|_2\le \e^{-\ll_1 t},\ \ \ t\ge 0.\end{equation}
Next, according to \cite[Theorem 3.3.14 and 3.3.15]{Wbooka}, \eqref{A10} implies the super Poincar\'e inequality
\beq\label{SPC} \mu(f^2)\le r \hat\EE(f,f) + c_1 (1+r^{-\ff d 2})\mu(|f|)^2,\ \ r>0, f\in \D(\hat\EE)\end{equation}
for some constant $c_1>0$, which further yields
\beq\label{EX1} \|P_t\|_{1\to\infty}\le c_2 (1\land t)^{-\ff d 2},\ \ t>0\end{equation}
for some constant $c_2>0$. Noting that $B\in {\bf B}^\aa$ implies
\beq\label{BAA} B(r)\ge k 1_{\{r\ge \ll_1\}}r^\aa,\ \ \ r\ge 0\end{equation} for  some constant $k>0$,
by \eqref{LT} we find a constant $c_3>0$ such that
$$\E[(S_t^B)^{-d}]=\ff{\E\int_0^\infty r^{d-1}\e^{-r S_t^B}\d r}{\int_0^\infty r^{d-1}\e^{-r}\d r}
\le c_3(1\land t)^{-\ff\d \aa},\ \ t>0.$$
Combining this with \eqref{LT}, \eqref{SB}, \eqref{EXP} and \eqref{EX1}, we find   constants $c_4,c_5\ge 1$ such that
\beg{align*} &\|P_t^B-\mu\|_{1\to\infty} \le \E[\|P_{S_t^B}-\mu\|_{1\to\infty}] \le c_4\E\Big[\big\{1+(S_t^B)^{-\ff d 2}\big\}\e^{-\ll_1 S_t^B}\Big]\\
&\le 2 c_4 \big(\E[1+(S_t^B)^{-d}]\big)^{\ff 1 2}
 \big(\E[\e^{-2\ll_1 S_t^B}]\big)^{\ff 1 2}\le c_5 (1\land t)^{-\ff d{2\aa}} \e^{-B(2\ll_1)t/2},\ \ t>0.\end{align*}
By the interpolation theorem, this and $\|P_t^B-\mu\|_p\le 2$ for $p\in [1,\infty]$ imply \eqref{UAA} for
some constants $c,\ll>0$. In particular, for $B(r)=r$ and  $Z=0$ or $Z\ne 0$, \eqref{UAA} implies  to
\beq\label{UAA'} \|\hat P_t-\mu\|_{p\to q}\lor \|P_t-\mu\|_{p\to q}\le k  t^{-\ff {d(q-p)}{2pq}}\e^{-\ll t},\ \ t>0, q\ge p\ge 1.\end{equation}

(b) To prove \eqref{A12}, we first prove that for some decreasing $c: (1,\infty)\to (0,\infty),$
\beq\label{A121} |\nn \hat P_tf|\le \ff {c(p)}{\ss t} (\hat P_t |f|^p)^{\ff 1 p},\ \ t\in (0,1], f\in C_{b,L}(M).\end{equation}
 By H\"older's inequality and   $f= f^+-f^-$, it suffices to prove for $p\in (1,2]$ and $f\ge 0$. Moreover, by first using $f+\vv$ replacing $f$ for $\vv>0$ then letting $\vv\downarrow 0$, we may and do assume that $\inf f>0.$

 By \eqref{A120} we have $\hat P_{t-s}f\in \D(\hat L)\cap C_{b,L}(M)$ for $s\in [0,t)$. By the chain rule, we obtain
 \beg{align*} &\ff{\d}{\d s}  \hat P_s (\hat  P_{t-s}f)^p  = \hat P_s \hat L (\hat  P_{t-s}f)^p - \hat P_s \big\{
 p(\hat P_{t-s}f)^{p-1} \hat L \hat P_{t-s} f\big\}\\
 &= p(p-1) \hat P_s\big\{ (\hat  P_{t-s} f)^{p-2} |\nn \hat P_{t-s} f|^2\big\},\ \ s\in [0,t).\end{align*}
So, for $p\in (1,2]$, we have
\beq\label{F2} \beg{split} &I:=  \hat P_t f^p -(\hat P_t f)^p= \int_0^t \ff{\d}{\d s}  \hat P_s (\hat  P_{t-s}f)^p \d s\\
&= p(p-1)\int_0^t \hat P_s\big\{ (\hat  P_{t-s} f)^{p-2} |\nn \hat P_{t-s} f|^2\big\}\d s.\end{split}\end{equation}
By the H\"older/Jensen inequalities, we obtain
\beg{align*}& \big[\hat P_s |\nn \hat P_{t-s}f|^p\big]^{\ff 2 p} \le \big[\hat P_s\{(P_{t-s} f)^{p-2} |\nn \hat P_{t-s}f|^2 \} \big] \big\{\hat P_s (\hat P_{t-s} f)^p\big\}^{\ff{2-p} p}\\
&\le  \big[\hat P_s\{(P_{t-s} f)^{p-2} |\nn \hat P_{t-s}f|^2 \} \big]  \big(\hat P_t f^p\big)^{\ff{2-p} p}.\end{align*}
 Combining this with \eqref{F2} and \eqref{A120} where we may assume that $k(p)$ is decreasing in $p$ due to Jensen's inequality, we find increasing  $C: (1,\infty)\to (0,\infty)$  such that
\beg{align*} &I\ge p(p-1)\int_0^t \big(\hat  P_s|\nn \hat  P_{t-s} f|^p\big)^{\ff 2 p} (\hat P_t f^p)^{\ff{p-2}p}\d s\\
&\ge C(p) \int_0^t |\nn \hat P_t f|^2 (\hat P_tf^p)^{\ff {p-2}p} \d s= C(p) t |\nn\hat P_t f|^2 (\hat P_tf^p)^{\ff {p-2}p}.\end{align*}
This implies \eqref{A121} for  some decreasing $c: (1,\infty)\to (0,\infty)$.

Next,    we  intend to prove that for some   constant $c>0$,
\beq\label{AANN}   \|\nn P_t f\|_\infty \le c\|f\|_\infty t^{-\ff 1 2},\ \ t\in (0,1], f\in \B_b(M).\end{equation}
For any $x\ne y\in M$, let
$$h_t(x,y):= \sup_{\|g\|_\infty\le 1} \ff{|P_t g(x)-P_tg(y)|}{\rr(x,y)},\ \ t\ge 0.$$
By \eqref{00} and \eqref{A121}, we find a constant $c_1>0$ such that
$$h_t(x,y)\le c_1 t^{-\ff 1 2}  +\int_0^t c_1 (t-s)^{-\ff 1 2}  h_s(x,y)\d s,\ \ t\in (0,1].$$
By the generalized Gronwall inequality, see \cite{Ye}, this implies \eqref{AANN}.

Moreover, by \eqref{A121},   the $L^p$-contraction of $P_t$ and $\hat P_t$, and the Duhamel's formula
$$P_tf= \hat P_tf+\int_0^t \hat P_s(Z P_{t-s}f)\d s,$$ we obtain
\beg{align*} &\|\nn P_t f\|_p\le \|\nn \hat P_tf\|_p+ \int_0^t \|\nn \hat P_s(Z P_{t-s}f)\|_p\d s\\
&\le c(p) t^{-\ff 1 2}\|f\|_p +\int_0^t c(p)\|Z||_\infty s^{-\ff 1 2} \|\nn P_{t-s}f\|_p\d s,\ \ t>0.\end{align*}
When $f\in \B_b(M)$, by \eqref{AANN} and the generalized Gronwall  inequality, this imply \eqref{A12} for
$t\in (0,1].$

Finally, by \eqref{UAA'} for $p=q$ such that
  $$\|  P_t-\mu\|_p\le k \e^{-\ll t}.$$
  Combining this with the semigroup property and \eqref{A12} for $t\in (0,1]$,   for any $t>1$ we have
$$ \|\nn   P_t f\|_p =\|\nn   P_1( P_{t-1} f-1)\|_p\le c(p) \| P_{t-1}(f-1)\|_p
\le c(p) k \e^{-\ll (t-1)} \|f\|_p.$$
So, \eqref{A12}   also  holds for $t>1$ and some constant $\ll>0$.

(c)  Let $P_t^*$ be the $L^2(\mu)$-adjoint operator of $P_t$. By \eqref{RM}, $P_t^*$ is the diffusion semigroup generated by
$L^*:= \hat L-Z$, and satisfies
\beq\label{FM} P_t^*g= \hat P_tg -\int_0^t \hat P_s(Z P_{t-s}^*g)\d s,\ \ \ t>0,\ \ g\in L^2(\mu).\end{equation}
Let $g\in \D(\hat \EE)$ with $\|g\|_{\ff {2p}{2p-1}}\le 1$. We have
$$\|\nn P_t^*g\|_{\ff{2p}{2p-1}}^2\le \|\nn P_t^*g\|_2^2=\hat \EE(P_t^*g, P_t^*g) \le \hat \EE(g,g)<\infty,\ \ t\ge 0,$$
so that  \eqref{RM}, \eqref{A12} and \eqref{FM} yield that for some constant $c_1>0$,
$$\|\nn P_t^*g\|_{\ff {2p}{2p-1}}\le c_1t^{- \ff 1 2} + c_1\int_0^t s^{-\ff 1 2} \|\nn P_{t-s}^* g\|_{\ff {2p}{2p-1}}\d s<\infty,\ \ t\in (0,1].$$ By the generalized Gronwall inequality, see \cite{Ye}, we find a constant $c_2>0$ such that
$$\sup_{g\in \D(\hat\EE), \|g\|_{\ff {2p}{2p-1}}\le 1} \|\nn P_t^*g\|_{\ff{2p}{2p-1}}\le c_2 t^{-\ff 1 2},\ \ t\in (0,1].$$
Combining this with the semigroup property and  \eqref{UAA'}, we find a constant $c_3>0$ such that
$$\sup_{g\in \D(\hat\EE), \|g\|_{\ff {2p}{2p-1}}\le 1} \|\nn P_t^*g\|_{\ff{2p}{2p-1}}\le c_3 t^{-\ff 1 2}\e^{-\ll t},\ \ t>0.$$
Thus, by \eqref{RM}, for any $f\in \D(\hat \EE)\cap L^{2p}(\mu)$ we have
\beg{align*} &\|P_t (Zf)\|_{2p} = \sup_{g\in \D(\hat \EE), \|g\|_{\ff {2p}{2p-1}}\le 1}  \big|\mu((P_t^*g)(Zf))\big|\\
&= \sup_{g\in \D(\hat \EE), \|g\|_{\ff {2p}{2p-1}}\le 1} \big|\mu(f (ZP_t^*g))\big|\le c_3 \|Z\|_\infty t^{-\ff 1 2} \e^{-\ll t} \|f\|_{2p},\ \ t>0.\end{align*}
Therefore, \eqref{A12'} holds for some constants $c(p),\ll>0$.

(d) Noting that
$$(-\hat L)^{-(1+\ff{d(p-1)}{4p}-\kk)} =\ff 1 {\GG(1+\ff{d(p-1)}{4p}-\kk)} \int_0^\infty s^{\ff{d(p-1)}{4p}-\kk} \hat P_s\d s,$$
by   \eqref{A12} and \eqref{UAA'}, we find constants $c_1,c_2,c_3>0$
such that
\beg{align*} & \|\nn\hat L^{-1} f\|_{2p}= \|\nn (-\hat L)^{-(1+\ff{d(p-1)}{4p}-\kk)}(-\hat L)^{\ff{d(p-1)}{4p}-\kk}f\|_{2p}\\
  &\le \ff 1 {\GG(1+\ff{d(p-1)}{4p}-\kk)}  \int_0^\infty s^{\ff{d(p-1)}{4p}-\kk}\|\nn \hat P_{s/2}\{\hat P_{s/2}(-\hat L)^{\ff{d(p-1)}{4p}-\kk} f\}\|_{2p}\d s\\
  &\le c_1  \int_0^\infty s^{\ff{d(p-1)}{4p}-\kk-\ff 1 2}\e^{-\ll s/2} \|\hat P_{s/2}(-\hat L)^{\ff{d(p-1)}{4p}-\kk} f\|_{2p}\d s\\
  &\le c_2\int_0^\infty   s^{-(\kk+\ff 1 2)}\e^{-\ll s} \|(-\hat L)^{\ff{d(p-1)}{4p}-\kk} f\|_{2}\d s
\le c_3 \|(-\hat L)^{\ff{d(p-1)}{4p}-\kk} f\|_{2},\end{align*}
where the last step is due to $\ff 1 2+\kk<1$.
Thus,  \eqref{A11} holds.

\end{proof}

Next, we present some consequence of $(A_1)$ and $(A_2)$.

\beg{lem}\label{L00} We have the following assertions.
\beg{enumerate}\item[$(1)$] If $\eqref{A10}$ holds, then $(M,\rr)$ is bounded, i.e.
\beq\label{RD} D:=\sup_{x,y\in M}\rr(x,y)<\infty.\end{equation}
\item[$(2)$] If $\eqref{A13}$ holds, then  for any $p,q\in [1,\infty)$, there exists
a constant $c>0$ such that
\beq\label{XB1} \big(\E^\mu[\W_{2p}(\mu_t^B,\mu_{t,r}^B)^{2q}]\big)^{\ff 1 q}\le c r,\ \ r\in (0,1].\end{equation}
\item[$(3)$]  If  $\eqref{A120}$ and  $\eqref{B*}$ hold, then there exist constants $\kk_0,\kk_1>0$ such that
\beq\label{A21}\beg{split}& |\nn P_t \e^f|^2\le (P_t\e^f)P_t (|\nn f|^2\e^f) +\kk_1 t^{\theta}\|\nn f\|_\infty^2 (P_t\e^f)^2,\\
&\qquad \qquad \text{for} \ t\in [0,1],\ f\in C_{b,L}(M)\ \text{with}\  t\|\nn f\|_\infty^2\le \kk_0.\end{split}\end{equation}\end{enumerate}
\end{lem}

\beg{proof} (1)    According to \cite[Theorem 3.3.15(2)]{Wbooka}, \eqref{A10} implies  the super Poincar\'e inequality
$$\mu(f^2)\le r \hat \EE(f,f)+(1+r^{-\ff d 2})\mu(|f|)^2,\ \ r>0, f\in \D(\hat \EE),$$
which further implies \eqref{RD} due to \cite[Theorem 3.3.20]{Wbooka}.

(2) By Jensen's inequality, we only need to prove \eqref{XB1} for $q\ge p\ge 1.$
Recall that $\dd_{X_s^B}\hat P_r$ is the distribution of $\hat X_r$ with initial value $X_s^B$,   we have
$$\pi_t:= \ff 1 t \int_0^t \big\{\dd_{X_s^B}\times (\dd_{X_s^B}\hat P_r)\big\}\d s\in \C(\mu_t^B,\mu_{t,r}^B), $$
so that
$$ \W_{2p}(\mu_t^B,\mu_{t,r}^B)^{2p}\le \int_{M\times M} \rr(x,y)^{2p}\pi_t(\d x,\d y)=\ff 1 t \int_0^t \E^x[\rr(x,\hat X_r)^{2p}]\big|_{x=X_s^B} \d s.$$
Noting that $\E^\mu=\int_M\E^x\mu(\d x)$ and $\hat X_t$ is stationary with initial distribution $\mu$, by combining this with   Jensen's inequality and \eqref{A13}, we obtain
\beg{align*}
&\E^\mu[\W_{2p}(\mu_t^B,\mu_{t,r}^B)^{2q}]\le
\E^\mu\bigg[\ff 1 t \int_0^t \E^x[\rr(x,\hat X_r)^{2q}]\big|_{x=X_s^B} \d s\bigg]\\
&=\ff 1 t\int_0^t \E^\mu[\rr(\hat X_0,\hat X_r)^{2q}]\d s\le c(2q) r^q,\ \ t>0, r\in (0,1].\end{align*}
So, \eqref{XB1} holds.

(3) Let $f\in C_{b,L}(M).$ By \eqref{AANN}, we have $(P_{t-s}\e^f)^{2m}\in \D(L)\cap C_{b,L}(M)$ for
$s\in [0,t)$, so that the chain rule implies
\beg{align*}& \ff{\d}{\d s} P_s (P_{t-s} \e^f)^{2m}= P_s L (P_{t-s} \e^f)^{2m} - P_s\big\{ 2m
(P_{t-s} \e^f)^{2m-1} LP_{t-s} \e^f\big\}\\
& = 2m(2m-1) P_s\big\{(P_{t-s} \e^f)^{2m-2} |\nn P_{t-s} \e^f|^2\big\},\ \ s\in [0,t).\end{align*}
By combining this with \eqref{A120} for $p=2$ and Jensen's inequality, we find a constant $c_1>0$ such that
 \beg{align*}& P_t\e^{2mf} -(P_t\e^f)^{2m} =\int_0^t \ff{\d}{\d s} P_s (P_{t-s} \e^f)^{2m} \d s\\
 &=\int_0^t P_s \big\{2m(2m-1) (P_{t-s} \e^f)^{2m-2} |\nn P_{t-s} \e^f|^2\big\}\d s\\
 &\le c_1 \|\nn f \|_\infty^2 \int_0^t P_s\big\{(P_{t-s}\e^f)^{2m-2} P_{t-s}\e^{2f}\big\}\d s
 \le c_1 t \|\nn f\|_\infty^2 P_t \e^{2mf}. \end{align*}
 Taking $\kk_0= \ff 1 {2c_1}$ such that  $t \|\nn f\|_\infty^2\le \kk_0 $  implies
 $c_1 t \|\nn f\|_\infty^2\le\ff 1 2,$ we derive
 \beq\label{21}P_t\e^{2mf}\le 2(P_t\e^f)^{2m}.\end{equation}
 Combining this with \eqref{B*}, we obtain \eqref{A21} for some constant $\kk_1>0$.
\end{proof}

  Noting that \eqref{A1-3} and \eqref{GPP} imply
\beq\label{GPP2} \|(-\hat L)^{\bb} (f_{t,r}^B-1)\|_2^2=\ff 1 t \sum_{i=1}^\infty \ll_i^{2\bb}\e^{-2\ll_i r} \psi_i^B(t)^2,\ \ r,t>0, \bb\in\R,\end{equation}
to bound $\E^\nu[\W_p(\mu_{t,r}^B,\mu)^{2q}]$ from above using  \eqref{Ledoux'} and \eqref{A11}, we   estimate  $\E^\nu[|\psi_i^B|^{2q}]$ as follows.

\beg{lem}\label{BLN0} Assume \eqref{A10} and let $B\in {\bf B}^\aa$ for some $\aa\in [0,1].$  Then: 
\beg{enumerate}\item[$(1)$] For any $q\in [1,\infty),$ there exists a constant $c(q)>0$ such that
\beq\label{BLN10} \sup_{t>0}\E^\nu[|\psi_i^B(t)|^{2q}] \le c(q) \|h\|_\infty\ll_i^{\ff{d(q-1)}{2}-q\aa}, \
\ \ i\in \mathbb N, \nu= h\mu.\end{equation}
\item[$(2)$]   For any $q\in [1,\infty)$ and
$k\in (\ff d {2\aa \i(q)}, \infty]\cap [1,\infty]$, there exists a constant $c(q,k)>0$ such that
 \beq\label{BLN1}  \E^\nu[|\psi_i^B(t)|^{2q}]
 \le c(q,k) \|h\|_k(1\land t)^{-\ff d{2\aa k}}\ll_i^{\ff{d(q-1)}{2}-q\aa},
\ \ i\in \mathbb N, \nu= h\mu, t>0.\end{equation}
Moreover, if $\i(q)>\ff d {2\aa}$, then  there exists a constant $c(q)>0$ such that
\beq\label{BLN1*}  \sup_{\nu\in \scr P} \E^\nu[|\psi_i^B(t)|^{2q}]
 \le c(q) (1\land t)^{-\ff {d}{2\aa}}
 \ll_i^{\ff{d(q-1)} 2  -q\aa},
\ \ i\in \mathbb N,  t>0.\end{equation}\end{enumerate}\end{lem}

\beg{proof} (1) Let $h_{i,\aa}(t):= \min\big\{(\ff 1 2\land t)^{-\ff 1 {2\aa}},\ \ll_i^{\ff 1 2}\big\}.$
When $\aa>0$,   for any $k>0$ there exist  constants $a_1,a_2>0$   such that
\beq\label{XX}\beg{split} &\int_0^\infty h_{i,\aa}(t)\e^{-kt}\d t
\le \int_0^{\ll_i^{-\aa}} \ll_i^{\ff 1 2}\d t
+ \int_{\ll_i^{-\aa}}^\infty \Big(t^{-\ff 1 {2\aa}}+ 2^{\ff 1 {2\aa}}\Big)\e^{-kt}\d t\\
&\le \ll_i^{\ff 1 2-\aa} + a_1 \ll_i^{(\ff 1 {2}-\aa)^+}\big[1+1_{\{\aa=\ff 1 2\}}\log(1+\ll_i)\big]\\
&\le a_2 \ll_i^{(\ff 1 2-\aa)^+} \big[1+1_{\{\aa=\ff 1 2\}}\log(1+\ll_i)\big],\ \ \ i\in \mathbb N.\end{split}\end{equation}
When $\aa=0$ we have $h_{i,\aa}(t)=\ll_i^{\ff 1 2}$ so that this estimate holds as well.

 We first prove the following estimate  for some constants $k_1,k_2>0$:
\beq\label{PBT} \|P_t^B\phi_i-\e^{-B(\ll_i)t}\phi_i\|_{2q}\le k_1\|Z\|_\infty \ll_i^{\ff{d(q-1)}{4q}-1}   h_{i,\aa}(t)\e^{-k_2 t},\ \ t>0, i\in \mathbb N, q\in [1,\infty].\end{equation}
By
\eqref{A1-3} and \eqref{UAA'}, we find  constants $c_1,c_2>0$ such that
\beq\label{M2} \|\phi_i\|_{2q}=\inf_{s>0} \|\hat P_s \phi_i\|_{2q}\e^{\ll_i s} \le c_1 \inf_{s\in (0,1]}
s^{-\ff{d(q-1)}{4q}}\e^{\ll_i s}\le c_2 \ll_i^{\ff{d(q-1)}{4q}},\ \ i\in \mathbb N, q\in [1,\infty].\end{equation}
By \eqref{00}, \eqref{SB} and \eqref{A1-3}, we obtain
\beq\label{PTP} P_t^B \phi_i= \E\bigg[\e^{-\ll_i S_t^B} \phi_i + \int_0^{S_t^B}\e^{-\ll_i (S_t^B-s)}P_s(Z\phi_i)\,\d s\bigg],\ \ t>0.\end{equation}
Combining this with \eqref{LT},   \eqref{A12'} and \eqref{M2}, we derive
\beq\label{V}   \|P_t^B \phi_i-\e^{-B(\ll_i)t}\phi_i\|_{2q}\le   c_2 c(q) \|Z\|_\infty\ll_i^{\ff{d(q-1)}{4q}}\E\int_0^{S_t^B}
\e^{-\ll_i(S_t^B-s)} s^{-\ff 1 2} \e^{-\ll s}\d s. \end{equation}
Noting that $\ll_i\ge\ll_1\ge \ll$ implies
$$-\ll_i(S_t^B-s)-\ll s \le -\ff {\ll_i}2 (S_t^B-s) -\ff{\ll} 2 (S_t^B-s)-\ll s=
-\ff{\ll_i}2 (S_t^B-s)-\ff\ll 2 S_t^B-\ff{\ll} 2 s,$$
by the FKG inequality, we find a constant $c_3>0$ such that
\beq\label{V*1}\beg{split}& \int_0^{S_t^B}\e^{-\ll_i(S_t^B-s)} s^{-\ff 1 2} \e^{-\ll s}\d s
\le \e^{-\ll S_t^B/2}\int_0^{S_t^B}\e^{-\ll_i(S_t^B-s)/2} s^{-\ff 1 2}\e^{-\ll s/2}\d s\\
&\le \e^{-\ll S_t^B/2}\bigg(\int_0^{S_t^B}\e^{-\ll_i(S_t^B-s)/2}\d s\bigg) \ff 1 {S_t^B} \int_0^{S_t^B}
 s^{-\ff 1 2} \e^{-\ll s}\d s
\le \ff{c_3}{\ll_i} \e^{-\ll S_t^B/2} (S_t^B)^{-\ff 1 2}.\end{split}\end{equation}
Moreover, by \eqref{LT} and \eqref{BAA},
we find constants $c_4,c_5,c_6>0$ such that
\beq\label{V*2}\beg{split}
&\E\big[(S_t^B)^{-\ff 1 2} \e^{-\ll S_t^B/2}\big]= \E\bigg[\ff{\e^{-\ll S_t^B/2}}{\GG(1/2)}  \int_0^\infty u^{-\ff 1 2} \e^{-u S_t^B}\d u\bigg]\\
&=  \ff 1 {\GG(1/2)}  \int_0^\infty u^{-\ff 1 2} \e^{-B(u+\ll/2)t}\d u\le\ff 1 {\GG(1/2)}  \int_0^\infty u^{-\ff 1 2} \e^{-B(u)t/2-B(\ll/2)t/2}\d u \\
&\le c_4\e^{-B(\ll/2)t/2}\bigg[ \int_0^{\ll_1}u^{-\ff 1 2}\d u +  \int_{\ll_1}^\infty  u^{-\ff 1 2}
\e^{-ku^\aa t }\d u\bigg]\le c_5 \Big(\ff 1 2\land t\Big)^{-\ff 1 {2\aa}}
\e^{-c_6 t}.\end{split}\end{equation}
Combining this with \eqref{V} and \eqref{V*1}, we find constants $k_1,k_2>0$ such that
\beq\label{3.31*} \big\|P_t^B\phi_i-\e^{-B(\ll_i)t} \phi_i\big\|_{2q} \le k_1 \|Z\|_\infty \ll_i^{\ff{d(q-1)}{4q}-1} \Big(\ff 1 2\land t\Big)^{-\ff 1 {2\aa}} \e^{-k_2 t}.\end{equation}
On the other hand, by \eqref{A12} and \eqref{M2}, we find constants $c_1',c_2'>0$ such that
\beq\label{3.32*} \|\nn \phi_i\|_{2q} =\inf_{s>0} \|\nn \hat P_s \phi_i\|_{2q}\e^{\ll_i s}
\le \inf_{s>0} c_1' s^{-\ff 1 2} \e^{\ll_i s} \|\phi_i\|_{2q} \le c_2'\ll_i^{\ff{d(q-1)}{4q}+\ff 1 2}.\end{equation}
So, instead of \eqref{V}, this and \eqref{UAA'} imply
\beq\label{3.33*}   \big\|P_t^B\phi_i-\e^{-B(\ll_i)t} \phi_i\big\|_{2q}\le\|Z\|_\infty \|\nn \phi\|_{2q} \E\int_0^{S_t^B} \e^{-\ll_i(S_t^B-s)-\ll s}\d s.\end{equation}
By $\ll_i\ge \ll$ and \eqref{LT}, we obtain
$$\E\int_0^{S_t^B} \e^{-\ll_i(S_t^B-s)-\ll s}\d s\le \E\bigg[\e^{-\ll S_t^B/2}\int_0^{S_t^B} \e^{-\ll_i(S_t^B-s)/2}\d s\bigg] \le \ff 2{\ll_i} \e^{-B(\ll/2)t}.$$
Combining this with \eqref{3.32*} and \eqref{3.33*}, we find constants $k_1,k_2>0$ such that
 $$ \big\|P_t^B\phi_i-\e^{-B(\ll_i)t} \phi_i\big\|_{2q}]\le k_1 \|Z\|_\infty\ll_i^{\ff{d(q-1)}{4q}-\ff 1 2}
 \e^{-k_2t}.$$
 This together with \eqref{3.31*} implies \eqref{PBT}.

Next, we   prove \eqref{BLN10} for $q\in \mathbb N$. By \cite[(2.14)]{WW} for $f=\phi_i$, we find a constant $k_0>0$ such that
\beq\label{M00} \E^\nu\big[|\psi_i^B(t)|^{2q}\big]\le k_0 \bigg(\ff 1 t \int_0^t \d s_1\int_0^{s_1} \big\{\E^\nu[|\phi_i P^B_{s_1-s} \phi_i|^q](X_s^B)\big\}^{\ff 1 q} \d s\bigg)^q.\end{equation}
By $\nu= h\mu$ and the Markov property, we obtain
\beq\label{NBM} \beg{split}& \E^\nu\big[|\phi_iP_{s_1-s}^B\phi_i|^q(X_s^B)\big]=\mu\big(h P_s^B|\phi_iP_{s_1-s}^B\phi_i|^q\big)\le \|h\|_k\big\|P_s^B|\phi_iP_{s_1-s}^B\phi_i|^q\big\|_{\ff k {k-1}}\\
&\le \|h\|_k \|P_s^B\|_{1\to\ff k{k-1}}\|\phi_iP_{s_1-s}^B\phi_i\|_{q}^q
\le \|h\|_k \|P_s^B\|_{1\to\ff k{k-1}}\|\phi_i\|_{2q}^q\|P_{s_1-s}^B\phi_i\|_{2q}^q,\ \ k\in [1,\infty].\end{split}\end{equation}
Taking $k=\infty$ and combining with  \eqref{PBT}, \eqref{M2}  and  \eqref{M00}, we find constants $k_1,k_2>0$ such that
$$\E^\nu \big[|\psi_i^B(t)|^{2q}\big]\le k_1 \|h\|_\infty \ll_i^{\ff{d(q-1)}2} \bigg(\ff 1 t \int_0^t \d s_1
\int_0^{s_1} \big(\e^{-B(\ll_i)(s_1-s)} +\ll_i^{-1} h_{i,\aa}(s_1-s)\e^{-k_2(s_1-s)}\big)\d s\bigg)^q.$$
Combining this with \eqref{BAA}, which together with \eqref{XX} implies
\beq\label{*DD} \int_0^\infty \Big[\e^{-B(\ll_i)t}+\ll_i^{-1}h_{i,\aa}(t)\e^{-k_2t}\Big]\d t
\le c \ll_i^{-\aa},\ \ i\in \mathbb N\end{equation} for some constant $c>0$,
we derive \eqref{BLN10} for $q\in \mathbb N.$

Finally, for any $q\in (1,\infty)$, let $\i(q)$ be the integer part of $q$. By \eqref{BLN10} for $\i(q)$ and $1+\i(q)$ replacing $q$ which have just been proved, and using H\"older's inequality, we find a constant $c(q)>0$ such that
\beq\label{HD}\beg{split} &\E^\nu \big[|\psi_i^B(t)|^{2q}\big]\le \big(\E^\nu \big[|\psi_i^B(t)|^{2\i(q)}\big]\big)^{\i(q)+1-q} \big(\E^\nu \big[|\psi_i^B(t)|^{2+2\i(q)}\big]\big)^{q-\i(q)}\\
&\le c(q) \|h\|_\infty \ll_i^{\ff 1 2\{d(\i(q)-1)(\i(q)+1-q)+ d\i(q) (q-\i(q))\}-\aa\{\i(q)(\i(q)+1-q)+(1+\i(q))(q-\i(q))\}}\\
&= c(q) \|h\|_\infty \ll_i^{\ff{d(q-1)}2-q\aa},\ \ \ t>0, i\in\mathbb N.\end{split}\end{equation}
Then \eqref{BLN10} is proved.

(2) By the same reason leading to \eqref{HD}, we only need to prove \eqref{BLN1} and \eqref{BLN1*}  for $q\in \mathbb N$ so that $\i(q)=q$.
Let $k\in (\ff d{2\aa q}, \infty]\cap [1,\infty]$. By   \eqref{PBT}, \eqref{M00} and \eqref{NBM},
we find a constant  $c_1>0$ such that
\beq\label{M11}\beg{split} &\E^\nu \big[|\psi_i^B(t)|^{2q}\big]
 \le c_1 \|h\|_k \ll_i^{\ff{d(q-1)}2} \\
 &\qquad\times \bigg(\ff 1 t \int_0^t \d s_1 \int_0^{s_1} (1\land s)^{-\ff d{2kq\aa}}
\Big[\e^{-B(\ll_i)(s_1-s)} + \ff{h_{i,\aa}(s_1-s)}{\ll_i \e^{k_2 (s_1-s)}}\Big]\d s\bigg)^q.\end{split}\end{equation}
By the FKG inequality and \eqref{*DD}, we find a constant $c_2>0$ such that
\beq\label{LIN}\beg{split} & \int_0^{s_1} (1\land s)^{-\ff d{2kq\aa}}
\big(\e^{-B(\ll_i)(s_1-s)} + \ll_i^{-1}h_{i,\aa}(s_1-s) \e^{-k_2 (s_1-s)}\big)\d s\\
&\le \bigg(\ff 1 {s_1} \int_0^{s_1} (1\land s)^{-\ff d{2kq\aa}}\d s\bigg)
\int_0^{s_1}\big(\e^{-B(\ll_i)(s_1-s)} + \ll_i^{-1}h_{i,\aa} \e^{-k_2 (s_1-s)}\big)\d s\\
&\le c_2 (1\land s_1)^{-\ff d{2kq\aa}}\ll_i^{-\aa},\ \ \ t>0, i\in \mathbb N.\end{split}\end{equation}
This together with \eqref{M11} yields
$$\E^\nu \big[|\psi_i^B(t)|^{2q}\big]\le c\|h\|_k (1\land t)^{-\ff d{2 \aa k} }
\ll_i^{\ff{d(q-1)}2-q\aa}, \ \ s_1>0, i\in \mathbb N$$
for some constant $c>0$. Therefore, \eqref{BLN1} holds for $q\in \mathbb N$.

It remains to prove \eqref{BLN1*}  for $\ff d{2\aa}<  q\in \mathbb N.$ By \eqref{UAA}, \eqref{PBT} and \eqref{M2},   we find   constants $c_1>0$ such that
\beg{align*}&\sup_{\nu\in \scr P}\big( \E^\nu \big[|\phi_i P_{s_1-s}^B \phi_i|^q(X_s^B)\big]\big)^{\ff 1 q}
= \sup_{\nu\in \scr P} \Big\{\nu\big(P_s^B|\phi_i P_{s_1-s}^B\phi_i|^q\big)\Big\}^{\ff 1 q}\\
&\le \|P_s^B\|_{1\to\infty}^{\ff 1 q} \|\phi_iP_{s_1-s}^B \phi_i\|_q\le c_1(1\land s)^{-\ff d{2q\aa}} \|\phi_i\|_{2q}\|P_{s_1-s}^B\phi_i\|_{2q}\\
&\le c_1(1\land s)^{-\ff d{2q\aa}} \ll_i^{\ff{d(q-1)}{2q}} \big(\e^{-B(\ll_i)(s_1-s)} + \ll_i^{-1} h_{i,\aa}(s_1-s) \e^{-k_2(s_1-s)}\big).\end{align*}
Combining this with \eqref{M00} and \eqref{LIN} for $k=1$, we get  \eqref{BLN1*}  for $\ff d{2\aa}< q\in \mathbb N.$
\end{proof}

\beg{lem}\label{LY} Assume $(A_1)$, $(A_2)$. Let $a_0\in [0,\infty), q\in [1,\infty), \bb\in\R.$
\beg{enumerate}
\item[$(1)$] For any  $k\in (\ff d{2\aa \i(q)}, \infty]\cap [1,\infty]$  where $k=\infty$ if $\aa=0$, and for any $ R\in [1,\infty)$, there exists a constant $c>0$ such that
 \beq\label{LY1} \beg{split} &\sup_{t\ge 1, \nu\in \scr P_{k,R}} t^q\E^\nu\Big[\big\|(a_0-\hat L)^{\ff 1 2} (-\hat L)^{\bb-\ff 1 2}(f_{t,r}^B-1)\big\|_2^{2q}\Big]\\
 & \le c \big[r^{-(2\bb+\ff{d'}2 +\ff{d(q-1)}{2q}-\aa)^+}
   +1_{\{\aa   = 2\bb   + \ff{ d'}2+\ff{d(q-1)}{2q}\}} \log(1+r^{-1})\big]^q,\ \ r\in (0,1].\end{split}\end{equation}
 \item[$(2)$] If $\i(q)>\ff d{2\aa}$, then there exists a constant $c>0$ such that
 \beq\label{LY2} \beg{split} &\sup_{t\ge 1, \nu\in \scr P}t^q \E^\nu\Big[\big\|(a_0-\hat L)^{\ff 1 2} (-\hat L)^{\bb-\ff 1 2}(f_{t,r}^B-1)\big\|_2^{2q}\Big]\\
 & \le c \big[r^{-(2\bb+\ff{d'}2 +\ff{d(q-1)}{2q}-\aa)^+}
   +1_{\{\aa   = 2\bb   + \ff{ d'}2+\ff{d(q-1)}{2q}\}} \log(1+r^{-1})\big]^q,\ \ r\in (0,1].\end{split}\end{equation}
 \end{enumerate}\end{lem}

 \beg{proof}  By \eqref{BAA}, \eqref{GPP2}, H\"older's inequality   and \eqref{BLN1}, we find a constant $c_1>0$ such that
 \beg{align*}&I:=\sup_{t\ge 1, \nu\in \scr P_{k,R}} t^q\E^\nu\Big[\big\|(a_0-\hat L)^{\ff 1 2} (-\hat L)^{\bb-\ff 1 2}(f_{t,r}^B-1)\big\|_2^{2q}\Big]\\
 &= \sup_{t\ge 1, \nu\in \scr P_{k,R}}\E^\nu \bigg[\bigg|\sum_{i=1}^\infty \Big(\ff{\ll_i+a_0}{\ll_i}\Big)\ll_i^{2\bb}\e^{-2\ll_i r}\psi_i^B(t)^2\bigg|^q\bigg]\\
 &\le\Big(\ff{\ll_1+a_0}{\ll_1}\Big)^q\sup_{t\ge 1, \nu\in \scr P_{k,R}}\Big(\sum_{i=1}^\infty \ll_i^\theta \e^{-2\ll_ir}\Big)^{q-1}
 \sum_{i=1}^\infty \ll_i^{2q\bb-\theta(q-1)}\e^{-2\ll_ir} \E^\nu[|\psi_i^B(t)|^{2q}]\\
 &\le c_1 \Big(\sum_{i=1}^\infty \ll_i^\theta \e^{-2\ll_ir}\Big)^{q-1}
 \sum_{i=1}^\infty \ll_i^{2q\bb-\theta(q-1)+\ff{d(q-1)}2-q\aa}\e^{-2\ll_ir},\ \ \ \theta\in \R.\end{align*}
 Taking
 \beq\label{TH} \theta:= 2\bb +\ff{d(q-1)}{2q}-\aa,\end{equation}
 so that $\theta=  2q\bb-\theta(q-1)+\ff{d(q-1)}2-q\aa,$ and noting that
 $$\ll_i^{\theta^+} \e^{-\ll_i r} \le \sup_{s>0} s^{\theta^+}\e^{-s r} \le c  r^{-\theta^+},\ \ r\in (0,1]$$
 holds for some constant $c>0$ depending on $\theta^+$, we find a constant $c_2>0$ such that
\beg{align*} &I\le c_1 \Big(\sum_{i=1}^\infty \ll_i^\theta \e^{-2\ll_ir}\Big)^{q}
\le c_2r^{-q\theta^+}\Big(\sum_{i=1}^\infty \ll_i^{-\theta^-}\e^{-\ll_i r}\Big)^q.\end{align*}
On the other hand, by \eqref{A1-1} and the integral transform $t=rs^{\ff 2 {d'}}$,  we find constants $c_3,c_4, c_5>0$ such that
\beq\label{CCV}  \beg{split}& \sum_{i=1}^\infty\ll_i^{-\theta^-}\e^{-\ll_i r}\le c_3 \int_1^\infty s^{-\ff{2\theta^-}{d'}} \e^{-c_3 r s^{\ff 2 {d'}}}\d s 
 = \ff{c_3 d'} 2 \int_r^\infty r^{\theta^- -\ff{d'}2} t^{\ff{d'}2 -\theta^--1}\e^{-c_3 t}\d t\\
 &\le c_5 \Big\{r^{-(\ff{d'}2-\theta^-)^+}+1_{\{\ff{d'}2=\theta^-\}}\log(1+r^{-1})\Big\},\ \ \ \ r\in (0,1].\end{split}\end{equation}
Thus, we find a constant $c_6>0$ such that
\beg{align*}&I\le c_6 r^{-q \theta^+} \big[r^{-(\ff{d'}2-\theta^-)^+}+1_{\{\ff{d'}2=\theta^-\}}\log(1+r^{-1})\big]^q\\
&= c_6\big[r^{-(\ff{d'}2+\theta)^+}+ 1_{\{\ff{d'}2+\theta=0\}}\log(1+r^{-1})\big]^q,\ \ r\in (0,1].\end{align*}
This together with \eqref{TH} implies \eqref{LY1}.

When $\i(q)>\ff d{2\aa}$, \eqref{LY2} can be proved in the same way by using \eqref{BLN1*} replacing \eqref{BLN1}.

 \end{proof}

We are now ready to  show  that  as $t\to\infty$, $\E[\Xi_r^B(t)]$ converges to
\beq\label{ETB} \eta_{Z,r}^B:= \sum_{i=1}^\infty \ff{2\e^{-2\ll_i r}}{\ll_i} {\bf V}_B(\phi_i),\ \ \ r>0.\end{equation}

\beg{lem}\label{L0} Assume $(A_1)$, and let $R\in [1,\infty)$. \beg{enumerate}
\item[$(1)$] There exists a constant $c>0$ such that
\beq\label{VB} \eta_{Z,r}^B\le c\sum_{i=1}^\infty \ll_i^{-1-\aa} \e^{-2r\ll_i},\ \ r>0.\end{equation}
Consequently, $\eta_Z^B<\infty$ provided $\sum_{i=1}^\infty \ll_i^{-1-\aa}<\infty.$
\item[$(2)$] There exists a constant $c>0$ such that
\beq\label{BL-0} \sup_{\nu\in \scr P_{\infty,R}}\big|\E^\nu[\Xi_r^B(t)]-\eta_{Z,r}^B\big|\le \ff c t \sum_{i=1}^\infty \ll_i^{-1-\aa}\e^{-2\ll_i r},\ \ \ t,r>0.\end{equation}
Consequently, when $\sum_{i=1}^\infty \ll_i^{-1-\aa}<\infty,$  there exists a constant $c'>0$ such that
\beq\label{BH1} \sup_{\nu\in \scr P_{\infty,R}} \E^\nu[\Xi^B(t)]\le \eta_Z^B+\ff{c'}t<\infty,\ \ t>0.\end{equation}
\item[$(3)$]   For any $k>\ff d {2\aa}$, there exists a constant $c>0$ such that
    \beq\label{BH2} \sup_{t\ge 1,\nu\in \scr P_{k,R}} \E^\nu[\Xi_r^B(t)]\le c_2 \big\{r^{-(\ff{d'}2 -1-\aa)^+}
+1_{\{d'=2(1+\aa)\}} \log (1+r^{-1})\big\},\ \  r\in (0,1].\end{equation}
    \end{enumerate}
\end{lem}

\beg{proof} (1) By \eqref{PBT} for $q=1$ and \eqref{*DD},   we
find constants $c_1,c_2 >0$ such that
$${\bf V}_B(\phi_i):=\int_0^\infty \mu(\phi_iP_t^B\phi_i)\d t\le c_1 \int_0^\infty (\e^{-B(\ll_i)t}+\ll_i^{-1}h_{i,\aa}(t)\e^{-k_2 t})\d t\le  c_2\ll_i^{-\aa},\ \ i\ge 1.$$
This together with \eqref{ETB} implies \eqref{VB}. By the dominated convergence theorem with   $r\to 0$, the claimed consequence   follows from \eqref{VB}.

(2) By \eqref{GPP} and the Markov property, we obtain
\beq\label{PM}\beg{split} &\E^\nu[\psi_i^B(t)^2] =\ff 2 t \int_0^t \d s_1\int_0^{s_1} \E^\nu[\phi_i(X_{s_1}^B)\phi_i(X_s^B)]\d s\\
&=\ff 2 t \int_0^t \d s_1\int_0^{s_1} \nu\big(P_s^B\{\phi_i P_{s_1-s}^B\phi_i\}\big)\d s.\end{split}\end{equation}
Since $\nu\in \scr P_{\infty,R}$ implies $\nu=h\mu$ with $\|h-1\|_\infty\le R+1$, by    \eqref{UAA} and \eqref{PBT}, we find   constants  $c_1,c_2 >0$ such that
\beg{align*} &\sup_{\nu\in \scr P_{\infty,R}} \big|\nu\big(P_s^B\{\phi_i P_{s_1-s}^B\phi_i\}\big)-\mu(\phi_iP_{s_1-s}^B\phi_i)\big|
 =\sup_{\nu\in \scr P_{\infty,R}} \big|\mu\big(\{(P_s^B)^* h-1\} \phi_iP_{s_1-s}^B\phi_i)\big|\\
&\le  c_1\e^{-c_2 s}  \|\phi_iP_{s_1-s}^B\phi_i\|_1 \le c_1 \e^{-c_2 s} \|P_{s_1-s}^B \phi_i\|_2
 \le  c_1 \e^{-c_2 s} \big(\e^{-B(\ll_i)(s_1-s)} +\ll_i^{-1}h_{i,\aa}(s_1-s)\e^{-c_2 (s_1-s)}\big).\end{align*}
Combining this with \eqref{PM} and \eqref{*DD} for $k_2=\ff 1 2c_2$, and noting that
$c_2 s+c_2(s_2-s)\ge \ff 1 2 c_2 s_1+ \ff 1 2 c_2(s_1-s),$
we find a constant $c_3>0$ such that
\beq\label{01}\beg{split} &\sup_{\nu\in \scr P_{\infty,R}}\bigg|t\E^\nu[\Xi_r^B(t)]- \sum_{i=1}^\infty \ff{2\e^{-2\ll_i r}}{\ll_i t}\int_0^t \d s_1 \int_0^{s_1} \mu(\phi_i P_{s_1-s}^B\phi_i)\d s\bigg|\\
&\le  \sum_{i=1}^\infty \ff{2\e^{-2\ll_i r}}{\ll_i t}\int_0^t\d s_1 \int_0^{s_1} c_1 \e^{-c_2 s} \big(\e^{-B(\ll_i)(s_1-s)} +\ll_i^{-1}h_{i,\aa}(s_1-s)\e^{-c_2 (s_1-s)}\big)\d s\\
&\le \ff{c_3} t \sum_{i=1}^\infty \ll_i^{-1-\aa}\e^{-2\ll_i r},\ \ \ t,r>0. \end{split}\end{equation}
Similarly, by \eqref{PBT} for $q=1$ and \eqref{*DD}, we find constants $c_4,c_5,c_6>0$ such that
\beg{align*} &\bigg|\int_0^{s_1} \mu\big(\phi_iP_{s_1-s}^B\phi_i\big)\d s-{\bf V}_B(\phi_i)\bigg|\le
\int_{s_1}^\infty \big|\mu\big(\phi_iP_{s}^B\phi_i\big)\big|\le \int_{s_1}^\infty\|P_s^B\phi_i\|_2\d s\\
&\le c_4  \int_{s_1}^\infty\big(\e^{-B(\ll_i)s}+ \ll_i^{-1} h_{i,\aa}(s_1-s)\e^{-k_2s}\big)\d s\le
c_6 \ll_i^{-\aa} \e^{-c_5 s_1},\ \ \ s_1>0, i\in \mathbb N. \end{align*}
This together with \eqref{01} implies \eqref{BL-0} for some constant $c>0.$

(3) Let $k>\ff d {2\aa}.$ By \eqref{BLN1} for $q=1$, and \eqref{CCV} for $\theta^-=1+\aa$, we find   constants $c_1, c_2>0$ such that
\beg{align*} &\sup_{t\ge 1,\nu\in \scr P_{k,R}} \E^\nu[\Xi_r^B(t)] = \sup_{t\ge 1,\nu\in \scr P_{k,R}}\sum_{i=1}^\infty
\ff{\e^{-2\ll_i r}}{\ll_i}  \E^\nu[\psi_i^B(t)^2] \\
&\le c_1 \sum_{i=1}^\infty\ff{\e^{-2\ll_i r}}{\ll_i^{1+\aa}}\le c_2 \big\{r^{-(\ff{d'}2 -1-\aa)^+}
+1_{\{d'=2(1+\aa)\}} \log (1+r^{-1})\big\},\ \ \ r\in (0,1].\end{align*}
Then the proof is finished.
\end{proof}

Finally, to get rid of  the term $\scr M(f_{t,r}^B)$ from \eqref{AM}, we present one more lemma.

\beg{lem}\label{BLN4}   Assume $(A_1)$ and $(A_2)$ with  $d'<2(1+\aa)$. Then the following assertions hold.
\beg{enumerate}\item[$(1)$] There exists a constant $c>0$ and $\si\in (0,1)$ such that
\beq\label{XB4'} \E^\mu[\mu(|f_{t,r}^B-1|^2)]\le c t^{-1} r^{-\si},\ \ t\ge 1, r\in (0,1].\end{equation}
\item[$(2)$] There exists a constant $\gg>1$ such that
\beq\label{XB2'} \lim_{t\to\infty} \E^\mu\big[\mu\big(|\scr M(f_{t,t^{-\gg}}^B)^{-1}-1|^q\big)\big]=0,\ \ q\in [1,\infty).\end{equation}
\item[$(3)$]   We have $q_\aa>1$  and
\beq\label{XB3} \sup_{t,r>0}t^q\E^\mu\big[\mu\big(|\nn \hat L^{-1} (f_{t,r}^B-1)|^{2q}\big)\big]<\infty,\ \ q\in [1,q_\aa).\end{equation}
\end{enumerate}\end{lem}

\beg{proof} (1) By \eqref{GPP}, \eqref{BLN10} and \eqref{BAA}, we find   constants $c_1,c_2>0$ such that
$$t\E^\mu[\mu(|f_{t,r}^B-1|^2)]\le c_1 \sum_{i=1}^\infty \e^{-2\ll_i r} \E^\mu[\psi_i^B(t)^2]
\le c_2 \sum_{i=1}^\infty \e^{-2\ll_i r} \ll_i^{-\aa},\ \ t,r>0.$$
Since $d'<2(1+\aa)$ implies $\ff {d'}2-\aa<1$, combining this with \eqref{CCV} we derive \eqref{XB4'} for
any $\si\in (\ff {d'}2 -\aa, 1).$

(2) Let $\si\in (0,1)$ be in \eqref{XB4'} and take $\theta\in (0, q^{-1}(1-\si))$ for fixed $q\in [1,\infty).$ According to \cite[Lemma 3.2]{WW}, for any $\eta\in (0,1)$ and $\dd(\eta):= |(1-\eta)^{-\ff 1 2}-\ff 2 {\eta+2}|,$ we have
\beg{align*} &\E^\mu\big[|\scr M(f_{t,r}^B(y))^{-1}-1|^q \big]\\
&\le \dd(\eta) + (1+\theta^{-1} r^{-\ff\theta 2})^q \E^\mu[1_{\{|f^B_{t,r}(y)-1|>\eta\}}],\ \
t\ge 1, r\in (0,1], y\in M.\end{align*}
Next, by \eqref{XB4'} and Chebyshev's inequality, we obtain
$$\int_M\E^\mu[1_{\{|f^B_{t,r}(y)-1|>\eta\}}]\mu(\d y) \le \eta^{-2} \E^\mu[\mu(|f_{t,r}^B-1|^2)]\le c \eta^{-2} t^{-1} r^{-\si}.$$
Putting these two estimates together, we find a function $C: (0,1)\to (0,\infty)$ such that
$$ \E^\mu\big[\mu\big(|\scr M(f_{t,r}^B)^{-1}-1|^q\big)\big]\le \dd(\eta)+C(\eta) t^{-1}r^{-(\si+\ff {\theta q} 2)},\ \ t\ge 1, r,\eta\in (0,1).$$
Noting that $\si\in (0,1)$ and $\theta\in (0, q^{-1}(1-\si)) $ imply $\theta':= \si+\ff{\theta q} 2\in (0,1)$, by taking $r=t^{-\gg}$ for $\gg\in (1, \ff 1 {\theta'}),$ and letting first $t\to\infty$ then $\eta\to 0$, we derive \eqref{XB2'}.

(3) By   \eqref{QA}, $d'<2(1+\aa)$ implies  $q_\aa>1$.
  For any $q\in [1,q_\aa)$, we have
$$\ff{d(q-1)}{2q}-1<\aa -\ff{d'}2 -\ff{d(q-1)}{2q}.$$
So,  there exists $\kk\in (0,\ff 1 2)$ such that
$\bb:= \ff{d(q-1)}{4q}-\kk $ satisfies
\beq\label{NB} 2\bb < \aa-\ff{d'}2 -\ff{d(q-1)}{2q}.\end{equation}
By \eqref{A11}, \eqref{GPP}, \eqref{BAA} and \eqref{LY1}, we find constants $c_1,c_2>0$ such that
\beg{align*} & t^q\E^\mu\big[\mu\big(|\nn \hat L^{-1} (f_{t,r}^B-1)|^{2q}\big)\big]\le c_1 t^q
\E^\mu\big[\|(-\hat L)^{\ff{d(q-1)}{4q}-\kk} (f_{t,r}^B-1)\|_2^{2q} \big]\\
&\le c_2    r^{-(2\bb+\ff{d'}2 +\ff{d(q-1)}{2q}-\aa)^+} =c_2,\end{align*}
where the last step follows from \eqref{NB}. Then \eqref{XB3} holds.

\end{proof}

\subsection{Proof of Theorem \ref{TU1}}

We first consider the stationary case where the initial distribution is the invariant measure $\mu$, then
extend to more general setting by using an approximation argument.
To this end, we need the following further modification of the empirical measure:
\beq\label{FT}\d\tt\mu_{t,r}^B= \tt f_{t,r}^B\d\mu,\ \ \ \tt f_{t,r}^B:= (1-r)f_{t,r}^B+ r,\ \ \ r\in (0,1], t>0.\end{equation}
By \eqref{Ledoux}, we have
\beq\label{XB1'} \W_2(\mu_{t,r}^B,\tt\mu_{t,r}^B)^2\le 4 r \Xi_r^B(t),\ \ \ t>0, r\in (0,1].\end{equation}

\beg{prp}\label{P1} Assume $(A_1)$ and $(A_2)$ with $d'<2(1+\aa)$. Then
\beq\label{U02} \lim_{t\to\infty}\E^\mu\Big[\big|\{t\W_2(\mu_t^B,\mu)^2-\Xi^B(t)\}^+\big|^q\Big]=0,\ \ q\in [1,q_\aa).\end{equation}
\end{prp}

\beg{proof} Let $t\ge 1$ and take $r_t=t^{-\gg}$ for $\gg>1$ in \eqref{XB2'}. By \eqref{BXR} and \eqref{FT}, we have
$$t\mu(|\nn\hat L^{-1}(\tt f_{t,r_t}^B-1)|^2)= (1-r_t)^2\Xi_{r_t}^B(t),$$ so that   \eqref{AM}, \eqref{XB3} and \eqref{XB2'} yield
\beg{align*} & \lim_{t\to\infty}\E^\mu\Big[\big|\{t\W_2(\tt \mu_{t, r_t}^B,\mu)^2-(1-r_t)^2\Xi_{r_t}^B(t)\}^+\big|^q\Big]\\
&= \lim_{t\to\infty}\E^\mu\Big[\big|\{t\W_2(\tt \mu_{t, r_t}^B,\mu)^2-t\mu(|\nn\hat L^{-1}(\tt f_{t,r_t}^B-1)|^2)\}^+\big|^q\Big]\\
&\le  \lim_{t\to\infty}t^q \E^\mu\Big[\big\{ \mu\big(|\nn\hat L^{-1}(\tt f^B_{t,r_t}-1)|^2|\scr M(\tt f_{t,r_t}^B)^{-1}-1|\big)\big\}^q \Big]\\
&\le \lim_{t\to\infty}t^q \E^\mu\Big[ \mu\big(|\nn\hat L^{-1}(\tt f_{t,r_t}^B-1)|^{2q}|\scr M(\tt f_{t,r_t}^B)^{-1}-1|^q \big)\Big]\\
&\le \lim_{t\to\infty}\Big(t^{q'} \E^\mu\Big[ \mu\big(|\nn\hat L^{-1}(\tt f_{t,r_t}^B-1)|^{2q'}\big)\Big]\Big)^{\ff q {q'}} \Big(\E^\mu \Big[\mu\big(|\scr M(\tt f_{t,r_t}^B)^{-1}-1|^{\ff{qq'}{q'-q}} \big)\Big]\Big)^{\ff{q'-q}{q'}}\\
&=0,\ \  \ q'\in (q,q_\aa).\end{align*}
Noting that \eqref{XIB} and \eqref{BXR} imply
$$\Xi^B(t)\ge \Xi_{r_t}^B(t)\ge (1-r_t)^2\Xi_{r_t}^B(t),$$ we derive
\beq\label{U01} \lim_{t\to\infty}\E^\mu\Big[\big|\{t\W_2(\tt \mu_{t, r_t}^B,\mu)^2- \Xi^B(t)\}^+\big|^q\Big]
=0,\ \ q\in [1,q_\aa).\end{equation}
On the other hand, noting that
  $q\in [1, \ff d{(d+d'-2-2\aa)^+})$ implies $1+\aa-\ff{d'}2-\ff{d(q-1)}{2q}>0$,
by  \eqref{BXR}, \eqref{Ledoux'}  and \eqref{LY1} with $a_0=0$ and $\bb=-\ff 1 2,$  we obtain
\beq\label{UIB}\beg{split} &\sup_{r>0,t\ge 1} 4^{-q} t^q \E^\mu[\W_2(\mu_{t,r}^B,\mu)^{2q}]
\le  \sup_{r>0,t\ge 1}t^q\E^\mu\big[\|(-\hat L)^{-\ff 1 2} (f_{t,r}^B-1)\|_2^{2q}\big]\\
&=\sup_{r>0,t\ge 1} \E^\mu\big[\Xi^B_r(t)^q\big]
<\infty,\ \  \ 1\le q< \ff d{(d+d'-2-2\aa)^+}.\end{split}\end{equation}
Since   $q\in [1,q_\aa)$  implies $d'<2(1+\aa)$ and  $q\in [1, \ff d{(d+d'-2-2\aa)^+}),$
by combining this with \eqref{XB1}, \eqref{XB1'}, the triangle inequality and $r_t=t^{-\gg}$ with $\gg>1,$ we find a constant $c>0$ such that
\beq\label{U02'} \beg{split} &\lim_{t\to\infty} t^q \E^\mu[\W_2(\mu_t^B, \tt\mu_{t,r_t}^B)^{2q}]\\
&\le 2^q  \lim_{t\to\infty} t^q \E^\mu[\W_2(\mu_{t,r_t}^B, \tt\mu_{t,r_t}^B)^{2q}+ \W_2(\mu_{t,r_t}^B, \mu_t^B)^{2q}]\\
&\le c\lim_{t\to\infty} (r_tt)^q\big(1+\E^\mu[\Xi_{r_t}^B(t)^q]\big)=0.\end{split}\end{equation}
This together with \eqref{U01}  and \eqref{UIB} implies \eqref{U02}.
\end{proof}

Next, we consider arbitrary initial distribution $\nu\in \scr P$. Let
$$\nu_\vv:= \nu P_\vv^B,\ \ \vv\in (0,1).$$
By \eqref{A10}, there exists a constant $c>0$ such that
\beq\label{BNU*} \nu_\vv\le c\vv^{-\ff d {2\aa}}\mu,\ \ \E^{\nu_\vv} \le c \vv^{-\ff d {2\aa}}\E^\mu,\ \ \
\vv\in (0,1).\end{equation}
Let
\beq\label{BME} \mu_t^{B,\vv}:= \ff 1 t\int_\vv^{t+\vv}\dd_{X_s^B}\d s,\ \ t,\vv>0.\end{equation}
By the Markov property, $\mu_t^{B,\vv}$ is the empirical measure with initial distribution $\nu_\vv$, so that
for any nonnegative measurable function $F$ on $\scr P$,
\beq\label{BG} \E^\nu [F(\mu_t^{B,\vv})]= \E^{\nu_\vv} [F(\mu_t^B)],\ \ t,\vv>0.\end{equation}
To estimate $\W_2(\mu_t^{B,\vv},\mu)$, we take
\beq\label{XIBE}\Xi^{B,\vv}(t):= \sum_{i=1}^\infty \ff 1 {\ll_i} \psi_i^{B,\vv}(t)^2,\ \ \
\psi_i^{B,\vv}(t):= \ff 1 {\ss t}\int_\vv^{t+\vv} \phi_i(X_s^B)\d s.\end{equation}

\beg{prp}\label{P2} Assume $(A_1)$ and $(A_2)$.
\beg{enumerate} \item[$(1)$] If $\aa>0$ such that $d+d'<2+4\aa$, then for any   $q\in [1, \ff d{(d+d'-2-2\aa)^+})$,
\beq\label{BNM} \lim_{\vv\downarrow 0}  \sup_{t\ge 1,\nu\in \scr P} \E^\nu\big[|t\W_2(\mu_t^B,\mu)^2
- t\W_2(\mu_t^{B,\vv},\mu)^2|^{2q}\big]=0.\end{equation}
\item[$(2)$] If $\aa>0$ and $d+d'\ge 2+4\aa$, then for any $k>\ff{d+d'-2-2\aa}{2\aa}$ and $q\in [1, \ff d{(d+d'-2-2\aa)^+}),$
\beq\label{BNM'} \lim_{\vv\downarrow 0}   \sup_{t\ge 1,\nu\in \scr P_{k,R}} \E^\nu\big[|t\W_2(\mu_t^B,\mu)^2
- t\W_2(\mu_t^{B,\vv},\mu)^2|^{2q}\big]=0,\ \ R\in [1,\infty).\end{equation}
\item[$(3)$] For any $q\in [1,\infty)$, $k=\infty$ or $k\in (\ff d {2\aa \i(q)},\infty]\cap [1,\infty]$,    there exists a constant $c>0$ such that
\beq\label{BN2}  \sup_{\nu\in \scr P_{k,R}}\E^\nu\big[|\psi_i^{B,\vv}(t)^2-\psi_i^B(t)^2|^q\big]\le c R
\vv^{\ff q 2} t^{-\ff q 2} \ll_i^{\ff{d(q-1)}2-q\aa},\ \ i\in\mathbb N, t\ge 1, \vv\in (0,1).\end{equation}
Moreover, if $\i(q)>\ff d {2\aa}$, then  there exists a constant $c>0$ such that
\beq\label{BN2'}  \sup_{\nu\in \scr P}\E^\nu\big[|\psi_i^{B,\vv}(t)^2-\psi_i^B(t)^2|^{q}\big]\le c \vv^{\ff q 2} t^{-\ff q 2} \ll_i^{\ff{d(q-1)}2-q\aa},\ \ i\in\mathbb N, t\ge 1, \vv\in (0,1).\end{equation}
\end{enumerate}
\end{prp}

\beg{proof} (1) By $d+d'<2+4\aa$, we have
$\ff d {2\aa}< \ff d{(d+d'-2-2\aa)^+}.$ So, it suffices to prove for $1\le q\in (\ff d {2\aa}, \ff d{(d+d'-2-2\aa)^+}).$ It is easy to see that
$$\pi_{t,\vv}:=\ff 1 t \int_\vv^t \dd_{(X_s^B,X_s^B)}\d s+\ff 1 t \int_0^\vv\dd_{(X_s^B,X_{t+s}^B)}\d s\in \C(\mu_t^B,\mu_t^{B,\vv}),\ \ t>\vv\ge 0.$$ So,
  \eqref{RD} implies that for any $t>\vv\ge 0$ and $ p\in [1,\infty),$
\beq\label{BWE} \W_p(\mu_t^B,\mu_t^{B,\vv})^p\le\int_{\R^d\times \R^d}|x-y|^p\pi_{t,\vv}(\d x,\d y)= \ff 1 t \int_0^\vv \rr(X_s^B,X_{s+t}^B)^p\d s\le \ff{\vv D^p}t.\end{equation}
On the other hand, by  \eqref{Ledoux}, \eqref{UIB}, \eqref{BNU*} and \eqref{BG}, we  find a map
$$c: \ \Big[1,\ff d{(d+d'-2-2\aa)^+}\Big)\to (0,\infty)$$ such that
\beq\label{U*} \beg{split} &\sup_{t\ge 1,\nu\in\scr P} t^q\E^\nu[\W_2(\mu_t^{B,\vv},\mu)^{2q}]= \sup_{t\ge 1, r\in (0,1],\nu\in\scr P} t^q \E^{\nu_\vv}[\W_2(\mu_{t,r}^{B},\mu)^{2q}]\\
&\le 4^q \sup_{t\ge 1, r\in (0,1],\nu\in\scr P}  \E^{\nu_\vv}[\Xi_{r}^B(t)^q]\le c(q)\vv^{-\ff d{2\aa}},\ \ \vv\in (0,1),\ q\in \Big[1,\ff d{(d+d'-2-2\aa)^+}\Big).\end{split}\end{equation}
Combining this with \eqref{BWE}, \eqref{BNU*} and $q>\ff{d}{2\aa}$,  we find constants $c_1,c_2 >0$ such that
\beg{align*} &\lim_{\vv\downarrow 0}  \sup_{t\ge 1,\nu\in \scr P} \E^\nu\big[|t\W_2(\mu_t^B,\mu)^2
-t \W_2(\mu_t^{B,\vv},\mu)^2|^{2q}\big]\\
&\le \lim_{\vv\downarrow 0} \sup_{t\ge 1,\nu\in \scr P} \E^\nu\big[|t\W_2(\mu_t^B,\mu_t^{B,\vv})^2 + 2 t\W_2(\mu_t^B,\mu_t^{B,\vv})\W_2(\mu_t^{B,\vv},\mu)|^{2q}\big]\\
&\le \lim_{\vv\downarrow 0}\sup_{t\ge 1,\nu\in \scr P} c_1\Big( \vv^{2q} +\vv^{q} \sup_{\nu\in \scr P} \E^{\nu_\vv}[\W_2(\mu_t^B,\mu)^{2q}]\Big)
\le \lim_{\vv\downarrow 0}   c_2\vv^{q-\ff d{2\aa}}=0.\end{align*}
So, \eqref{BNM} holds.

(2) Let  $\aa>0$, $d+d'\ge 2+4\aa$ and $k>\ff{d+d'-2-2\aa}{2\aa}$. It suffices to prove \eqref{BNM'} for
$1\le q\in (\ff d{2\aa k}, \ff d{(d+d'-2-2\aa)^+}).$
By the same reason leading to \eqref{BNU*}, we find a constant $c>0$ such that
$$\sup_{\nu\in \scr P_{k,R}}\E^{\nu_\vv}\le c \vv^{-\ff d{2\aa k}}\E^\mu,\ \ \vv\in (0,1).$$
Hence, as shown above that $q>\ff d{2\aa k}$ implies
$$ \lim_{\vv\downarrow 0}  \sup_{t\ge 1,\nu\in \scr P_{k,R}} \E^\nu\big[|t\W_2(\mu_t^B,\mu)^2
-t \W_2(\mu_t^{B,\vv},\mu)^2|^{2q}\big]\le \lim_{\vv\downarrow 0}   c_2\vv^{q-\ff d{2\aa k}}=0. $$

(3) Let $\psi_i^B$ and $\psi_i^{B,\vv}$ be in \eqref{XIB} and \eqref{XIBE}. We have
\beq\label{SPP1} \beg{split} &\big|\psi_i^B(t)^2-\psi_i^{B,\vv}(t)^2\big|\\
&=\ff 1 t \bigg|\int_0^\vv \big\{\phi_i(X_{t+s}^B)-\phi_i(X_s^B)\big\}\d s\bigg|\cdot
\bigg|\int_0^t \big\{\phi_i(X^B_{s+\vv})+\phi_i(X_s^B)\big\}\d s\bigg|\\
&\le \ff{\ss \vv}{\ss t} \big(|\psi_i^{B,t}(\vv)|+|\psi_i^B(\vv)|\big)
\big(|\psi_i^{B,\vv}(t)|+|\psi_i^B(t)|\big).\end{split}\end{equation}
Since    \eqref{BG} implies  $\{\nu_t: \nu\in \scr P, t\ge 1\}\subset \scr P_{\infty,R}$
for some $R>0$,  by \eqref{BLN10} and \eqref{BG} we find a constant $c_1>0$ such that
\beq\label{SPP2}\sup_{\vv>0,t\ge 1,\nu\in \scr P} \E^\nu[\psi_i^{B,t}(\vv)^{2q}]=\sup_{\vv>0,t\ge 1,\nu\in \scr P} \E^{\nu_t}[\psi_i^{B}(\vv)^{2q}]
\le c_1 \ll_i^{\ff{d(q-1)}2 -q\aa}.\end{equation}
Moreover, by \eqref{BLN10} for $k=\infty$, \eqref{BLN1} for $\aa>0$ and $k\in (\ff d {2\aa\i(q)},\infty]\cap [1,\infty]$,  and the fact that $\nu\in\scr P_{k,R}$ implies
$\nu_\vv\in \scr P_{k,R}$ for $\vv>0$,   we find a constant $c_2>0$ such that
\beq\label{NNV}\sup_{\vv>0,t\ge 1,\nu\in \scr P_{k,R}}
 \E^\nu[\psi_i^{B,\vv}(t)^{2q}+\psi_i^B(\vv)^{2q}]\le c_2 \ll_i^{\ff{d(q-1)}2 -q\aa},\ \ i\ge 1.\end{equation}
Combining these estimates we derive \eqref{BN2}.

Finally, when $\i(q)>\ff d {2\aa}$, \eqref{BN2'} can be proved in the same way by using
\eqref{BLN1*} in place of \eqref{BLN1}.
\end{proof}

We are now ready to prove Theorem \ref{TU1}.

\beg{proof}[Proof of Theorem \ref{TU1}]

(1) It suffices to prove for $q\in (\ff d{2\aa},q_\aa).$   By $\aa>\aa(d,d'),$ we have
$$\ff d {2\aa}<\ff{2(1+\aa)-d'}{(d+d'-2-2\aa)^+}.$$
So, either $\i(q)>\ff d {2\aa}$ or $\i(q)<\ff{2+2\aa-d'}{(d+d'-2-2\aa)^+}.$ Below we consider these two situations respectively.

$(1_a)$ Let $\i(q)>\ff d {2\aa}.$
 By \eqref{BG}, \eqref{BNU*}  and \eqref{U01}, we obtain
\beq\label{Y02'}\beg{split}& \lim_{t\to \infty }\sup_{\nu\in\scr P} \E^\nu\Big[
\Big|\big\{t\W_2(\mu_t^{B,\vv},\mu)^2- \Xi^{B,\vv}(t)\big\}^+\Big|^q\Big]\\
&=\lim_{t\to \infty }\sup_{\nu\in\scr P} \E^{\nu_\vv}\Big[
\Big|\big\{t\W_2(\mu_t^{B},\mu)^2- \Xi^{B}(t)\big\}^+\Big|^q\Big]\\
&\le \lim_{t\to \infty }c\vv^{-\ff d {2\aa}} \E^{\mu}\Big[
\Big|\big\{t\W_2(\mu_t^{B},\mu)^2- \Xi^{B}(t)\big\}^+\Big|^q\Big]=0,\ \ \vv\in (0,1).\end{split}\end{equation}
Next, by \eqref{XIB} and \eqref{XIBE},
$$|\Xi^B(t)-\Xi^{B,\vv}(t)|\le \sum_{i=1}^\infty \ff 1 {\ll_i} |\psi_i^{B,\vv}(t)^2-\psi_i^B(t)^2|.$$
Combining this with \eqref{BN2'}, when $\i(q)>\ff d {2\aa}$ we find a constant $k_1>0$ such that
\beg{align*} &\sup_{t\ge 1,\nu\in \scr P}  t^{\ff q 2} \E^\nu\big[ |\Xi^B(t)-\Xi^{B,\vv}(t)|^q\big]\\
&\le \Big(\sum_{i=1}^\infty \ll_i^{-\theta}\Big)^{q-1} \sum_{i=1}^\infty
\ll_i^{\theta(q-1) -q}\sup_{t\ge 1,\nu\in \scr P}  t^{\ff q 2} \E^\nu\big[|\psi_i^{B,\vv}(t)^2-\psi_i^B(t)^2|^q\big]
\\
&\le k_1 \vv^{\ff q 2} \Big(\sum_{i=1}^\infty \ll_i^{-\theta}\Big)^{q-1}
 \sum_{i=1}^\infty\ll_i^{\theta(q-1)+\ff{d(q-1)}2 -q-\aa q },\ \ \theta\in\R,\  \vv\in (0,1).\end{align*}
 Taking
 \beq\label{THA} \theta= 1+\aa-\ff{d(q-1)}{2q}\end{equation}
 such that $-\theta=\theta(q-1)+\ff{d(q-1)}2 -q-\aa q,$ we arrived at
 \beq\label{TA}\sup_{t\ge 1,\vv\in (0,1),\nu\in \scr P} \vv^{-\ff q 2} t^{\ff q 2} \E^\nu\big[ |\Xi^B(t)-\Xi^{B,\vv}(t)|^q\big]\le k_1\Big(\sum_{i=1}^\infty\ll_i^{-\theta}\Big)^q.\end{equation}
  Noting that \eqref{QA},
    \eqref{THA} and $q<q_\aa$ imply
    $\ff{2\theta}{d'} >1$,   by \eqref{A1-1} we find a  constant  $k_2 >0$ such that
 $$\sum_{i=1}^\infty\ll_i^{-\theta}\le k_2\sum_{i=1}^\infty i^{-\ff{2\theta}{d'}}<\infty.$$
 Thus,
 \beq\label{N1-0}\sup_{t\ge 1,\vv\in (0,1),\nu\in \scr P} \vv^{-\ff q 2} t^{\ff q 2} \E^\nu\big[ |\Xi^B(t)-\Xi^{B,\vv}(t)|^q\big]<\infty.\end{equation}
 Consequently,
 $$\lim_{t\to\infty} \sup_{\vv\in (0,1),\nu\in \scr P}  \E^\nu\big[ |\Xi^B(t)-\Xi^{B,\vv}(t)|^q\big]=0.$$
 Combining this with \eqref{BNM} and \eqref{Y02'}, we derive \eqref{BU01}.

 $(1_b)$ Let $\i(q)<\ff{2+2\aa-d'}{(d+d'-2-2\aa)^+}.$ Since  $q>\ff d{2\aa}$ implies $\i(q)+1>\ff d {2\aa}$, by \eqref{BLN1*} and H\"older's inequality, we find a constant
 $c_1>0$ such that
 \beg{align*}\sup_{t\ge 1,\nu\in \scr P} \E^\nu\big[|\psi_i^B(t)|^{2q}\big]
 \le \sup_{t\ge 1,\nu\in \scr P} \big(\E^\nu\big[|\psi_i^B(t)|^{2\{\i(q)+1\}}\big]\big)^{\ff q{\i(q)+1}}
 \le c_1 \ll_i^{\ff{dq\i(q)}{2(\i(q)+1)}-q\aa},\ \ \ i\in\mathbb N.\end{align*}
 By the  calculations leading to \eqref{TA}, we derive the same estimate for
 $$\theta:=1+\aa-\ff{d\i(q)}{2(\i(q)+1)}.$$
We have  $\theta>\ff{d'}2$ due to $\i(q)<\ff{2+2\aa-d'}{(d+d'-2-2\aa)^+}.$
 Hence, \eqref{N1-0} holds, which together  with \eqref{BNM} and \eqref{Y02'} imply \eqref{BU01}.

 (2) Let $\aa\in (0,1], q\in [1,q_\aa)$ and $k\in (\ff d{2\aa\i(q)},\infty]\cap [1,\infty].$
 By using \eqref{BN2} in place of \eqref{BN2'}, the proof of \eqref{N1-0} implies
 $$\sup_{t\ge 1,\vv\in (0,1),\nu\in \scr P_{k,R}} \vv^{-\ff {q} 2} t^{\ff {q} 2} \E^\nu\big[ |\Xi^B(t)-\Xi^{B,\vv}(t)|^{q}\big]<\infty,\ \ R\in [1,\infty),$$
 so that
 $$\lim_{t\to\infty} \sup_{\vv\in (0,1),\nu\in \scr P_{k,R}}  \E^\nu\big[ |\Xi^B(t)-\Xi^{B,\vv}(t)|^q\big]=0,\ \ R\in [1,\infty).$$
 This together with   \eqref{BNM} and \eqref{Y02'} implies  \eqref{BU02}.

 When $k=\infty,$ \eqref{BU02} follows from \eqref{U02} and that
 $\E^\nu\le \|h\|_\infty\E^\mu$ for $\nu= h\mu.$

\end{proof}

\subsection{Proof of Theorem \ref{TU2}}
(1) Let $d'<2(1+\aa)$. By Lemma \ref{L0},    \eqref{BH1} holds with $\eta_Z^B<\infty.$
Combining this with \eqref{U02} and $\E^\nu\le \|h\|_\infty\E^\mu$ for $\nu= h\mu$,
we obtain
\beq\label{BCB}  \limsup_{t\to\infty} \sup_{\nu\in \scr P_{\infty,R}} t\E^\nu\big[\W_2(\mu_t^B,\mu)^2\big]
 \le \limsup_{t\to\infty}\sup_{\nu\in \scr P_{\infty,R}} \E^\nu \big[\Xi^B(t)\big]\le \eta_Z^B,\ \ R\in [1,\infty).  \end{equation}
 On the other hand, by \eqref{BNU*}, for any $\vv\in (0,1)$ there exists $R\in [1,\infty)$
 such that $\nu_\vv\in \scr P_{\infty,R}$ holds for all $\nu\in \scr P$. So,
 \eqref{BCB} together with \eqref{BG} implies
 \beq\label{*1}  \limsup_{t\to\infty} \sup_{\nu\in \scr P} t\E^\nu\big[\W_2(\mu_t^{B,\vv},\mu)^2\big]=\limsup_{t\to\infty} \sup_{\nu\in \scr P} t\E^{\nu_\vv}\big[\W_2(\mu_t^{B},\mu)^2\big]\le \eta_Z^B,\ \ \vv\in (0,1).\end{equation}
 Combining this with \eqref{BWE} for $p=2$ and using the triangle inequality,
 we find a constant $c>0$ such that
 \beg{align*} &\limsup_{t\to\infty} \sup_{\nu\in \scr P} t\E^\nu\big[\W_2(\mu_t^{B},\mu)^2\big]\\
 &\le \limsup_{t\to\infty} \sup_{\nu\in \scr P} t\E^\nu\big[(1+\vv^{\ff 1 2})\W_2(\mu_t^{B,\vv},\mu)^2+ (1+\vv^{-\ff 1 2}) \W_2(\mu_t^{B,\vv},\mu^B_t)^2\big]\\
 &\le (1+\vv^{\ff 1 2})\eta_Z^B+ c (\vv+\vv^{\ff 1 2}),\ \ \ \vv\in (0,1).\end{align*}
 Letting $\vv\downarrow 0$ we obtain \eqref{BU02'}.

 (2) Let $d'\ge 2(1+\aa)$. By \eqref{BXR}, \eqref{Ledoux'} and \eqref{LY1} with $a_0=0, q=1$ and $\bb=-\ff 1 2$, we find a constant $c_1>0$ such that
 \beq\label{LY3} \beg{split} &\E^\mu[\W_2(\mu_{t,r}^B,\mu)^2]\le  4  \E^\mu[ \mu(|\nn \hat L^{-1}(f_{t,r}^B-1)|^2)]
 =  4  \E^\mu[\mu(|(-\hat L)^{-\ff 1 2} (f_{t,r}^B-1)|^2)]\\
 &\le \ff {c_1} t \big(r^{1+\aa-\ff{d'}2}+1_{\{d'=2+2\aa\}} \log (1+r^{-1})\big),\ \ t\ge 1,r\in (0,1].\end{split}\end{equation}
 Combining this with \eqref{XB1} and the triangle inequality, we find a constant $c_2>0$ such that
 \beg{align*} &\E^\mu[\W_2(\mu_{t}^B,\mu)^2]\le 2 \E^\mu[\W_2(\mu_{t,r}^B,\mu)^2]+ 2 \E^\mu[\W_2(\mu_{t,r}^B,\mu_t^B)^2]\\
 &\le \ff {c_1} t \big\{r^{1+\aa-\ff{d'}2}+1_{\{d'=2+2\aa\}} \log (1+r^{-1})\big\}
 +c_2 r,\ \ t\ge 1, r\in (0,1].\end{align*}
 Taking $r=t^{-\ff 2{d'-2\aa}}$ when $d'>2+2\aa$, and $r=t^{-1}$ when $d'=2+2\aa$, we
 find a constant $c_3>0$ such that
  $$ \E^\mu[\W_2(\mu_{t}^B,\mu)^2]\le
  \beg{cases} c_3t^{-1}\log (1+t),\ &\text{if}\ d'=2(1+\aa),\\
 c_3 t^{-\ff{2}{d'-2\aa}},\ &\text{if}\ d'>2(1+\aa).\end{cases} $$
 By combining this with \eqref{BNU*} and \eqref{BG}, we find a constant $c_4>0$ such that
\beq\label{BWE'}\beg{split}& \sup_{\nu\in \scr P} \E^\nu[\W_2(\mu_{t}^{B,1},\mu)^2] =
  \sup_{\nu\in \scr P} \E^{\nu_1}[\W_2(\mu_{t}^{B},\mu)^2]
  \le c_4 \E^{\mu}[\W_2(\mu_{t}^{B},\mu)^2]\\
  &\le c_3c_4 \Big\{1_{\{d'=2(1+\aa)\}}  t^{-1}\log (1+t)+t^{-\ff{2}{d'-2\aa}}\Big\},\ \ t\ge 1. \end{split}
  \end{equation}
Noting that the triangle inequality  implies $$\E^\nu[\W_2(\mu_{t}^{B},\mu)^2]\le
 2\W_2(\mu_{t}^{B,1},\mu)^2+ 2\W_2(\mu_{t}^{B,1},\mu_t^B)^2,$$
 we deduce   \eqref{BUa} from \eqref{BWE} and \eqref{BWE'}.

\subsection{Proof of Theorem \ref{TU3}}

    By approximating $\mu_t^B$ using $\mu_t^{B,1}$ as in \eqref{BWE} and \eqref{BWE'}, we only need to prove for $\nu=\mu$.

 By \eqref{Ledoux'}, \eqref{A11} and \eqref{LY1}, for any $\kk\in (0,\ff 1 2)$ we find constants $c_1,c_2>0$ such that for any $t\ge 1$ and $r\in (0,1],$
 \beq\label{LY5} \beg{split} &\big(\E^\mu[\W_{2p}(\mu_{t,r}^B,\mu)^{2q}]\big)^{\ff 1 q}
 \le c_1 \big(\E^\mu[\|(-\hat L)^{\ff{d(p-1)}{4p}-\kk}(f_{t,r}^B-1)\|_2^{2q}\big)^{\ff 1 q}\\
 &\le \ff{c_2}t \Big\{r^{-\big(\ff{d(p-1)}{2p}-2\kk+\ff{d'}2+\ff{d(q-1)}{2q} -\aa\big)^+}+1_{\big\{\ff{d(p-1)}{2p}-2\kk+\ff{d'}2+\ff{d(q-1)}{2q} -\aa=0\big\}}\log(1+r^{-1})\Big\}.\end{split}\end{equation}
 Below we prove assertions (1)-(3) in Theorem \ref{TU3} respectively.

(1) Let $\gg_{\aa,p,q}<0$. We may take $\kk\in (0,\ff 1 2)$ such that
 $$\ff{d(p-1)}{2p}-2\kk+\ff{d'}2+\ff{d(q-1)}{2q} -\aa<0,$$
  so that \eqref{LY5} implies
  $$\big(\E^\mu[\W_{2p}(\mu_{t,r}^B,\mu)^{2q}]\big)^{\ff 1 q}\le \ff{c_2} t,\ \ r\in (0,1], t\ge 1.$$
  By Fatou's lemma for   $r\to 0$, we obtain \eqref{BUb} for $\nu=\mu$.

  (2) Let   $\gg_{\aa,p,q}\ge 0$.  For any $\gg>\gg_{\aa,p,q},$ we find
  $\kk\in (0,\ff 1 2)$ such that
  $$ \ff{d(p-1)}{2p}-2\kk+\ff{d'}2+\ff{d(q-1)}{2q} -\aa \le\gg,$$
  so that \eqref{LY5}, \eqref{XB1}  and the triangle inequality imply
  $$\big(\E^\mu[\W_{2p}(\mu_{t}^B,\mu)^{2q}]\big)^{\ff 1 q}\le
  c_1\big\{t^{-1}r^{-\gg}+ r\big\},\ \ r\in (0,1], t\ge 1$$ for some constant $c_1>0$.
  Taking $r=t^{-\ff 1 {1+\gg}}$ we obtain \eqref{BUd} for $\nu=\mu.$

   (3) Let   \eqref{RS} hold and $\gg_{\aa,p,q}\ge 0$. By \eqref{Ledoux'} we find constants $c_1,c_2>0$ such that
  \beq\label{NE1} \beg{split} &\W_{2p}(\mu_{t,r}^B,\mu)^{2p} \le c_1 \|\nn (-\hat L)^{-1}(f_{t,r}^B-1)\|_{2p}^{2p}
  \le c_2 \|(a_0-\hat L)^{\ff 1 2}  (-\hat L)^{-1}(f_{t,r}^B-1)\|_{2p}^{2p}.\end{split}\end{equation}
On the other hand,
by the Sobolev embedding theorem, \eqref{A10} implies that for any constants $k_2\ge k_1>-\infty$ and
$q_1\ge q_2\ge 1$ with
\beq\label{QK} \ff 1 {q_1}=\ff 1 {q_2} +\ff{k_1-k_2}d,\end{equation}
there exists a constant $C>0$ such that
$$ \|(-\hat L)^{\ff {k_1}2}f\|_{q_1}\le C\|(-\hat L)^{\ff{k_2} 2}f\|_{q_2},\ \ \mu(f)=0.$$
Taking
$$k_1=-2,\ \  k_2=\ff{d(p-1)}{2p}-2,\ \  q_1=2p,\ \  q_2=2$$ such that \eqref{QK} holds,
we find a constant $c_2>0$  such that
$$ \|(a_0-\hat L)^{\ff 1 2}  (-\hat L)^{-1}(f_{t,r}^B-1)\|_{2p}\le c_2  \|(a_0-\hat L)^{\ff 1 2}(-\hat L)^{\ff{d(p-1)}{4p}-1} (f_{t,r}^B-1)\|_2.$$
Combining this with \eqref{LY1} and \eqref{NE1}, we find a constant  $c_3>0$ such that
  $$\big(\E^\mu[\W_{2p}(\mu_{t,r}^B,\mu)^{2q}]\big)^{\ff 1 q}\le
  \ff{c_3}t \big\{r^{-\gg_{\aa,p,q}}+ 1_{\{\gg_{\aa,p,q}=0\}}\log(1+r^{-1})\big\},
  \ \ t\ge 1, r\in (0,1].$$
 By this together with \eqref{XB1}  and the triangle inequality, we find a constant $c_4>0$
  such that
  $$\big(\E^\mu[\W_{2p}(\mu_{t}^B,\mu)^{2q}]\big)^{\ff 1 q}
  \le c_4 \Big\{t^{-1}r^{-\gg_{\aa,p,q}}+ t^{-1} 1_{\{\gg_{\aa,p,q}=0\}}\log(1+r^{-1})+ r\big\},\ \ t\ge 1,r\in (0,1].$$
  Taking $r=t^{-\ff 1 {1+\gg_{\aa,p,q}}},$ we obtain \eqref{BUd}.

  \section{Proofs of Theorems \ref{TL1} and \ref{TL2}}

We will follow the line of \cite{WZ} to estimate  the lower bound of $\W_2(\mu_t^B,\mu)$ by using an idea of \cite{AM}.
For any $f\in \D(L)$ with $\|f\|_\infty+\|\nn f\|_\infty+\|Lf\|_\infty<\infty$, let
$$T_t^\si f:=-\si \log \hat P_{\ff{\si t} 2} \e^{\si^{-1}f},\ \ \si>0,t\in [0,1].$$

\beg{lem}\label{LB2} Assume that $(M,\rr)$ is a geodesic space. If $\eqref{A120}$ and $\eqref{B*}$ hold, then there exist  constants $k_1,k_2>0$ such that
$ \|\nn f\|_\infty^2\le k_1\si$ implies
\beg{equation}\label {B5} \beg{split}&T_1^\si f(y)-f(x)\le \ff 1 2 \rr(x,y)^2+\si\|\hat Lf\|_\infty^2 +k_2  \si^{\theta}\|\nn f\|_\infty^2,\\
&\mu(f-T_1^\si f)\le \ff 1 2 \mu(|\nn f|^2)+k_2\si^{-1}\|\nn f\|_\infty^4.\end{split} \end{equation}
\end{lem}

\beg{proof}   Let
  $k_1=2\kk_0$, then $ \|\nn f\|_\infty^2\le k_1\si $ implies
$$ \ff{t\si }2 \|\nn\si^{-1}f\|_\infty^2\le \ff 1 2 \si^{-1}\|\nn f\|_\infty^2\le \kk_0,\ \ t\in [0,1],$$
so that \eqref{21} holds for $m=1$ and $(\ff{\si t} 2,\si^{-1}f)$ replacing $(t,f),$
which together  with
\eqref{A120} yields
\beq\label{*ZZ} \beg{split}&|\nn T_t^\si f|^2 =\ff{\si^2|\nn \HAT P_{\ff{t\si}2} \e^{\si^{-1}f}|^2}{(\HAT P_{\ff{t\si}2} \e^{\si^{-1}f})^2} \le \ff{(1+k(2))\si^2   \HAT P_{\ff{t\si}2} (|\nn \si^{-1}f|^2 \e^{2\si^{-1}f})}{(\HAT P_{\ff{t\si}2} \e^{\si^{-1}f})^2}  \\
 &\le \ff{(1+k(2)) \|\nn f\|_\infty^2   \HAT P_{\ff{t\si}2} (  \e^{2\si^{-1}f})}{(\HAT P_{\ff{t\si}2} \e^{\si^{-1}f})^2}\le 2 (1+k(2)) \|\nn f\|_\infty^2=:c_1\|\nn f\|_\infty^2. \end{split}\end{equation}
 Next, by  \eqref{A21}, we find a constant $c_2>0$ such that
   \beq\label{ASC} \beg{split}&\hat LT_t^\si f=-\ff{\si \hat L\HAT P_{\ff{t\si}2} \e^{\si^{-1}f}}{\HAT P_{\ff{t\si}2}\e^{\si^{-1}f}} +\ff{\si|\nn \HAT P_{\ff{t\si}2}\e^{\si^{-1}f}|^2}{(\HAT P_{\ff{t\si}2}\e^{\si^{-1}f})^2} \\
   &=\ff{-\hat P_{\ff{t\si}2} \{\e^{\si^{-1}f}\hat Lf\}}{\HAT P_{\ff{t\si}2}\e^{\si^{-1}f}}
    +\ff{\si|\nn \HAT P_{\ff{t\si}2}\e^{\si^{-1}f}|^2-\si(\HAT P_{\ff{t\si}2}\e^{\si^{-1}f})\HAT P_{\ff{t\si}2}(|\nn \si^{-1}f|^2\e^{\si^{-1}f})}{(\HAT P_{\ff{t\si}2}\e^{\si^{-1}f})^2}\\
    &\le \|\hat L f\|_\infty + c_2 \si^{\theta-1} \|\nn f\|_\infty^2, \ \ \si,t\in (0,1], \|\nn f\|_\infty^2\le k_1\si. \end{split}\end{equation}
Moreover, for any two points $x,y\in M$,  let $\gg: [0,1]\to M$ be the minimal geodesic from $x$ to $y$ with
$$|\dot\gg_t|:=\limsup_{s\to t} \ff{\rr(\gg_t,\gg_s)}{|t-s|}=\rr(x,y),\ \ \text{a.e.}\ t\in [0,1].$$
 So,
\beq\label{KL1}\limsup_{s\to t} \ff{|f(\gg_t)-f(\gg_s)|}{|t-s|}\le |\nn f(\gg_t)| \rr(x,y),\ \ \ t\in [0,1]. \end{equation}
 By   the backward Kolmogorov equation and the chain rule, we have
\beq\label{KL}\pp_t T_t^\si f= -\ff{\si^2\hat L\HAT P_{\ff{t\si}2 }\e^{\si^{-1}f}}
{2\HAT P_{\ff {t\si}2}\e^{\e^{\si^{-1}f}}}=\ff\si 2 \hat LT_t^\si f -\ff 1 2 |\nn T_t^\si f|^2.\end{equation} This together with  \eqref{ASC} and \eqref{KL1} yields
  \beg{align*}  &\ff{\d }{\d t} T_t^\si f (\gg_t) = (\pp_t T_t^\si f) (\gg_t) +
  \ff{\d }{\d t} T_s^\si f (\gg_t)\Big|_{s=t}\\
 &\le \ff \si 2 \hat LT_t^\si f(\gg_t)-\ff 1 2 |\nn T_t^\si f(\gg_t)|^2 +|\nn T_t^\si f(\gg_t)|\rr(x,y)  \\
 &\le  \ff 1 2 \big[\rr(x,y)^2+ \si  \|\hat Lf\|_\infty + c_2  \si^\theta  \|\nn f\|_\infty^2\big],\ \      \ t\in [0,1],\ \|\nn f\|_\infty^2\le k_1\si.  \end{align*}
 Integrating over $t\in [0,1]$ and noting that $T_0^\si f=f$, we derive   the first inequality in \eqref{B5}.

 On the other hand,  by \eqref{KL}, $T_0^\si f=f$ and  $\mu(\hat LT_t^\si f)=0$, we obtain
 \beq\label{PD1}\beg{split}& \mu(f-T_1^\si f)=-\int_M\d\mu\int_0^1(\pp_t T_t^\si f)\d t\\
 &=\int_0^1 \d t \int_M\Big\{\ff 1 2 |\nn T_t^\si f|^2 -\ff\si 2\hat LT_t^\si f\Big\}\d\mu =\ff 1 2 \int_0^1 \mu(|\nn T_t^\si f|^2) \d t.\end{split} \end{equation}
 Moreover, by \eqref{KL} and
 the integration by parts formula, we obtain
 \beg{align*}&\ff{\d}{\d s} \mu(|\nn T_s^\si f|^2) =-\ff{\d}{\d s} \int_M (T_s^\si f)\hat L T_s^\si f \d\mu \\
 &=-\int_M (\hat LT_s^\si f)\pp_sT_s^\si f\d\mu-\int_M (T_s^\si f)\hat L(\pp_sT_s^\si f)\d\mu\\
 &=-2 \int_M (\hat LT_s^\si f)\pp_s T_s^\si f\,\d\mu= - 2 \int_M (\hat LT_s^\si f)\Big(\ff \si 2\hat LT_s^\si f -\ff 1 2 |\nn T_s^\si f|^2\Big)\d\mu\\
 &\le \ff 1 {4\si}\|\nn T_s^\si f\|_\infty^4,\ \ s\in (0,1], t\in [0,1].  \end{align*}
  This together with \eqref{*ZZ}   implies
 \beq\label{*X3}\beg{split} &\mu(|\nn T_t^\si f|^2)-\mu(|\nn f|^2) \le \ff{c_1^2}{4 \si}  \|\nn f\|_\infty^4,\ \ t\in [0,1], \si \in (0,1],\|\nn f\|_\infty^2\le k_1\si \end{split}\end{equation}
    Substituting this into \eqref{PD1}, we derive the second estimate in \eqref{B5}.

\end{proof}

\beg{proof}[Proof of Theorem \ref{TL1}]

Similarly to the proof of Theorem \ref{TU1} using Proposition  \ref{P1} and the approximation argument with Proposition  \ref{P2},  the assertions Theorem \ref{TL1}(1) and (2)  follow from
\beq\label{ABC}  \lim_{t\to\infty}\E^\mu \big[\big|\{t h(0) \W_2(\mu_{t}^B,\mu)^2- \Xi^B(t)\}^-\big|^q\big]=0,\ \ q\in [1,q_\aa).\end{equation}
Moreover, according to the proof of Theorem \ref{TU2}(1), Theorem \ref{TL1}(3) is implied by Theorem \ref{TL1}(1) and (2). So, it remains to verify \eqref{ABC}.
The main idea for the proof of \eqref{ABC} goes back to
\cite{AM,WZ}, but we have to make suitable modifications for the present  situation.
Let $\hat f_{t,r}=(-\hat L)^{-1}(1-f_{t,r}^B)$.

Firstly, by \eqref{UAA'}, we find a constant $c_1>0$ such that
\beq\label{FINF}\|\hat f_{t,r}\|_\infty \le  \int_0^\infty \|\hat P_s(f_{t,r}^B-1)\|_\infty\d s \le c_1\|f_{t,r}^B-1\|_\infty \int_0^\infty \e^{-\ll_1 s}\d s =\ff{c_1}{\ll_1} \|f_{t,r}^B-1\|_\infty.\end{equation}
By \eqref{A12} and \eqref{UAA}, we find a constant $c_2>0$ such that
$$\|\nn \hat P_s g\|_\infty\le c_2 (s\land 1)^{-\ff 1 2}\e^{-\ll_1 s}\|g\|_\infty,\ \ s>0,
g\in \B_b(M).$$
Noting that $\hat f_{t,r}:=  (-\hat L)^{-1}(f_{t,r}^B-1)=\int_0^\infty \hat P_s (f_{t,r}^B-1)\d s,$ we obtain
$$\|\nn \hat f_{t,r}\|_\infty \le \int_0^\infty \|\nn \hat P_s(f_{t,r}^B-1)\|_\infty\d s
\le c_2\|f_{t,r}^B-1\|_\infty \int_0^\infty (1\land s)^{-\ff 1 2}\e^{-\ll_1 s}\d s.$$
Combining this  with \eqref{FINF} and  $|\hat L \hat f_{t,r}|=|f_{t,r}^B-1|$, we find a constant $c>0$ such that
\beq\label{LN1}\|\hat L \hat f_{t,r}\|_\infty+ \|\hat f_{t,r} \|_\infty+\|\nn \hat f_{t,r}\|_\infty \le c  \|f_{t,r}^B-1\|_\infty,\ \ t,r>0.\end{equation}

Next, let
 \beq\label{ANN} \beg{split} &C_1(\hat f_{t,r},\si):=  \si \|\hat L \hat f_{t,r}\|_\infty^2+k_2 \si^{\theta} \|\nn \hat f_{t,r}\|_\infty^2,\\
  &C_2(\hat f_{t,r},\si):=   k_2 \si^{-1}  \|\nn \hat f_{t,r}\|_\infty^4,\\
  &B_{t,r}(\si):=\{\|f_{t,r}^B-1\|_\infty^2\le k_1 c^{-1}\si^{1+\theta} \},\ \ \si\in (0,1],t,r>0. \end{split}\end{equation}
By \eqref{B5}, the integration by parts formula, and the Kantorovich dual formula, we obtain
  \beg{align*} & C_1(\hat f_{t,r},\si)+\ff 1 2\W_2(\mu_{t,r}^B,\mu)^2\ge  \mu(T_1^\si \hat f_{t,r})-\mu_{t,r}^B(\hat f_{t,r})\\
&= \mu\big((f_{t,r}^B-1) (-\hat L)^{-1} (f_{t,r}^B-1)\big)-\mu(\hat f_{t,r})+\mu(T_1^\si \hat f_{t,r})\\
&\ge \mu(|\nn \hat f_{t,r}|^2)-\ff 1 2 \mu(|\nn \hat f_{t,r}|^2)-C_2(\hat f_{t,r},\si)=\ff 1{2t}\Xi_r^B(t) -C_2(\hat f_{t,r},\si),
\ \ \si\in (0,1],t,r>0.\end{align*}
This together with \eqref{LN1} and \eqref{ANN}   yields
\beq\label{BL3}  1_{B_{t,r}(\si)}\big\{t\W_2(\mu_{t,r}^B,\mu)^2- \Xi_r^B(t)\big\}\ge - c_3t\si^{1+2\theta},\ \ \si\in (0,1], t,r>0   \end{equation}  for some constant $c_3>0$.
On the other hand, by \eqref{A10} we find a constant  $c_4>0$  such that
$$\|f_{t,r}^B-1\|_\infty^2= \|\hat P_{\ff r 2 }(f_{t,\ff r 2}^B-1)\|_\infty^2
\le c_4r^{-\ff d 2} \|f_{t,\ff r 2}^B-1\|_2^2,\ \ r\in (0,1],$$
so that by \eqref{XB4'} with $1$ replacing $\si\in (0,1)$, there exists a constant $c_5>0$ such that
\beg{align*} \E^\mu [\|f_{t,r}^B-1\|_\infty^2]\le c_4 r^{-\ff d 2} t^{-1} \sum_{i=1}^\infty \e^{-\ll_i r}\E^\mu[\psi_i^B(t)^2]  
 \le c_5r^{-\ff d 2-1} t^{-1},\ \ \ t\ge 1, r\in (0,1].\end{align*}
Hence, we find   a constant $c_6>0$ such that
\beq\label{PR} \P^\mu(B_{t,r}(\si)^c)=\P^\mu(\|f_{t,r}^B-1\|_\infty^2> k_1c^{-1}\si^{1+\theta})\le c_6 r^{-\ff d 2 -1}t^{-1}\si^{-(1+\theta)},\ \ t\ge 1, \si, r\in (0,1].\end{equation}
Taking $\si=\si_t:=t^{-\ff 1 {1+3\theta/2}}$ in \eqref{BL3} and \eqref{PR}, we arrive at
$$ \limsup_{t\to\infty}\P^\mu\big(\{t\W_2(\mu_{t,r}^B,\mu)^2- \Xi_r^B(t)\}^-\ge\vv\big)
\le \limsup_{t\to\infty}\P^\mu(B_{t,r}(\si_t)^c)=0,\ \ \vv>0, r\in (0,1].$$
On the other hand, by \eqref{XIB}, \eqref{BXR},   \eqref{BLN10} and $\sum_{i=1}^\infty \ll_i^{-1-\aa}<\infty$ due to \eqref{A1-1} and $d'<2(1+\aa)$, we obtain
\beg{align*}&\lim_{r\to 0} \sup_{t\ge 1} \E^{\mu} \big[ |\Xi_r^B(t)-\Xi^B(t)|\big]\le
\lim_{r\to 0}\sum_{i=1}^\infty \ff{1-\e^{-2\ll _i r}}{\ll_i}\sup_{t\ge 1} \E^\mu [\psi_i^B(t)^2]\\
&\le c(1)\lim_{r\to 0}
\sum_{i=1}^\infty \big(1-\e^{-2\ll _i r}\big)\ll_i^{-1-\aa}=0.\end{align*}
Moreover,  \eqref{B22} and \eqref{UIB} imply
$$\lim_{r\to 0} \sup_{t\ge 1} \E^\mu[\{th(0)\W_2(\mu_t^B,\mu)^2-t\W_2(\mu_{t,r}^B,\mu)^2\}^-]=0.$$
Therefore, for any $\vv>0$,
\beg{align*} &\limsup_{t\to\infty}\P^\mu\big(\{th(0)\W_2(\mu_{t}^B,\mu)^2- \Xi^B(t)\}^-\ge 3\vv\big)\\
&\le \lim_{r\to 0} \limsup_{t\to\infty} \P^\mu\big(\{t\W_2(\mu_{t,r}^B,\mu)^2- \Xi_r^B(t)\}^-\ge \vv \big)\\
&+ \lim_{r\to 0}\sup_{t\ge 1} \Big[\P^\mu\big(|\Xi^B(t)- \Xi_r^B(t)|\ge \vv \big)
  + \P^\mu\big(\{th(0)\W_2(\mu_t^B,\mu)^2-t\W_2(\mu_{t,r}^B,\mu)^2\}^- \ge \vv \big)\Big]=0.
  \end{align*}
Combining this with \eqref{UIB} and applying the dominated convergence theorem, we prove   \eqref{ABC}.

 \end{proof}

 To prove Theorem \ref{TL2}, we need the following lemma.

 \beg{lem}\label{PH} Assume $\eqref{A10}$.
 \beg{enumerate} \item[$(1)$] Let ${\bf V}= {\bf V}_B$ for $B(\ll)=\ll$. We have
\beq\label{ZN1} {\bf V}(\phi_i)=  \ll_i^{-1}-\ll_i^{-2} {\bf V}(Z\phi_i),\ \ \ i\ge 1.\end{equation}
\item[$(2)$] If    $\eqref{A120}$ holds and $B\in{\bf B}^\aa\cap{\bf B}_{\aa'}$ for some $\aa,\aa'\in [0,1],$ then there exist  constants $c_1,c_2 >0$ such that
 \beq\label{ZN2} {\bf V}_B(\phi_i)\ge c_1 \ll_i^{-\aa'} - c_2 \ll_i^{-1-(\ff 1 2\land\aa)}\big[1+1_{\{\aa=\ff 1 2\}}\log(1+\ll_i)\big],\ \ \ i\ge 1.\end{equation}\end{enumerate}
 \end{lem}
 \beg{proof} (1) By $\hat L \phi_i=  -\ll_i\phi_i $ and using the Kolmogorov  equation,   we obtain
 \beq\label{PS} P_s\phi_i= -\ff 1 {\ll_i}  P_s \hat L \phi_i =  -\ff 1 {\ll_i} \ff{\d}{\d s}  P_s    \phi_i  + \ff 1 {\ll_i} P_s (Z   \phi_i).\end{equation}
 This together with  \eqref{RM} and $\mu(\phi_i^2)=1$ implies
 \beq\label{O2} \beg{split} & \int_0^{t} \mu( \phi_iP_{s}\phi_i )\d s=   -\ff 1 {\ll_i} \int_0^{t} \Big(\ff{\d }{\d s} \mu( \phi_i  P_s\phi_i)
 - \mu\big( \phi_i(Z P_{s}\phi_i) \big)\Big\}\d s \\
 &= \ff 1 {\ll_i}\Big(1- \mu( \phi_iP_{t} \phi_i)\Big) -\ff 1 {\ll_i} \int_0^{t} \mu\big( \{Z \phi_i\}  P_s \phi_i\big)\d s\\
 &=  \ff 1 {\ll_i}\Big( 1-\mu(\phi_iP_{t}\phi_i)\Big) +\ff 1 {\ll_i^2} \int_0^{t}\Big\{ \ff{\d } {\d s}  \mu\big(\{Z \phi_i\}  P_s \phi_i\big)
- \mu\big(\{Z \phi_i\}P_s\{Z  \phi_i\}\big)\Big\}\d s\\
 &= \ff 1 {\ll_i}\Big(1- \mu(\phi_iP_{t}\phi_i)\Big) +\ff 1 {\ll_i^2} \mu(\{Z  \phi_i\}P_{t}\phi_i)    -\ff 1{\ll_i^2} \int_0^{t}\mu\big(\{Z \phi_i\}P_s\{Z  \phi_i\}\big) \d s.  \end{split}\end{equation}
   By   \eqref{L-3} and $\|P_t-\mu\|_2\le \e^{-\ll_1 t}$,    we may let $t\to\infty$ to derive
 $${\bf V}(\phi_i)=\ff 1{\ll_i}- \ff 1 {\ll_i^2} {\bf V}(Z\phi_i),\ \ i\ge 1.$$

 (2) Let \eqref{A12} hold. By \eqref{A1-3}, \eqref{LT} and \eqref{PTP} we obtain
 \beq\label{F1} {\bf V}_B(\phi_i):= \int_0^\infty\mu(\phi_i P_t^B\phi_i)\d t
 = \ff 1 {B(\ll_i)} +\int_0^\infty \bigg\{ \E\int_0^{S_t^B} \e^{-\ll_i(S_t^B-s)}\mu\big(\phi_i P_s(Z\phi_i)\big)\d s\bigg\}\d t.\end{equation}
 By \eqref{A1-3} and \eqref{00} for $(P_t^*,-Z)$ replacing $(P_t,Z)$, we derive
 $$P_s^*\phi_i =\hat P_s \phi_i-\int_0^s P_r^*\{Z\hat P_{s-r} \phi_i\}\d r= \e^{-\ll_i s}\phi_i -\int_0^s \e^{-\ll_i(s-r)}P_r^*(Z\phi_i)\d r.$$
 This together with \eqref{RM} implies
 $$\mu(\phi_i P_s(Z\phi_i))= \mu((P_s^*\phi_i)  Z\phi_i) =-\int_0^s \e^{-\ll_i(s-r)}\mu((Z\phi_i)P_r(Z\phi_i))\d r.$$
 Combining this with  \eqref{A12'} and noting that \eqref{A1-3} yields 
  $$\|Z\phi_i\|_2^2\le \|Z\|_\infty \|\nn \phi_i\|_2^2=\|Z\|_\infty^2\ll_i,$$
  we derive 
 \beq\label{YY0} \beg{split}  \mu(\phi_i P_s(Z\phi_i))\le c(1) \ll_i^{\ff 1 2} \int_0^s \e^{-\ll_i(s-r)}r^{-\ff 1 2}\e^{-\ll r}\d r,\ \ s\ge 0, i\in \mathbb N.\end{split}\end{equation}
  By \eqref{V*1} for $s$ replacing $S_t^B$, we find a constant $a_1>0$ such that
  $$\int_0^s \e^{-\ll_i(s-r)-\ll r}r^{-\ff 1 2}  \d r \le a_1 \ll_i^{-1} s^{-\ff 1 2} \e^{-\ll s/2},\ \ s>0.$$
  Combining this with \eqref{YY0}, \eqref{V*1} and \eqref{V*2},  we find constants $a_2,a_3,a_4>0$ such that
  \beq\label{XD1}\beg{split} &I_i(t):= \E\int_0^{S_t^B} \e^{-\ll_i(S_t^B-s)}\mu(\phi_iP_s(Z\phi_i))\d s
  \le \ff{ a_1}{\ss{ \ll_i}} \E\int_0^{S_t^B} \e^{-\ll_i(S_t^B-s)}s^{-\ff 1 2} \e^{-\ll s/2}\d s\\
  &\le a_2 \ll_i^{-\ff 3 2} \E\Big[(S_t^B)^{-\ff 1 2} \e^{-\ll S_t^B/4}\Big] \le a_3 \ll_i^{-\ff 3 2}\Big(\ff 1 2\land  t\Big)^{-\ff 1 {2\aa} } \e^{-a_4 t},\ \ t>0, i\ge 1.\end{split}\end{equation}
  On the other hand, noting that
  $$\big|\mu\big((Z\phi_i)P_r(Z\phi_i)\big)\big|\le \|Z\|_\infty^2 \|\nn \phi_i\|_2^2 \e^{-\ll_1 r}=\ll_i
  \|Z\|_\infty^2 \e^{-\ll_1 r}, $$
  by \eqref{YY0} and \eqref{LT}, we find   constants $a_5,a_6>0$ such that
  \beg{align*}&I_i(t)\le \ll_i\|Z\|_\infty^2 \E\int_0^{S_t^B} \e^{-\ll_i(S_t^B-s)}\d s\int_0^s \e^{-\ll_i(s-r)-\ll_1 r}\d r\\
  &\le a_5 \E\int_0^{S_t^B} \e^{-\ll_i(S_t^B-s)-\ll_1 s/2}\d s \le a_5 \ll_i^{-1} \e^{-a_6t}.\end{align*}
Combining this with \eqref{XD1} and \eqref{XX},  we find a constant $a_7>0$ such that
\beg{align*}&\int_0^\infty I_i(t)\d t \le (a_3\lor a_5) \ll_i^{-\ff 3 2 }\int_0^\infty h_{i,\aa}(t)
\e^{-(a_4\land a_6)t}\d t\\
&\le a_7 \ll_i^{-1-(\ff 1 2\land \aa)}\big[1+1_{\{\aa=\ff 1 2\}}\log(1+\ll_i)\big].\end{align*}
This together with \eqref{F1} and $B\in {\bf B}_{\aa'}$ implies \eqref{ZN2}.
\end{proof}

 \beg{proof}[Proof of Theorem \ref{TL2}] Since $\aa'\in [0,1]$ and $\aa\in [0,\aa']\cap (\aa'-1,\aa']$ imply
 $2+(\ff 1 2\land \aa)>1+\aa',$ by \eqref{ZN2} we find   constants $a_0>1$ and $a_1,a_2\ge 0$ such that 
 $$\ff{{\bf V}_B(\phi_i)}{\ll_i}\ge a_1\ll_i^{-1-\aa'}- a_2 \ll_i^{-a_0(1+\aa')},\ \ \ i\in \mathbb N.$$
 Combining this with \eqref{A2'} and  
 $d'=2(1+\aa),$    we find constants $a_1',a_2',  a_3,a_4>0$ such that
 \beq\label{K0}\beg{split} &\eta^B_{Z,r}= \sum_{i=1}^\infty\e^{-2\ll_i r} \ff{{\bf V}_B(\phi_i)}{\ll_i}
 \ge \sum_{i=1}^\infty\e^{-2\ll_i r} \Big( a_1 \ll_i^{-1-\aa'} - a_2 \ll_i^{-a_0(1+\aa')} \Big)\\
 &\ge  \sum_{i=1}^\infty  \e^{-2r c_2i^{\ff 1{1+\aa'}}}\Big(a_1'   i^{-1} -a_2' 
 i^{-a_0}\Big)
  \ge a_3\log(1+r^{-1})-a_4,\ \ r\in (0,1]. \end{split}\end{equation}
Next, by \eqref{BL3}, we obtain
 \beq\label{Z1} \beg{split}&\E^\mu[\W_2(\mu_{t,r}^B,\mu)^2]\ge \E^\mu[1_{B_\si}\W_2(\mu_{t,r}^B,\mu)^2]\ge t^{-1}\E^\mu[1_{B_{t,r}(\si)}\Xi_r^B(t)]-c_3\si^{1+2\theta}\\
 &\ge t^{-1} \E^\mu[ \Xi_r^B(t)]- t^{-1} \E^\mu[1_{B_{t,r}(\si)^c}\Xi_r^B(t)]-c_3\si^{1+2\theta},\ \ t\ge 1, r,\si\in (0,1].\end{split}\end{equation}
 By \eqref{BL-0}, \eqref{A2'}   and $d'=2(1+\aa')$, we find constants $k_1,k_2 >0$ such that
  \beg{align*}& \E^\mu[\Xi_r^B(t)]= \sum_{i=1}^\infty
  \ff{\e^{-2\ll_i r}}{\ll_i}\E^\mu[\psi_i^B(t)^2]
 \ge \eta_{Z,r}^B - k_1 t^{-1}\sum_{i=1}^\infty \e^{-2r\ll_i} {\ll_i^{-(1+\aa)}}\\
 &\ge \eta_{Z,r}^B-k_1t^{-1}\sum_{i=1}^\infty \e^{-2rc_2 i^{\ff 1 {1+\aa'}}} c_2^{-1-\aa}i^{-\ff{1+\aa}{1+\aa'}} \ge \eta_{Z,r}^B-k_2t^{-1} r^{-k_2},\ \ t\ge 0, r\in (0,1].\end{align*}
 Combining this with \eqref{K0} and \eqref{Z1},  we derive
\beq\label{Z3}\beg{split} & \E^\mu[\W_2(\mu_{t,r}^B,\mu)^2]\ge a_3 t^{-1}\log(1+r^{-1})-a_4 t^{-1}-k_2t^{-2} r^{-k_2}\\
&\qquad  - t^{-1} \E^\mu[1_{B_{t,r}(\si)^c}\Xi_r^B(t)]-c_3\si^{1+2\theta},\ \ t\ge 1, r,\si\in (0,1].\end{split}\end{equation}
 On the other hand, by \eqref{LY1} and \eqref{PR}, we find constants $k_3,k_4>0$ such that
 \beg{align*}&\E^\mu[\Xi_r^B(t)^2]= \E^\mu[\|(-\hat L)^{-\ff 1 2} (f_{t,r}^B-1)\|_2^4]\le k_3 r^{-k_4},\\
 &\P^\mu(B_{t,r}(\si)^c)\le k_3 r^{-k_4} t^{-1} \si^{-(1+\theta)}\ \ t\ge 1, r\in (0,1],\si\in (0,1],\end{align*}
 where $\theta>0$ is a constant. Thus,
 $$\E^\mu[1_{B_{t,r}(\si)^c}\Xi_r^B(t)]\le  \ss{\P^\mu(B_{t,r}(\si)^c)   \E^\mu[\Xi_r^B(t)^2]]}
 \le k_3 r^{-k_4} t^{-\ff 1 2}\si^{-\ff{ 1+\theta }2},\ \ t\ge 1, r,\si\in (0,1].$$
 Combining this with \eqref{Z3}, we arrive at
  \beg{align*} & \E^\mu[\W_2(\mu_{t,r}^B,\mu)^2]\ge a_3 t^{-1}\log(1+r^{-1})-a_4t^{-1}-k_2 t^{-2}r^{-k_2}\\
  &\qquad -k_3 r^{-k_4} t^{-\ff 32}\si^{-\ff{ 1+\theta }2} -c_3\si^{1+2\theta},\ \ t\ge 1, r,\si\in (0,1].\end{align*}
  Taking $\si= t^{-\ff 1 {1+2\theta}}$ and $r= t^{-\ff 1 {k_2}\land\ff{\theta}{ 2 k_4(1+2\theta)}}$, obtain
$$\E^\mu[\W_2(\mu_{t,r}^B,\mu)^2]\ge \ff{a_3}{1+2\theta} t^{-1}\log t - (k_2+k_3+c_3)t^{-1},\ \ t\ge 1.$$
  Therefore, \eqref{BLa} holds for some constants $c,t_0>0$.

 \end{proof}

 Finally, to prove Theorem \ref{TL3}, we present one more lemma.

\beg{lem} \label{BL} Let $(E,\rr)$ be a Polish space. Let
   $X_t$ be a  continuous time Markov process on $E$   such that the associated semigroup $P_t$ satisfies
\beq\label{EXX} \|P_t-\mu\|_{2}\le c_1\e^{-\ll_1 t},\ \ t\ge 0\end{equation}   for
 some constants $c_1,\ll_1>0$ and a probability measure $\mu$ on $E$. If there exists $\phi\in C_{b,L}(E)$ such that $\mu(\phi)=0$ and
 $${\bf V}(\phi):=\int_0^\infty \mu(\phi P_t\phi) \d t>0,$$
 then  there exist  constants $c,t_0>0$ such that $\mu_t:=\ff 1 t\int_0^t \dd_{X_s}\d s$ satisfies
\beq\label{W1}   \E^\mu[\W_1(\mu_t,\mu)]\ge c t^{-\ff 1 2},\ \ t\ge t_0.\end{equation}
If $M\subset E$ such that $\mu(M)=1$ and
$\|P_t\|_{1\to 2}<\infty$ for $t>0$, then \eqref{W1} holds for $\inf_{\nu\in \scr P}\E^\nu$ replacing $\E^\mu$, where $\scr P$ is the set of all probability measures on $M$. \end{lem}

\beg{proof}     By \cite[Theorem 2.1(c)]{Wu}, we have
\beq\label{CL} \lim_{t\to\infty } \ss t\E^\mu[|\mu_t(\phi)|]= \big(2\pi {\bf V}(\phi)\big)^{-\ff 1 2}
\int_{-\infty}^\infty |r| \e^{-\ff{r^2}{2{\bf V}(\phi)}}\d r.\end{equation}
So, by the Kantorovich dual formula, there exist constants $c,t_0>0$ such that
$$\ss t \E^\mu[\W_1(\mu_t,\mu)]\ge  \ss t\E^\mu[|\mu_t(\phi)-\mu(\phi)|]
\ge c,\ \ t\ge t_0.$$

Next, let $M\subset E$ such that $\mu(M)=1$ and
$\|P_t\|_{1\to 2}<\infty$ for $t>0$. Then
$$\{\nu P_1: \nu\in \scr P\}\subset \big\{\nu\in \scr P: \nu=h\mu, \|h\|_2\le \|P_1\|_{1\to 2}\big\},$$
so that   \cite[Theorem 2.1(c)]{Wu}, $\hat \mu_t:= \ff 1 t \int_{1}^{t+1}\dd_{X_s}\d s$ satisfies
$$\lim_{t\to\infty }\inf_{\nu\in \scr P}  \ss t\E^\nu[|\hat\mu_t(\phi)|]>0.$$
Noting that
$$|\hat \mu_t(\phi)-\mu_t(\phi)|\le  \|\phi\|_\infty t^{-1},\ \ t>1,$$
we obtain
$$\lim_{t\to\infty }\inf_{\nu\in \scr P}  \ss t\E^\nu[|\mu_t(\phi)|]>0.$$
By Kantorovich's dual formula and $\mu(\phi)=0$, this implies \eqref{W1} for $\inf_{\nu\in \scr P}\E^\nu$ replacing $\E^\mu$.

\end{proof}

\beg{proof}[Proof of Theorem \ref{TL3}] (1) By \eqref{A10} and \eqref{A120}, for any $i\ge 1$  we have
$\|\phi_i\|_\infty<\infty$ and there exist constants $c_1,c_2>0$ such that
\beg{align*}&\|\nn\phi_i\|_\infty= \e^{\ll_i} \|\nn \hat P_1 \phi_i\|_\infty\le c_1 \e^{\ll_i} \big\|(P_1|\nn\phi_i|^2)^{\ff 1 2}\big\|_\infty\\
&\le c_2 \e^{\ll_i} \|\nn\phi_i\|_2= c_2\e^{\ll_i}\ss{\ll_i}<\infty.\end{align*}
So, $\phi_i\in C_{b,L}(M)$, and hence is uniquely extended to $\phi_i\in C_{b,L}(E)$ for $E:=\bar M$.
On the other hand, by \eqref{ZN2} we have ${\bf V}_B(\phi_i)>0$ for large enough $i$.
Moreover,  \eqref{UAA} implies
  $\|P_t^B\|_{1\to 2}<\infty$ for $t>0$. So, the first assertion follows from Lemma \ref{BL} with $E=\bar M.$

(2) Let $B\in {\bf B}_\aa$ for some $\aa\in [0,1]$ with $d''>2(1+\aa)$.
 By \eqref{A13} for $p=1$ we find a constant $c_1>0$ such that
  $$  \E^\mu[\rr(\hat X_t,\hat X_0)]\le c_1 t^{\ff 1 2},\ \ t\ge 0.$$
  Combining this with \eqref{00} and \eqref{A120}, we find a constant $c_2>0$ such that
  \beg{align*} &\E^\mu[\rr(X_t,X_0)]=\int_M   P_t \rr(x,\cdot)(x)\mu(\d x)\\
  &= \E^\mu[\rr(\hat X_t,\hat X_0)]+ \int_M\mu(\d x) \int_0^t P_s\{Z\hat P_{t-s}\rr(x,\cdot)\}(x)\d s\\
  &\le c_1 t^{\ff 1 2}+ c_2\int_0^t (t-s)^{-\ff 1 2}\d s\le c_3 t^{\ff 1 2},\ \ t\ge 0.\end{align*}
  According to   the proof of \cite[Theorem 1.1(2)]{WW}, this and  \eqref{MU} imply
  \beq\label{LLB} \E^\mu[\W_1(\mu_t^B,\mu)]\ge c_4 t^{-\ff 1 {d''-2\aa}},\ \ \ t\ge t_1\end{equation}
 for some constants $c_4,t_1>0$.

 Finally, by \eqref{UAA}, we find a constant $t_2>0$ such that $\|P_{t_2}-\mu\|_{1\to\infty}\le \ff 1 2$, so that
 $$\nu_{t_2}:= \nu P_{t_2}^B \ge \ff 1 2 \mu,\ \ \ \nu\in \scr P.$$
 Let $\mu_t^{B,t_2}$ be in \eqref{BME} for $\vv=t_2$. Then   the Markov property and   \eqref{LLB} yield
 $$\inf_{\nu\in \scr P}\E^\nu[\W_1(\mu_t^{B,t_2},\mu)]=\inf_{\nu\in \scr P}\E^{\nu_{t_2}} [\W_1(\mu_t^{B},\mu)]\ge \ff 1 2 c_4 t^{-\ff 1 {d''-2\aa}},\ \ \ t\ge t_1.$$
 Combining this with \eqref{BWE}, the triangle inequality and $d''-2\aa>2$, we find constants
 $c_5, c, t_0>t_1$ such that
 $$\inf_{\nu\in \scr P}\E^\nu[\W_1(\mu_t^{B},\mu)]\ge \ff 1 2 c_4 t^{-\ff 1 {d''-2\aa}}- c_5 t^{-1}
 \ge c t^{-\ff 1 {d''-2\aa}},\ \ \ t\ge t_0.$$

 \end{proof}

  \section{Some concrete models}

In this part, we apply our general results to some typical models including:
1) the (reflecting) subordinated diffusion process on a compact manifold; 2) the subordinated conditional diffusion process on a bounded  open domain; 3)
the subordinated Wright-Fisher diffusion process; 4) the subordinated subelliptic diffusion process on $\SU(2)$. It is also possible to consider  more general hypoelliptic diffusion processes studied in \cite{FN,FB2} under  the generalized curvature-dimension conditions.
For simplicity, throughout this section, we take
$$B\in {\bf B}^\aa\cap {\bf B}^\aa\  \text{for\ some}\  \aa\in (0,1].$$

  \subsection{Subordinated (Reflecting) diffusion process}
In this part, we consider the model stated in Introduction, for which all conditions in Theorems \ref{TU1}-\ref{TL3} are satisfied for $d=d'=d''=n$.

 Indeed, \eqref{A10} and \eqref{A2'} are well known (see \cite{Chavel, Davies}), \eqref{A120} follows from \cite[Lemma 2.1]{W05}, and \eqref{B22} with $h(r)=\kk \e^{Kr}$ is implied  by \cite[(3.36), (3.37)]{WZ}, where $\kk\ge 1$ and $K\ge 0$ are constants with $\kk=1$ when $\pp M$ is empty or convex.
  Moreover, the following lemma confirms other conditions.

\beg{lem}\label{XPP} $\eqref{B*}$, $(A_2)$ and \eqref{RS} hold. \end{lem}

\beg{proof} (1) Let $l_t$ be the local time of $\hat X_t$ on $\pp M$ if $\pp M$ exists, and let $l_t=0$ otherwise. By \cite[(2.1)]{W05} and the proof of \cite[Lemma 2.1]{W05}, there exist constants $c_1,K,\dd>0$ such that
\beq\label{K1} |\nn \hat P_t f(x)|\le \E^x[|\nn f(\hat X_t)|\e^{Kt+\dd l_t}],\ \ t\ge 0, x\in M, f\in C_b^1(M),\end{equation}
\beq\label{K2} \sup_{x\in M} \E^x[l_t^2]\le c_1 t,\ \ \ \sup_{x\in M}\E^x[\e^{\ll l_t}]<\infty,\ \ \ll, t\ge 0.\end{equation}
By the Schwarz  inequality, \eqref{K1} implies
\beg{align*}& |\nn \hat P_t \e^f (x)|\le \hat P_t|\nn \e^f|(x) + \E^x \big[|\nn\e^f(\hat X_t)|(\e^{Kt+\dd l_t}-1)\big]\\
&\le \big\{\hat P_t \e^f(x)\big\}^{\ff 1 2}\big\{\hat P_t(|\nn f|^2\e^f)(x)\big\}^{\ff 1 2}
+\|\nn f\|_\infty (\hat P_t\e^{2f})^{\ff 1 2} \big(\E^x[\e^{2Kt+2\dd l_t}-1]\big)^{\ff 1 2}.\end{align*}
On the other hand, by \eqref{K2} we find a constant $c_2>0$ such that
\beg{align*}& \E^x[\e^{2Kt+2\dd l_t}-1]\le \E^x[(2K+2\dd l_t)\e^{2Kt+2\dd l_t}]\\
&\le \big(\E^x[(2K+2\dd l_t)^2]\big)^{\ff 1 2}\big(\E^x[\e^{4Kt+4\dd l_t}]\big)^{\ff 1 2}\le c_2 t^{\ff 1 2},\ \ t\in [0,1].\end{align*}
Therefore, \eqref{B*} holds for $\theta=\ff 1 2$ and $m=1$.

(2) When $\pp M$ is empty or convex, we have
$$\<\nn \rr(\hat X_0,\cdot), {\bf N}\>(\hat X_t) \d l_t\le 0,$$
where ${\bf N}$ is the inward unit normal vector on $\pp M$. On the other hand, by the Laplacian comparison theorem, there exists a constant $c_0>0$ such that
$$\hat L \rr(\hat X_0,\cdot)^2\le c_0.$$ So, by It\^o's formula, we find a constant $c_1>0$ such that
$$\d \rr(\hat X_0,\hat X_t)^2 \le c_1 \d t + 2\ss 2 \rr(\hat X_0,\hat X_t)\d B_t,\ \ t\ge 0,$$
where $B_t$ is the one-dimensional Brownian motion.
Thus, for any $p\ge 1,$ there exists a constant $c(p)>0$ such that
$$\d \rr(\hat X_0,\hat X_t)^{2 p} \le c(p)\rr(\hat X_0,\hat X_t)^{2 (p-1)} \d t + d M_t $$
holds for some martingale $M_t$, so that
$$\E[\rr(\hat X_0,\hat X_t)^{2 p}]\le c(p)\int_0^t \E[\rr(\hat X_0,\hat X_s)^{2 (p-1)}]\d s,\ \ t\ge 0.$$
Consequently, for $p=1$ we get
$$\E^x[\rr(\hat X_0,\hat X_t)^{2}]\le c(1) t,$$
and by inducting in $p\in \mathbb N$,   we derive \eqref{A13} for $p\in 2\mathbb N$. Therefore, \eqref{A13} holds for all $p\ge 1$ due to Jensen's inequality.

When $\pp M$ is non-convex, as explained in the proof of  \cite[Proposition 3.2.7]{Wbook}, there exists a function $1\le \phi\in C_b^\infty(M)$ such that $\pp M$ is convex under the metric
$$\<\cdot,\cdot\>':= \phi^{-1} \<\cdot,\cdot\>,$$
and $\hat L=\phi^{-2}\DD'+Z',$ where $\DD'$ is the Laplacian induced by the new metric and $Z'$ is a $C_b^1$ vector field. Let $\rr'$ be the Riemannian distance induced by the new metric, we have
$\rr\le \|\phi\|_\infty \rr'$, so that   the above argument for convex $\pp M$ leads to
$$\E[\rr(\hat X_0,\hat X_t)^{2p}]\le \|\phi\|_\infty^{2p} \E[\rr'(\hat X_0,\hat X_t)^{2p}]
\le k(p) t^p,\ \ t\in [0,1].$$

(3)
To verify \eqref{RS}, we follow the line of \cite{AC}. As explained in the end of page 12 in \cite{AC}, see also the proof of \cite[Theorem 1.5]{AC}, under  \eqref{A12},   it remains to verify the volume doubling condition and scaled Poincar\'e inequalities on balls. More precisely, we only need to find a distance $\tt\rr$ and constants $c_1,c_2,c_3>0$ such that
$$c_1\tt\rr\le \rr\le c_2\tt\rr,$$
and   the balls  $\tt B(x,r):= \{y\in M: \tt\rr(x,y)\le r\}$ for all $x\in M$ and $r>0$ satisfy
\beq\label{VL} \mu(\tt B(x,2 r))\le c_3 \mu(\tt B(x,r)),\end{equation}
 \beq\label{SP} \mu(1_{\tt B(x,r)}f^2)\le c_3 r^2 \mu(1_{\tt B(x,r)}|\nn f|^2),
 \ \ \ f\in C_b^1(M),  \mu(1_{\tt B(x,r)}f)=0.\end{equation}
Since for the present model we have $\tt D:=\|\tt\rr \|_\infty <\infty$ and
$kr^d \le \mu(\tt B(x,r))\le K r^d$   for some constants $K>k>0$ and all $r\in [0,\tt D]$, \eqref{VL} holds true.
To verify \eqref{SP}, by   the conformal change of metric used in step (2),  we may and do assume that
$\pp M$ is either empty or convex. In this case, there exists a constant $r_0>0$ such that
$B(x,r):=\{y\in M:\rr(x,y)\le r\}$ is convex for all $x\in M$ and $r\in (0,r_0].$
Then we take $\tt\rr:= \rr\land r_0$, so that $\tt B(x,r)$ is convex for  all $r>0$, since $\tt B(x,r)=B(x,r)$ for $r<r_0$ and $\tt B(x,r)=M$ for $r\ge r_0$.  Thus, \eqref{SP} follows from \cite[Theorem 1.4]{W94}.

\end{proof}

According to the above observations, we  conclude that
 all assertions in Theorems \ref{TU1}-\ref{TL3} hold for $d=d'=d''=n$, i.e.
\beq\label{5A1}\beg{split} &q_\aa= \ff{2n}{(3n-2-2\aa)^+},\  \ \ \aa(d,d')=\aa(n):= \ff 1 2 \ss{1+2n(n-1)}-\ff 12, \\ &\gg_{\aa,p,q}=\ff n 2 (3-p^{-1}-q^{-1})-1-\aa,\ \ \ p,q\in [1,\infty).\end{split} \end{equation}
In this case, $\aa>\aa(n)$ implies $n<2(1+\aa)$ and $q_\aa>\ff n{2\aa},$ so that in Theorem \ref{TU1}(1) we only need $\aa>\aa(n).$
Below we summarize these results, which in particular imply Theorems \ref{TN} and \ref{T1} stated in Introduction, according to  $B(\ll)=\ll (\aa=1)$, \eqref{ZN1} and
$${\bf V}(Z\phi_i):= \mu((Z\phi_i) (-L)^{-1} (Z\phi_i))= \mu(|\nn L^{-1} (Z\phi_i)|^2).$$

\beg{thm} Let $q_\aa, \aa(n),\gg_{\aa,p,q}$  for $p,q\in [1,\infty)$ be in \eqref{5A1}. Then the following assertions hold for some constant $\kk\in [1,\infty)$, where $\kk=1$ when $\pp M$ is either empty or convex.
\beg{enumerate} \item[$(1)$] If  $\aa> \aa(n)$, then
$$\lim_{t\to\infty} \sup_{\nu\in\scr P} \E^\nu\big[\big|\{t\W_2(\mu_t^B,\mu)^2-\Xi^B(t)\}^+
+ \{t\kk\W_2(\mu_t^B,\mu)^2-\Xi^B(t)\}^-\big|^q\big]=0,\ \ q\in [1,q_\aa).$$
\item[$(2)$] If $n<2(1+\aa)$, then for any $q\in [1,q_\aa)$ and $k\in (\ff n {2\aa \i(q)},\infty]\cap [1,\infty]$,
 $$ \lim_{t\to\infty} \sup_{\nu\in \scr P_{k,R}} \E^\nu\Big[\Big|\{t\W_2(\mu_t^B,\mu)^2-\Xi^B(t)\}^++ \{t\kk\W_2(\mu_t^B,\mu)^2-\Xi^B(t)\}^-\Big|^q\Big]=0,\ \   R\in (0,\infty).$$
 \item[$(3)$] If $n<2(1+\aa)$, then
 $$\kk^{-1}\eta_Z^B\le \liminf_{t\to\infty} \inf_{\nu\in \scr P}t \E^\nu[\W_2(\mu_t^B,\mu)^2]
 \le  \limsup_{t\to\infty} \sup_{\nu\in \scr P}t \E^\nu[\W_2(\mu_t^B,\mu)^2]\le \eta_Z^B<\infty.$$
 \item[$(4)$] If $n=2(1+\aa)$, i.e. $(n,\aa)\in\{(3,\ff 1 2), (4,1)\}, $ then there exist constants $c, t_0>1$ such that
 $$c^{-1} t^{-1}\log t\le  \E^\nu[\W_2(\mu_t^B,\mu)^2]\le c t^{-1}\log t,\ \ t\ge t_0, \nu\in \scr P.$$
  \item[$(5)$] If $n>2(1+\aa)$, then there exist constants $c, t_0>1$ such that
  $$c^{-1} t^{-\ff 2{n-2\aa}}\le \big(\E^\nu[\W_1(\mu_t^B,\mu)]\big)^2\le \E^\nu[\W_2(\mu_t^B,\mu)^2]\le
  c t^{-\ff 2{n-2\aa}},\ \ t\ge t_0,\nu\in \scr P.$$
 \item[$(6)$] If  $p,q\in [1,\infty)$ with $n(3-p^{-1}-q^{-1})<2(1+\aa),$ then  there exist  constants $c,t_0>1$ such that
 $$c^{-1}t^{-1}\le \big(\E^\nu[\W_1(\mu_t^B,\mu)]\big)^2\le \big(\E^\nu[\W_{2p}(\mu_t^B, \mu)^{2q}]\big)^{\ff 1 q}\le c t^{-1},\ \ t\ge t_0, \nu\in \scr P.$$
  \item[$(7)$] If $p,q\in [1,\infty)$ such that $n(3-p^{-1}-q^{-1})\ge 2(1+\aa)$, then  there exists a constant $c>0$ such that for any $t\ge 1$,
      $$ \sup_{\nu\in \scr P} \big(\E^\nu[\W_{2p}(\mu_t^B, \mu)^{2q}]\big)^{\ff 1 q}\le \beg{cases}c t^{-1}\log(1+t),\ &\text{if}\ n(3-p^{-1}-q^{-1})=2(1+\aa),\\
        c t^{-\ff 1{1+\gg_{\aa,p,q}}},\ &\text{if}\ n(3-p^{-1}-q^{-1})>2(1+\aa).\end{cases}$$
 \end{enumerate}\end{thm}

 \subsection{Subordinated conditional diffusion process}

 Let   $M$   be a bounded connected   $C^{2}$    open domain in an $n$-dimensional  complete Riemannian manifold,    and let   $V\in C_b^2(M)$  be such that
   $\mu_0(\d x):=\e^{V(x)}\d x$   is a probability measure on   $M$. Consider the Dirichlet eigenproblem for $\hat L_0:= \DD+\nn V$  in   $M$:
 $$\hat L_0 h_i=-\theta_i h_i,\ \  \ h_i|_{\pp M}=0, \ \ \ i\ge 0,$$
where   $\{\theta_i\}_{i\ge 0}$   are listed in the increasing order counting multiplicities, and   $\{h_i\}_{i\ge 0}$  are the associated unitary eigenfunctions in $L^2(\mu_0)$  with   $h_0>0$. Let
  $$\hat L:= \hat L_0+2\nn \log h_0=\DD+\nn (V+2\log h_0),\ \ \mu(\d x):= h_0(x)^2\mu_0(\d x).$$
Then the diffusion process $\hat X_t$ generated by $\hat L$ is non-explosive in $M$, whose distribution coincides
with the  conditional distribution of the $\hat L_0$-diffusion process $\hat X_t^0$ under the condition that
 $$\tau:=\inf\{t\ge 0: \hat X_t^0\in\pp M\}=\infty,$$
  in the sense that for any $T>0$ and any $F\in C_b((C[0,T]; M)),$
$$\E[F(\hat X_{[0,T]})] =\lim_{m\to\infty} \E[F(\hat X_{[0,T]}^0)|\tau>m].$$
  Let $Z$ be a $C_b^1$-vector field on $M$ satisfying \eqref{RM}.

It is well known   that  $\{\ll_i:= \theta_i-\theta_0\}_{i\ge 0}$ are all eigenvalues of $-\hat L$ with unitary eigenfunctions $\{\phi_i:=h_ih_0^{-1}\}_{i\ge 0},$
and  that \eqref{A10}, \eqref{A2'} and \eqref{MU} hold  for
$$d=n+2,\ \ \ \ d'=d''=n,$$
see for instance  \cite{Chavel, OW}.
Next, by \cite[Lemma 4.6]{W2}, \eqref{B22} holds with
 $h(r)=\kk\e^{Kr}$ for some constants $\kk\ge 1$ and $K\ge 0$, where $\kk=1$ when $\pp M$ is convex.
   The following lemma confirms other conditions in Theorems \ref{TU1}-\ref{TL3}, except \eqref{RS} which is not yet verified.
 
 \beg{lem} For the present model, \eqref{A120} and $(A_2)$ hold. When $\pp M$ is convex or $n\le 3$,   \eqref{B*} is satisfied.
 \end{lem}
 \beg{proof}  According to the proof of  \cite[Lemma 4.6]{W2}, if $\pp M$ is convex then the Bakry-Emery curvature
  $\hat L$ is bounded from below by a constant $-K$, so that
 $$|\nn P_t f|\le \e^{Kt}P_t|\nn f|,$$
 which implies
 \eqref{A120}   for $k(p)= \e^{K}$ as well as \eqref{B*} for $\theta=1$. So, in the following we only prove these conditions for non-convex $\pp M$, and   verify $(A_2)$.

 (1) When $\pp M$ is non-convex, let $\rr'$ be the Riemannian distance induced by $\<\cdot,\cdot\>'$ introduced in the     proof of \eqref{XPP}. According to   the proof of  \cite[Lemma 4.6]{W2},  for any $x,y\in M$, there exists a coupling $(\hat X_t,\hat Y_t)$ of the  diffusion process generated by $\hat L$   starting from $(x,y)$, such that for some constant $c_1>0$ we have
 $$\d \rr'(\hat X_t,\hat Y_t)\le c_1 \rr'(\hat X_t,\hat Y_t)\d t +\d M_t,$$
 where $M_t$ is a martingale with $\d\<M\>_t\le c_1 \rr'(\hat X_t,\hat Y_t)\d t.$
 Thus, for any $q\in (1,\infty)$, there exists a constant $K(q)>0$  such that
 $$\big(\E [\rr'(\hat X_t,\hat Y_t)^q]\big)^{\ff 1 q}\le K(q) \rr'(x,y),\ \ t\in [0,1].$$
 Therefore, by $\rr'\le\rr\le \|\phi\|_\infty\rr'$,
\beg{align*} &|\nn \hat P_t f(x)|:=\limsup_{y\to x} \ff{|\hat P_tf(x)-\hat P_tf(y)|}{\rr(x,y)}
 \le \limsup_{y\to x} \E\Big[ \ff{|f(\hat X_t)-f(\hat Y_t)|}{\rr'(\hat X_t,\hat Y_t)}\cdot\ff{\rr'(\hat X_t,\hat Y_t)}{\rr(x,y)}\Big]\\
 &\le \limsup_{y\to x} \Big(\E\Big[\ff{ |f(\hat X_t)-f(\hat Y_t)|^p}{\rr'(\hat X_t,\hat Y_t)^p}\Big]\Big)^{\ff 1 p} \ff{(\E[\rr'(\hat X_t,\hat Y_t)^{\ff p{p-1}}])^{\ff {p-1}p} }{\rr(x,y)}
 \le  K\Big(\ff p{p-1}\Big)  \|\phi\|_\infty( P_t|\nn f|^p)^{\ff 1 p}.\end{align*}
So, \eqref{A120} holds.

 (2) Let $\pp M$ be non-convex and $n\le 3$. Since $\nn V\in C_b^1(M)$ and
 $$-\Hess_{\log h_0} = -\ff{\Hess_{h_0}}{h_0}+\ff{(\nn h_0)\otimes(\nn h_0)}{h_0^2} \ge -\ff{\|\nn^2h_0\|_\infty}{h_0},$$
 there exists a constant $c_1>0$ such that the Bakry-Emery curvature of $\hat L$ is bounded below by $-\ff {c_1}{2h_0},$ i.e.
 $$\GG_2(g):=\ff 1 2  \hat L|\nn g|^2 - \<\nn g, \nn \hat Lg\> \ge -\ff{c_1}2 h_0^{-1} |\nn g|^2,\ \ g\in C(M).$$
 So, by \eqref{UAA'} for $d=n+2$, and applying Jensen's inequality, we find a constant $c_2>0$ such that
 \beg{align*} & |\nn \hat P_t \e^f|^2 -\hat P_t|\nn \e^f|^2 =-\int_0^t\ff{\d}{\d s} \hat P_s|\nn \hat P_{t-s} \e^f|^2\d s\\
 &\le c_1 \int_0^t \hat P_s\big\{h_0^{-1} |\nn   \hat P_{t-s}\e^f|^2\big\}\d s\le
 c_1\int_0^t (\hat P_sh_0^{-\ff m{m-1}})^{\ff{m-1}m} (\hat P_s|\nn \hat P_{t-s}\e^f|^{2m})^{\ff 1 m}\d s\\
 & \le c_2\int_0^t s^{-\ff{(n+2)(m-1)}{2m}} \mu(h_0^{-\ff m{m-1}})^{\ff{m-1}m} (\hat P_s|\nn \hat P_{t-s}\e^f|^{2m})^{\ff 1 m}\d s\\
 &\le c_2\|\nn f\|_\infty^2(\hat P_t\e^{2m f})^{\ff 1 m} \int_0^t s^{-\ff{(n+2)(m-1)}{2m}} \mu(h_0^{-\ff m{m-1}})^{\ff{m-1}m} \d s.\end{align*}
Noting that
 $n+2\le 5$ and $\mu(h_0^{-r})=\mu_0(h_0^{2-r})<\infty$ for $r<3$,
 for any $m\in (\ff 3 2,\ff 5 3)$, we have
 $$\theta_1:= 1-\ff{(n+2)(m-1)}{2m}\in (0,1),\ \ \mu\big(h_0^{-\ff m{m-1}}\big)<\infty,$$
 so that for some  constant $c_3>0$ we have
 \beq\label{FN1} |\nn \hat P_t \e^f|^2 -\hat P_t|\nn \e^f|^2\le c_3 t^{\theta_1} \|\nn f\|_\infty^2(\hat P_t\e^{2m})^{\ff 1 m}.
 \end{equation}
 Similarly, by \eqref{A120} and its consequence
 $$|\nn \hat P_tg|\le \ff c {\ss t}(\hat P_t|\nn g|^2)^{\ff 1 2}, $$ we find a constant $c_4>0$ such that
 \beg{align*} & \hat P_t |\nn \e^f|^2- (\hat P_t\e^f)\hat P_t(|\nn f|^2 \e^f)=
 \int_0^t\ff{\d}{\d s}\hat P_{t-s} \big\{(\hat P_s\e^f)\hat P_s(|\nn f|^2\e^f)\big\}\d s\\
 &=-\int_0^t\hat P_{t-s} \<\nn\hat P_s\e^f,\nn\hat P_s(|\nn f|^2\e^f)\>\d s\\
& \le c_4 \|\nn f\|_\infty^2 \int_0^t   s^{-\ff 1 2}\hat P_{t-s}  (\hat P_s\e^{2f}) \d s
= 2c_4 t^{\ff 1 2} \|\nn f\|_\infty^2\hat P_t\e^{2f}.\end{align*}
This together with \eqref{FN1} implies \eqref{B*} for   $\theta=\theta_1\land \ff 1 2.$

(3) It remains to verify \eqref{A13}.
By the conformal change of metric as in the end of the proof of Lemma \ref{XPP}, we only consider the case where $\pp M$ is convex, so that
$$\<{\bf N}, \nn\rr(x,\cdot)\>|_{\pp M}\le 0,\ \ x\in M,$$
where ${\bf N}$ is the inward unit normal vector field of $\pp M$.
 Let $\rr_\pp$ be the distance to $\pp M$. It is well known (see for instance \cite{OW}) that
 $\nn h_0$ is inward normal on the boundary and $c_1:=\|h_0^{-1}\rr_\pp\|_\infty<\infty.$ So,
 \beq\label{XP1} \rr_\pp\le c_1 h_0,\ \  \<\nn h_0, \nn\rr(x,\cdot)\>|_{\pp M}\le 0,\ \ x\in M.\end{equation}
 Moreover, by the Hessian comparison theorem, there exists a constant $c_2>0$ such that
 \beq\label{XP2} \Hess_{\rr(x,\cdot)^2}(v,v)\le c_2 |v|^2,\ \ \ x\in M,\ v\in TM.\end{equation}
 We intend to show that these two estimates imply
 \beq\label{XP} \sup_{x,y\in M} \<\nn\log h_0, \nn \rr(x,\cdot)^2\>(y)\le c\end{equation}
 for some constant $c>0$. To see this, for any $x,y\in M$, let $z\in\pp M$ such that
 $\rr(y,z)= \rr_\pp (y)$. Let
 $$\gg: [0,1]\to M, \ \ \gg_0=z,\ \ \gg_1=y,\ \ |\dot\gg|=\rr_\pp(y)$$
 be the minimal geodesic from $z$ to $y$. Let $v_s=\nn h_0(\gg_s)$. We have
 $$v_0= a_0\dot\gg_0, \ \ \ \parallel_{0\to s} v_0= a_0\dot\gg_s,\ \ \ s\in [0,1],$$
 where $a_0:=\<{\bf N},\nn h_0\>(z)$ and $\parallel_{0,s}$ is the parallel displacement along the geodesic $\gg_s$.
 Since $h_0\in C_b^2(M)$, we find a constant
 $c_3>0$ such that
 $$|v_1- a_0 \dot \gg_1|= |v_1-\parallel_{0\to 1} v_0|\le c_3\rr_\pp(y).$$ Combining this with \eqref{XP1} and \eqref{XP2}, and noting that $|\nn \rr^2|\le 2\|\rr\|_\infty<\infty,$ we find  a constant $c_4>0$ such that
 \beg{align*} &\<\nn h_0, \nn\rr(x,\cdot)^2\>(y) = \<v_1, \nn\rr(x,\cdot)^2(\gg_1)\>
 \le a_0\<\dot\gg_1, \nn\rr(x,\cdot)^2(\gg_1)\>+  c_3
 \rr_\pp(y)\\
& =  a_0\<\dot\gg_0, \nn\rr(x,\cdot)^2(\gg_0)\>+a_0\int_0^1 \ff{\d}{\d s} \<\dot\gg_s, \nn\rr(x,\cdot)^2(\gg_s)\>\d s +c_3\rr_\pp(y)\\
 &\le a_0\int_0^1 \Hess_{\rr(x,\cdot)^2}(\dot\gg_s,\dot\gg_s)\d s +c_3\rr_\pp(y) \le a_0 c_2 \rr_\pp(y)^2+ c_3 \rr_\pp(y)\le c_4h_0(y).\end{align*}
 Therefore, \eqref{XP} holds for $c=c_4.$

By \eqref{XP} and It\^o's formula, we obtain
$$\d \rr(\hat X_0,\cdot)^2(\hat X_t)\le c \d t +2\ss 2 \rr(\hat X_0,\hat X_t)\d B_t,$$
where $B_t$ is the one-dimensional Brownian motion. This implies $(A_2)$ as explained in the proof of
Lemma \ref{XPP}(2).
 \end{proof}

We now conclude that
   all assertions in Theorems \ref{TU1}-Theorem \ref{TL2} hold, except Theorem \ref{TU2}(4) where the condition \eqref{RS} is to be verified for this model,   for $d=n+2, d'=d''=n$ and
 \beq\label{QU} \beg{split} &q_\aa:= \ff{2n+4}{3n+2-2\aa}\le 2,\ \ \ \aa(d,d')=\tt \aa(n):=\ff 1 2\ss{4+2n(n+2)}-1,\\
 &\ \ \ \gg_{\aa,p,q}:=\ff n 2 +\ff{n+2}2 (2-p^{-1}-q^{-1})-\aa-1.\end{split}\end{equation}
Noting that when $n=1$ the condition $\aa>\tt\aa (n)$ becomes $\aa>\ff 1 2 \ss{10}-1,$ which implies $1=n<2(1+\aa)$ and $\ff{6}{5-2\aa}=q_\aa>\ff n{2\aa}=\ff 1 {2\aa},$  while $q_\aa\le 2$ yields ${\rm i}(q)=1$ for $q\in [1,q_\aa)$,  we have the following result according to  Theorems \ref{TU1}-Theorem \ref{TL2}.

\beg{thm} Let $\hat L:=L_0+2\nn\log h_0$, and $L=\hat L+Z$ for some $C_b^1$-vector field $Z$ satisfying $\eqref{RM}$. The following assertions hold for
$q_\aa, \tt a(n)$ and, $\gg_{\aa,p,q} $ in \eqref{QU}, and a constant $\kk\ge 1$ with $\kk=1$ when $\pp M$ is convex.

 \beg{enumerate} \item[$(1)$] When $n=1$ and $\aa\in (\ff 1 2 \ss{10}-1, 1],$
$$\lim_{t\to\infty} \sup_{\nu\in\scr P} \E^\nu\big[\big|\{t\W_2(\mu_t^B,\mu)^2-\Xi^B(t)\}^+
+ \{t\kk\W_2(\mu_t^B,\mu)^2-\Xi^B(t)\}^-\big|^q\big]=0,\ \ q\in \Big[1,\ff{6}{5-2\aa}\Big).$$
\item[$(2)$] If $n<2(1+\aa)$, then for any $q\in [1,q_\aa)$ and $k\in (\ff {n+2} {2\aa},\infty]\cap [1,\infty]$,
 $$ \lim_{t\to\infty} \sup_{\nu\in \scr P_{k,R}} \E^\nu\Big[\Big|\{t\W_2(\mu_t^B,\mu)^2-\Xi^B(t)\}^++ \{t\kk\W_2(\mu_t^B,\mu)^2-\Xi^B(t)\}^-\Big|^q\Big]=0,\ \   R\in (0,\infty).$$
 \item[$(3)$] If $n<2(1+\aa)$, then
 $$\kk^{-1}\eta_Z^B\le \liminf_{t\to\infty} \inf_{\nu\in \scr P}t \E^\nu[\W_2(\mu_t^B,\mu)^2]
 \le  \limsup_{t\to\infty} \sup_{\nu\in \scr P}t \E^\nu[\W_2(\mu_t^B,\mu)^2]\le \eta_Z^B<\infty.$$
 \item[$(4)$] Let $n=2(1+\aa)$, i.e. $(n,\aa)\in \{(3,\ff 1 2), (4,1)\}.$ Then there exist constants $c, t_0>0$ such that
 $$\sup_{\nu\in \scr P} \E^\nu[\W_2(\mu_t^B,\mu)^2]\le c t^{-1}\log t,\ \ t\ge t_0.$$
 If $\pp M$ is convex  or $(n,\aa)=(3,\ff 1 2)$, then there exists a constant $c'>0$ such that
 $$\inf_{\nu\in \scr P} \E^\nu[\W_2(\mu_t^B,\mu)^2]\ge c' t^{-1}\log t,\ \ t\ge t_0.$$
  \item[$(5)$] If $n>2(1+\aa)$, then there exist constants $c, t_0>1$ such that
  $$c^{-1} t^{-\ff 2{n-2\aa}}\le \big(\E^\nu[\W_1(\mu_t^B,\mu)]\big)^2\le \E^\nu[\W_2(\mu_t^B,\mu)^2]\le
  c t^{-\ff 2{n-2\aa}},\ \ t\ge t_0, \nu\in \scr P.$$
 \item[$(6)$] If  $p,q\in [1,\infty)$ with $\gg_{\aa,p,q}<0,$ then  there exist  constants $c,t_0>1$ such that
 $$c^{-1}t^{-1}\le \big(\E^\nu[\W_1(\mu_t^B,\mu)]\big)^2\le \big(\E^\nu[\W_{2p}(\mu_t^B, \mu)^{2q}]\big)^{\ff 1 q}\le c t^{-1},\ \ t\ge t_0, \nu\in \scr P.$$
  \item[$(7)$] Let $p,q\in [1,\infty)$ with $\gg_{\aa,p,q}\ge 0.$ Then for any $\gg>\gg_{\aa,p,q},$ there exists a constant $c>0$ such that
      $$\sup_{\nu\in \scr P} \big(\E^\nu[\W_{2p}(\mu_t^B, \mu)^{2q}]\big)^{\ff 1 q}\le ct^{-\ff 1 {1+\gg}},\ \ t\ge 1.$$ If $\eqref{RS}$ holds, then there exists a constant $c>0$ such that
  $$ \sup_{\nu\in \scr P} \big(\E^\nu[\W_{2p}(\mu_t^B, \mu)^{2q}]\big)^{\ff 1 q}\le
        c t^{-\ff 1{1+\gg_{\aa,p,q}}} + c t^{-1}\log (1+t)1_{\{\gg_{\aa,p,q}=0\}}. $$
 \end{enumerate}\end{thm}

 \subsection{Subordinated Wright-Fisher diffusion process}

    Let   $a,b >\ff 1 4$  be two constants, and let
   $$\mu :=1_{[0,1]}(x) \ff{\GG(2a+2b)}{\GG(2a)\GG(2b)} x^{2a-1}(1-x)^{2b-1} \d x$$
   be the   Beta distribution  on  $M=[0,1]$. The Fisher-Wright diffusion process $\hat X_t$ is generated by
  $$\hat L:=  \ff 1 2  x(1-x)\ff{\d^2}{\d x^2} +\{a-(a+b)x\}\ff{\d}{\d x}. $$ Under the Riemannian metric
  $\<\pp_x,\pp_x\>= 2 \{x(1-x)\}^{-1},$  we have
   $$\GG(f, g)=\<\nn f,\nn g\>:= \ff 1 2 x (1-x)f'(x)g'(x),\ \ x\in M,$$
  $$\rr(x,y)= \ss 2 \int_x^y \{s(1-s)\}^{-\ff 1 2}\d s,\ \ 0< x\le y< 1.$$
Since ${\rm div}_\mu Z=0$ implies $Z=0$, we have $L=\hat L$.

\beg{lem}\label{WD} For the present model with $L=\hat L$, $(A_1)$ with $d=4(a\lor b)$ and $(A_2)$ hold,
      $\eqref{A2'}$ holds for $d'=2$, $(B)$ holds with $h(r)=\e^{Kr}$ for some constant $K>0$, and $\eqref{RS}$ holds.
\end{lem}

\beg{proof} Firstly, the condition  \eqref{A10} with $d= 4(a\lor b)$ is implied by \cite[Corollary 2.3]{FW}.
By \cite[(2.4)]{STA},   the Bakry-Emery curvature of $\hat L$ is bounded below by
$-K<0$ for some constant $K\ge 0$, so that
$$|\nn P_t f|\le \e^{Kt}P_t|\nn f|,$$  \eqref{A120} and $(B)$ hold for $\theta=m=1$,
$k(p)=\e^K$ and $h(r)=\e^{rK}$.

Next, for any $p\ge 1$ there exists a constant $c_1>0$ such that
\beg{align*}&\hat L\rr(X_0,\cdot)^{2p}(X_t)= 2p \rr(X_0,X_t)^{2(p-1)}L\rr(X_0,\cdot)(X_t)^2 +  p(p-1)\rr(X_0,X_t)^{2(p-1)}\\
& \le c_1 \rr(X_0,X_t)^{2(p-1)}, \end{align*} so that
\beq\label{PO}\E^\mu[ \rr(X_0,X_t)^{2p}]\le c_1\int_0^t \E^\mu[ \rr(X_0,X_s)^{2(p-1)}]\d s.\end{equation}
In particular, for $p=1$ we obtain \eqref{A13}, and for general $p\in \mathbb N$  it follows from
\eqref{PO} by the induction argument.

Moreover, we have $\ll_i=(a+b)i$ so that \eqref{A2'} holds for $d'=2.$ Indeed, according to the proof of \cite[Theorem 1.1]{FMW}, all eigenfunctions are polynomials. The trivial eigenvalue is $\ll_0=0$ with
$\phi_0=1$. For any $i\in \mathbb N$,   let
$$\phi_i(x):=\sum_{j=0}^i \aa_j x^j$$ with $\aa_i>0$ be the unitary eigenfunction for $\ll_i$. We have
$$-\ll_i \phi_i(x)=\hat L\phi_i(x).$$
Since the coefficients of $x^i$ in left hand and right hand sides  are
$-\ll_i\aa_i$ and $-i(a+b)\aa_i$ respectively, these two constants have to be equal each other, so that
$\ll_i=i(a+b).$

Finally, as explained in step (3) in the proof of Lemma \ref{XPP}, for \eqref{RS} it suffices to verify \eqref{VL} and \eqref{SP}. Since the curvature is bounded from below as indicated in the beginning of the proof, and since a one-dimensional ball is convex,
\eqref{SP} follows from \cite[Theorem 1.4]{W94}. So, it remains to \eqref{VL}. With the transform $x\to 1-x$, we only need to prove this condition for $x\in [0,\ff 1 2]$.
 Let $x\in [0,\ff 1 2]$ and $B(x,r):=\{y\in [0,1]: \rr(x,y)\le r\}$. Take, for instance,  $r_0=\ff 1 8 \rr(\ff 1 2,1)$ such that
 $$x_0:=\sup B( 1/2, 2r_0)\in (1/2, 1).$$ We have
 $$ c_0:=\inf_{x\in [0,\ff 1 2]} \mu(B(x,r_0))>0,$$ so that
 $$\mu(B(x,2r))\le 1 \le c_0^{-1}\mu(B(x,r)),\ \ \ r\ge r_0.$$
 Hence, we only need to consider $r\in (0,r_0).$ On the other hand,
    we find constants $c_2>c_1>0$ such that
 $$c_1^{-1} \big|\ss x -\ss y\big| \ge \rr(x,y)\ge c_2^{-1} \big|\ss x -\ss y\big|,\ \ x\in [0, 1/2], r\in (0,2r_0),$$
 so that for some constants $c_3,c_4>0$,
 $$\big[\big\{(\ss x-c_1r)^+\big\}^2, \big\{\ss x+ c_1 r\}^2\land x_0\big]\subset B(x,r)\subset \big[\big\{(\ss x-c_2r)^+\big\}^2, \big\{\ss x+ c_2 r\}^2\land x_0\big],\ \ r\in (0,2r_0).$$
 Noting that $0<1-x_0\le 1-s\le 1 $ for $s\in B(x,2r_0)$, we find  constants  $ c_4>c_3>0$ such that
 \beg{align*}&c_3 \big\{(x+c_3 r)^{2a} - x^{2a}\big\}\le \mu(B(x,r))\\
 &\le \mu(B(x,2r))\le c_4\big\{(x+c_4 r)^{2a} - [(x-c_4 r)^+]^{2a}\big\},\ \ x\in [0,1/2], r\in (0,r_0).\end{align*}
since  $\ff{(x+c_4r)^{2a}-\{(x-c_4r)^+\}^{2a}}{(x+c_3r)^{2a}-x^{2a}}$ is a continuous function of $(x,r)\in [0,\ff 1 2]\times [0,r_0],$ where when $r=0$ the function is understood as the limit $\ff {c_4}{c_3} $ as $r\to 0$,
we obtain
 $$\sup_{x\in [0,1/2], r\in (0,r_0)}\ff{\mu(B(x,2 r))}{\mu(B(x,r))}\le \sup_{x\in [0, \ff 1 2], r\in [0,r_0]} \ff{(x+c_4r)^{2a}-\{(x-c_4r)^+\}^{2a}}{(x+c_3r)^{2a}-x^{2a}}<\infty.$$
 Therefore, \eqref{VL} holds.
\end{proof}

In conclusion,
all assertions in Theorems \ref{TU1}-\ref{TL3}(1) hold  for $d=4(a\lor b)$ and $d'=2$ so that
\beq\label{QO1}\beg{split}& q_\aa:= \ff{4(a\lor b)}{4(a\lor b)-\aa},\ \ \ \aa(d,d')=\tt\aa:= (a\lor b)\big(\ss 5-1\big),\\
& \ \gg_{\aa,p,q}:= 2(a\lor b) (2- p^{-1}-q^{-1})-\aa.\end{split}\end{equation}
Noting that $\aa>\ff 4 3 (a\lor b)$ implies $q_\aa> \ff {4(a\lor b)}{2\aa}, 2(1+\aa)>d'=2$ and $\aa> \tt\aa,$
Theorems \ref{TU1}-\ref{TL3}(1) imply the following result.

\beg{thm} For the above  $L=\hat L$ and $\eta^B=\eta_Z^B$ for $Z=0$, the following assertions hold for
$q_\aa, \tt a$ and, $\gg_{\aa,p,q} $ in \eqref{QO1}.

 \beg{enumerate} \item[$(1)$] If $\aa>\ff 4 3(a\lor b),$   then
$$\lim_{t\to\infty} \sup_{\nu\in\scr P} \E^\nu\big[\big|t\W_2(\mu_t^B,\mu)^2-\Xi^B(t) \big|^q\big]=0,\ \ q\in [1,q_\aa).$$
\item[$(2)$] For any $q\in [1,q_\aa)$ and $k\in (\ff d {2\aa \i(q)},\infty]\cap [1,\infty]$, where we set $(\ff d {2\aa \i(q)},\infty]=\{\infty\}$ if $\aa=0$,
 $$ \lim_{t\to\infty} \sup_{\nu\in \scr P_{k,R}} \E^\nu\Big[\Big| t\W_2(\mu_t^B,\mu)^2-\Xi^B(t) \Big|^q\Big]=0,\ \   R\in (0,\infty).$$
 \item[$(3)$] $\eta^B<\infty$ and
 $ \limsup_{t\to\infty} \sup_{\nu\in \scr P}\big|t \E^\nu[\W_2(\mu_t^B,\mu)^2]-\eta^B\big|=0.$
 \item[$(4)$]   Let  $p,q\in [1,\infty)$ with $\gg_{\aa,p,q}<0.$ There exist  constants $c,t_0>1$ such that
 $$c^{-1}t^{-1}\le \big(\E^\nu[\W_1(\mu_t^B,\mu)]\big)^2\le \big(\E^\nu[\W_{2p}(\mu_t^B, \mu)^{2q}]\big)^{\ff 1 q}\le c t^{-1},\ \ t\ge t_0, \nu\in \scr P.$$
  \item[$(5)$] Let $p,q\in [1,\infty)$ with $\gg_{\aa,p,q}\ge 0.$ Then   there exists a constant $c>0$ such that
  $$ \sup_{\nu\in \scr P} \big(\E^\nu[\W_{2p}(\mu_t^B, \mu)^{2q}]\big)^{\ff 1 q}\le
        c t^{-\ff 1{1+\gg_{\aa,p,q}}} + c 1_{\{\gg_{\aa,p,q}=0\}} t^{-1}\log (1+t),\ \ t\ge 1. $$
 \end{enumerate}\end{thm}

\subsection{Subordinated subelliptic diffusions on $\SU(2)$}

Let $M=\SU(2)$ be the space of $2\times 2$, complex, unitary matrices with determinant $1,$ which is a $3$-dimensional compact Lie group, with Lie algebra ${\bf su}(2)$ and Riemannian metric $\<\cdot,\cdot\>$ given by
$${\bf su}(2):= {\rm span}\{U_1,U_2,U_3\},\ \ \ \<U_i,U_j\>=1_{\{i=j\}},\ \ 1\le i,j\le 3,$$
where for ${\bf i}=\ss{-1},$
$$U_1:= \left(\beg{matrix} 0 & 1\\ -1 & 0 \end{matrix}\right),\ \ U_2:= \left(\beg{matrix}  0 & {\bf i}\\ {\bf i} & 0 \end{matrix}\right),\ \ U_3:= \left(\beg{matrix}  {\bf i} & 0\\ 0 & -{\bf i} \end{matrix}\right).$$
For each $1\le i\le 3$, $U_i$ is understood as a left-invariant vector field defined as
$$U_if(x):=\lim_{\vv \downarrow 0} \ff{f(\e^{\vv U_i} x)-f(x)}{\vv},\ \ f\in C^1(\SU(2)).$$
 Then $[U_1,U_2]=2U_3,$ so that
$$\hat L:=U_1^2+U_2^2$$  satisfies H\"ormader's condition.  Moreover, $\hat L$
  is symmetric   in $L^2(\mu)$ where $\mu$ is the normalized Haar measure on $\SU(2)$, and the intrinsic distance $\rr$
induced by
$$\GG(f,g):= (U_1 f)(U_1g)+ (U_2f)(U_2g)$$
is the Carnot-Carath\'eodory distance. By Chow's theorem, $(M,\rr)$ is a compact geodesic   space.

To formulate the diffusion process $\hat X_t$ generated by $\hat L$, we use the cylindrical coordinates
introduced in \cite{[11]}:
$$\Big[0,\ff \pi 2\Big)\times [0,2\pi]\times [-\pi,\pi]\ni (r,\theta, z)\mapsto \e^{r(\cos \theta)U_1+ r(\cos\theta)U_2} \e^{z U_3}\in M:=\SU(2).$$
Under these coordinates, the diffusion process $\hat X_t:= (r_t,\theta_t,z_t)$ is constructed by solving the SDEs
\beq\label{SDE} \beg{split} &\d r_t= 2\cot(2r_t) \d t +\d B_t,\\
&\d (\theta_t, z_t)= \Big(\ff 2 {\sin \theta_t}, \tan r_t\Big)\d \tt B_t,\end{split}\end{equation}
where $(B_t,\tt B_t)$ is a two-dimensional Brownian motion, see \cite[Remark 2.2]{SU2}. The following lemma shows that conditions $(A_1), (A_2),\eqref{RS}$ and $\eqref{A2'} $ hold. However, due to the degeneracy of the diffusion, assumption $(B)$ may be invalid.

\beg{lem} Conditions $(A_1), (A_2), \eqref{RS}$ and $\eqref{A2'} $   hold  for $d =4$ and $d'=3.$

\end{lem}

\beg{proof} By \cite[Theorem 4.10]{SU2}, for any $p>1$ there exists a constant $c(p)>0$ such that
$$|\nn \hat P_tf|\le c(p)\e^{-2t} (P_t|\nn f|^p)^{\ff 1 p},\ \ \ t\ge 0.$$
So, \eqref{A120} holds.

According to \cite{FB2}, the generalized curvature-dimension condition $CD(\rr_1,\ff 1 2, 1,2)$ holds, so that
 \eqref{RS} is implied by \cite[Theorem 1.2]{FN}.

Let $\hat p_t$ be the heat kernel of $\hat P_t$ with respect to $\mu$.
By \cite[Proposition 3.1]{SU2} and the spectral representation of heat kernel, see also \cite{[8]},
all eigenvalues with multiplicities of $-\hat L$ are given by
$$\{\ll_i\}_{i\ge 0}=\big\{4k(k+|n|+1)+2|n|:\ \ n\in\mathbb Z, k\in \mathbb Z_+\big\}.$$
In particular, $\ll_1=2$. It is easy to see that for large $i\in\mathbb N$,
$$\#\big\{4k(k+|n|+1)+2|n|\le i:\ \ n\in\mathbb Z, k\in \mathbb Z_+\big\}$$ has order $i^{\ff 3 2}$, so that \eqref{A2'} holds for $d'=3$.

To verify \eqref{A10}, we use the cylindrical coordinates $(r,\theta, z)$, for which the identity matrix becomes ${\bf 0}:=(0,0,0).$  Let $\hat p_t$ be the heat kernel of $\hat P_t$ with respect to $\mu$. By \cite[Proposition 3.9]{SU2}, there exists a constant $c>0$ such that
$$\hat p_t({\bf 0},{\bf 0})\le c t^{-2},\ \ t\in (0,1].$$
By the left invariant of the heat kernel which follows from the same property of the generator $\hat L$, we
obtain
$$\|\hat P_t\|_{1\to\infty}=\sup_{x\in M} \hat p_t(x,x) \le c t^{-2},\ \ t\in (0,1].$$
This together with $\ll_1=2$ implies \eqref{A10} for $\ll=2.$

It remains to verify $(A_2)$. For any $(r,z)\in [0,\ff \pi 2)\times [-\pi,\pi]$, let $\theta(r,z)\in [-\pi,\pi]$ be the unique solution to the equation
$$\theta(r,z)-z= \ff{(\cos r)(\sin\theta(r,z))\arccos[(\cos\theta(r,z))\cos r]}{\ss{1-(\cos^2r)\cos^2\theta(r,z)}}.$$
By \cite[Remark 3.12]{SU2}, the distance of $x:=(r,\theta,z)$ to ${\bf 0}$ depends only on $(r,z)$, and there exists a constant $c_1>0$ such that
\beq\label{LOP}\rr({\bf 0}, x)^2=\ff{(\theta(r,z)-z)^2\tan^2 r}{\sin^2 \theta(r,z)}\le c_1\sin^2 r\le c_1 r^2.\end{equation}
On the other hand, letting $\hat X_0={\bf 0}$, by \eqref{SDE} and It\^o's formula, for any $p\ge 1$ we find a constant $c_1(p)>0$ such that
$$\d r_t^{2p}\le c(p) r_t^{2(p-1)}\d t +\d M_t$$ for some martingale $M_t$. So, we find a constant $c_2(p)>0$ such that
$$\E[r_t^{2p}]\le k(p) t^p,\ \ t\ge 0.$$
 Combining this with \eqref{LOP} and using the left invariance of the heat kernel, we obtain
 $(A_2)$.

\end{proof}

By the above lemma and that $M={\bf SU}(2)$  is a Polish space,  we conclude that
   all assertions in Theorems \ref{TU1}-Theorem \ref{TU3} and Theorem \ref{TL3}(1) hold    for $d=4, d'=3$ and
 \beq\label{QUN} \beg{split} &q_\aa:= \ff 8 {9-2\aa},\ \ \   \
  \gg_{\aa,p,q}:=\ff 1 2 -\aa+ 2(2-p^{-1}-q^{-1}).\end{split}\end{equation}
Noting that $1<q_\aa<2$ for $\aa\in (\ff 1 2,1]$, so that ${\rm i}(q)=1$ for $q\in [1,q_\aa)$,  we have the following result.

\beg{thm} Let $\hat L:=U_1^2+U_2^2$, and $L=\hat L+Z$ for some $C_b^1$-vector field $Z$ satisfying $\eqref{RM}$. The following assertions hold for
$q_\aa$ and $\gg_{\aa,p,q} $ in \eqref{QUN}.

 \beg{enumerate} \item[$(1)$]   If $\aa\in(\ff 1 2,1]$, then for any $q\in [1,q_\aa)$ and $k\in (\ff 2\aa, \infty]\cap [1,\infty]$,
 $$ \lim_{t\to\infty} \sup_{\nu\in \scr P_{k,R}} \E^\nu\Big[\Big|\{t\W_2(\mu_t^B,\mu)^2-\Xi^B(t)\}^+\Big|^q\Big]=0,\ \   R\in (0,\infty).$$
 \item[$(2)$] If $\aa\in (\ff 1 2,1]$, then
 $$  \limsup_{t\to\infty} \sup_{\nu\in \scr P}t \E^\nu[\W_2(\mu_t^B,\mu)^2]\le \eta_Z^B<\infty.$$
 \item[$(3)$] If  $\aa\in (0,\ff 1 2]$, then there exist constants $c, t_0>0$ such that
 $$\sup_{\nu\in \scr P} \E^\nu[\W_2(\mu_t^B,\mu)^2]\le ct^{-\ff{2}{3-2\aa}}+c 1_{\{\aa=\ff 1 2\}} t^{-1}\log t,\ \ t\ge t_0.$$
  \item[$(4)$]  If  $p,q\in [1,\infty)$ with $\gg_{\aa,p,q}<0,$ then  there exist  constants $c,t_0>1$ such that
 $$c^{-1}t^{-1}\le \big(\E^\nu[\W_1(\mu_t^B,\mu)]\big)^2\le \big(\E^\nu[\W_{2p}(\mu_t^B, \mu)^{2q}]\big)^{\ff 1 q}\le c t^{-1},\ \ t\ge t_0, \nu\in \scr P.$$
  \item[$(5)$] Let $p,q\in [1,\infty)$ with $\gg_{\aa,p,q}\ge 0.$ Then for any $\gg>\gg_{\aa,p,q}$, there exists a constant $c>0$ such that for any $t\ge 1$,
  $$ \sup_{\nu\in \scr P} \big(\E^\nu[\W_{2p}(\mu_t^B, \mu)^{2q}]\big)^{\ff 1 q}\le
        c t^{-\ff 1{1+\gg_{\aa,p,q}}} + c1_{\{\gg_{\aa,p,q}=0\}} t^{-1}\log (1+t). $$
 \end{enumerate}\end{thm}

 \paragraph{Acknowledgement.} The author  would like to thank Dr. Jie-Xiang Zhu and Bingyao Wu for corrections. 

\end{document}